\documentclass[11pt]{article}

\usepackage{graphicx}
\usepackage{amsmath,amssymb,amsfonts}
\usepackage{subfigure}
\usepackage[T1]{fontenc}
\usepackage{color}
\usepackage{url}
\usepackage{algorithm,algpseudocode}
\usepackage{enumitem}
\usepackage{mathabx} % for \widecheck symbol

\usepackage[hidelinks]{hyperref}  % not compatible with Springer macros
\usepackage{authblk}
\usepackage{multirow}
\usepackage{color}

    \usepackage{amsthm}

	\usepackage{todonotes} 			           %
\newcommand{\KK}{\boldsymbol{K}_{T}}
\newcommand{\KKO}{\boldsymbol{K}_{T0}}
\newcommand{\KKhat}{\widehat{\boldsymbol{K}_{T}}}
\newcommand{\II}{\boldsymbol{I}}
\newcommand{\BB}{\boldsymbol{B}}
\newcommand{\BBhat}{\widehat{\boldsymbol{B}}}
\newcommand{\DeltaK}{\Delta\boldsymbol{K}}
\newcommand{\DeltaKhat}{\widehat{\Delta\boldsymbol{K}}}
\newcommand{\GG}{\boldsymbol{G}}
\newcommand{\DD}{\boldsymbol{D}}

\newcommand{\svec}{\boldsymbol{s}}
\newcommand{\svecO}{\boldsymbol{s}^{(0)}}
\newcommand{\sk}{\boldsymbol{s}^{(k)}}
\newcommand{\skk}{\boldsymbol{s}^{(k+1)}}

\newcommand{\stilde}{\widetilde{\boldsymbol{s}}}
\newcommand{\shat}{\widehat{\boldsymbol{s}}}
\newcommand{\sstar}{\boldsymbol{s}^{*}}
\newcommand{\RR}{\boldsymbol{r}}
\newcommand{\ff}{\boldsymbol{f}}
\newcommand{\lbold}{\boldsymbol{l}}
\newcommand{\ds}{\displaystyle}
\newcommand{\zzero}{\boldsymbol{0}}
\newcommand{\uu}{\boldsymbol{u}}

\newcommand{\uhat}{\widehat{\boldsymbol{u}}}
\newcommand{\bssig}{\boldsymbol{\sigma}}
\newcommand{\SB}{\boldsymbol{S}_{B}}
\newcommand{\yy}{\boldsymbol{y}}
\newcommand{\bsrho}{\boldsymbol{\rho}}
\newcommand{\bslambda}{\boldsymbol{\lambda}}
\newcommand{\bslambdak}{\boldsymbol{\lambda}^{(k)}}
\newcommand{\bslambdakk}{\boldsymbol{\lambda}^{(k+1)}}
\newcommand{\bslambdaO}{\boldsymbol{\lambda}^{(0)}}
\newcommand{\bslambdatilde}{\widetilde{\boldsymbol{\lambda}}}

\begin{document}

\title{Inexact Newton combined approximations in the topology optimization of 
geometrically nonlinear elastic structures and compliant mechanisms}

    \author[1]{T. A. Senne}
        \author[2]{F. A. M. Gomes}
    \author[2]{S. A. Santos} 
    \affil[1]{ Institute of Science and Technology, Federal
           University of S\~ao Paulo, Avenida Cesare Mansueto Giulio Lattes, 1201, 12247-014, S\~ao Jos\'e dos Campos, Brazil. Email: {\tt senne@unifesp.br}}
    \affil[2]{Institute of Mathematics, University of Campinas, Rua Sergio
    Buarque de Holanda, 651, 13083-859, Campinas, S\~ao Paulo, Brazil
    \authorcr Emails: {\tt sasantos@unicamp.br, chico2@unicamp.br}}
    \renewcommand\Affilfont{\itshape\small}

    \date{\today}

\maketitle

\begin{abstract}
This work blends the inexact Newton method with iterative combined approximations (ICA) for solving topology optimization problems under the assumption of geometric nonlinearity. The density-based problem formulation is solved using a sequential piecewise linear programming (SPLP) algorithm. Five distinct strategies have been proposed to control the frequency of the factorizations of the Jacobian matrices of the nonlinear equilibrium equations. Aiming at speeding up the overall iterative scheme while keeping the accuracy of the approximate solutions, three of the strategies also use an ICA scheme for the adjoint linear system associated with the sensitivity analysis. The robustness of the proposed reanalysis strategies is corroborated by means of numerical experiments with four benchmark problems -- two structures and two compliant mechanisms. Besides assessing the performance of the strategies considering a fixed budget of iterations, the impact of a theoretically supported stopping criterion for the SPLP algorithm was analyzed as well.


\end{abstract}

%===========================================
\section{Introduction}\label{sec:intro}

One of the most common problems in topology optimization is the minimization of the mean compliance of  structures, under the constraints of static equilibrium and a prescribed volume of material~\cite{bsbook}. In several circumstances, the requirements upon the structure may impose an elastic but nonlinear relation between strains and displacements, meaning that the structure is under geometric nonlinearity. In such a case, the computation of the objective function of the topology optimization problem demands the solution of a system of nonlinear equations that models the nodal displacements of the structure due to the external forces. When it comes to compliant mechanisms, under the geometric nonlinearity assumption, addressing nonlinear systems of equations is mandatory as well.  The approximate solutions of these systems are usually computed by Newton's method, constituting, by far, the most expensive step in the process of achieving optimal configurations.  	

Several attempts have been made to reduce the effort spent to compute the nodal displacements of the structures and mechanisms~\cite{Kirsch2010}. Challenges such as parallel and numerical scalability have been pursued, specially when it comes to 3-D problems (see e.g.~\cite{Amir2015} and references therein). The use of preconditioned conjugate gradients has been exploited in the solution of the linear systems associated with the equilibrium equations~\cite{AmirSigmundLazarovSchevenelsl2012,AmirStolpeSigmund2010,LiuWuLi2014,LongEtal2019}. Meshless-based local reanalysis~\cite{Cheng2017} and a meshfree reduced basis approach for large deformations~\cite{Nguyen2021} have been analyzed as well. In the last years, educational codes handling geometric nonlinearity, based on the SIMP (Solid Isotropic Material with Penalization) method~\cite{Chen2019,Zhu2021}, and on the BESO (Bi-directional Evolutionary Structural Optimization) method \cite{Han2021}, also became available. The recent survey~\cite{Mukherjee2021} discusses additional aspects concerning high-performance computing, approximate reanalysis, reduced-order modeling, and multigrid methods.

Inspired on the \emph{combined approximations} approach \cite{AmirKirschSheinman2008,Bogomolny2010,kirschbook}, we have analyzed some strategies that merge the inexact Newton scheme with iterative combined approximations for minimizing the mean compliance of structures and for maximizing the displacement of compliant mechanisms. Our approach rests upon a sequential piecewise linear programming (SPLP) algorithm \cite{GS2014} that, aside from having a stopping criterion with theoretical support, has proven to be efficient and robust. In a previous work~\cite{SGS2019}, the approximate reanalysis technique has been applied in combination with the SPLP algorithm to solve three benchmark structure problems under small displacements. In the current work, besides addressing two structures, we have also applied the devised strategies to the design of two mechanisms. A detailed analysis of the results has been reported, assessing the scope and the impact of the methodology.

This text is organized as follows. The problem formulation is presented in~Section~\ref{sec:probform}. The elements of the reanalysis technique that encompasses the inexact Newton method with iterative combined approximations are detailed in Section~\ref{sec:inca}. The experimental setup is provided in Section~\ref{sec:numtests}, namely the test problems, the proposed strategies and the implementation details.
The analysis of the performance of the strategies with a fixed budget is done in Section~\ref{sec:perfstr}, whereas Section~\ref{sec:meshref} presents the results concerning the effect of the mesh refinement. 
The convergence performance of the SPLP with the nonlinear reanalysis is reported in Section~\ref{sec:convsplp}. The behavior of $\| \BB \|_2$ is assessed in Section~\ref{sec:normB}, and Section~\ref{sec:final} contains our final remarks.

%===========================================
\section{Problem formulation}\label{sec:probform}

Suppose that $\Omega$ is the design domain of a structure. In the context of topology optimization, we want the distribution of material on the domain $\Omega$ to be optimal in some sense, satisfying certain constraints. In this work, we consider that $\Omega$ is two-dimensional and, with the aim of making the problem computationally affordable, we discretize it into rectangular finite elements, using bilinear interpolating functions to approximately evaluate the displacements. 

Let $\Omega_{d}$ be the discretized domain, partitioned in disjointed subdomains $\Omega_i$, i.e. $\Omega_{d} = \cup_{i=1}^{n_{el}} \Omega_i$. We denote by $n_{el}$ the number of elements of $\Omega_{d}$ and by $n_{nd}$ the number of its nodes. We associate with each $\Omega_i$ a continuous variable $\rho_{i}$ representing the density of material, in such a way that $\rho_{i} = 1$ if the element is solid or $\rho_{i} = 0$ if the element is void. 
 
We address two topology optimization problems: the compliance minimization of rigid structures and the displacement maximization of compliant mechanisms. In both cases, the optimization problem can be stated as 
\begin{equation}
    \begin{array}{ll}
        \ds{\min_{\bsrho}} &  {\lbold}^{T}\uu  \\ 
							 \mbox{s.t.} &  \RR(\uu,\,\bsrho) = \, \zzero \\
					                 &  \sum^{n_{el}}_{i=1} v_{i}\rho_{i} = V^{*} \\ 
                           &  0 < \rho_{\min} \leq \rho_{i} \leq 1, \qquad i = 1,\dots,n_{el},
	  \end{array} 
	  \label{prob01}
\end{equation}
where $\bsrho \in \mathbb{R}^{n_{el}}$ is the vector of densities of each element of $\Omega_{d}$ (that are our design variables), $\uu \in \mathbb{R}^{2n_{nd}}$ is the vector of nodal displacements, $v_{i}$ is the volume of the $i$-th element, $V^{*}$ is the prescribed volume of the structure, and $\rho_{min}$ is a sufficiently small positive number used to avoid numerical instabilities. For the compliance minimization of rigid structures, $\lbold \in \mathbb{R}^{2n_{nd}}$ is the vector of nodal external forces, $\ff$, and in the case of the maximization of the displacements for compliant mechanisms it is a vector filled with zeros, except at the positions corresponding to the output points (the points in which the displacement is to be maximized), where it is equal to $-1$. We remark that, since we adopt the SIMP interpolation model \cite{Bendsoe89}, $\uu$~depends on the vector of densities $\bsrho$, so we can express $\uu \equiv \uu(\bsrho)$.

The system
\begin{equation}
    \RR(\uu,\,\bsrho) = \zzero\,
    \label{eq01b}
\end{equation}
represents the static equilibrium conditions of the structure (or compliant mechanism).
Under the assumption of geometric nonlinearity, we have 
\begin{equation}
    \RR(\uu,\,\bsrho) = \ff_{int}(\uu,\,\bsrho) - \ff,
    \label{eq01a}
\end{equation}
where $\ff_{int}(\uu,\,\bsrho) \in \mathbb{R}^{2n_{nd}}$ is the vector of internal nodal forces. We say that $\RR(\uu,\,\bsrho)$ in (\ref{eq01a}) is the \textit{residual} associated with the nonlinear system (\ref{eq01b}).

Usually, the Newton-Raphson method (that, from now on, will be referred as \textit{Newton's method}) is used to find an approximate solution for the nonlinear system (\ref{eq01b}). Given an initial point $\uu^{(0)}$ and keeping $\bsrho$ fixed, at each iteration of  Newton's method, we solve the linear system
\begin{equation}
    \KK(\uu^{(\ell)},\,\bsrho)\svec \, = \, -\RR(\uu^{(\ell)},\,\bsrho), 
    \label{eq02}
\end{equation}
where $\KK \equiv \KK(\uu,\,\bsrho) \in \mathbb{R}^{2n_{nd} \times 2n_{nd}}$ is the Jacobian matrix of the nonlinear system (\ref{eq01b}), also known as the \textit{global tangent stiffness matrix}, that is supposed to be nonsingular. The solution of the linear system~(\ref{eq02}), denoted by $\boldsymbol{s}^{(\ell)}$, is used to obtain the new approximate solution, $\boldsymbol{u}^{(\ell+1)} = \boldsymbol{u}^{(\ell)} + \boldsymbol{s}^{(\ell)}$. Usually, we consider that $\boldsymbol{u}^{(\ell)}$ is a good approximate solution to the original nonlinear system (\ref{eq01b}) whenever $\| \RR(\boldsymbol{u}^{(\ell)},\,\bsrho) \| < \varepsilon$, where $\varepsilon$ is a prescribed small positive number. 

It is important to mention that the definitions of the global tangent stiffness matrix $\KK$ and of the residual $\RR$ depend on the hyperelastic model for the material of the structure. In this work, we adopt the neo-Hookean material model of Simo-Ciarlet \cite{Lahuerta2013}, combined with the SIMP model, so we have
$$ 
  \KK(\uu, \bsrho)  =  \sum_{i=1}^{n_{el}} \rho^{p}_{i} \int_{\Omega_{i}}\!\! \GG^{T}\DD_i\GG\,\,d\Omega_{i}
  \mbox{ \ and \ }
  \ff_{int}(\uu,\bsrho) = \sum^{n_{el}}_{i=1} \rho^{p}_{i}\int_{\Omega_{i}}\!\! \GG^{T}\bssig_i\,\,d\Omega_{i}, 
$$
in which $p$ is the penalty parameter of the SIMP model ($p > 1$), $\GG \in \mathbb{R}^{4 \times 8}$ is the matrix of the derivatives of the shape function with respect to the displacements and, for the domain $\Omega_i$ of the $i$-th element of the discretized domain $\Omega_d$, $\DD_i \in \mathbb{R}^{4 \times 4}$ is the tangent stiffness modulus matrix and $\bssig_i \in \mathbb{R}^{4}$ is the first Piola-Kirchhoff stress tensor.

%===========================================
\section{Inexact Newton with Iterative Combined Approximations}\label{sec:inca}

The factorization of $\KK$ is the most expensive step in the process of solving the linear system (\ref{eq02}).  Although the Cholesky factorization is usually employed for symmetric matrices, the  $LDL^T$ factorization \cite{Golub} is required here, since $K_T$ is not necessarily positive definite. With the aim of saving computational effort, several reanalysis techniques were proposed in the literature. One of the most common reanalysis strategy, named \textit{Combined Approximations} (CA), will be summarized here to ground our presentation, following the lines presented by Amir, Kirsch and Sheinman \cite{AmirKirschSheinman2008} and by Senne, Gomes and Santos \cite{SGS2019}. 

First of all, let $\KKO$ be the Jacobian at the $i$-th iteration of Newton's method. Supposing that its factorization is available, it will be reused at the subsequent $m$ iterations $i+1,\,i+2,\,\dots,\,i+m$, and updated at the iteration $i+m+1$. Denoting by $\KK$ the Jacobian at the $(i+j)$-th iteration (with $1 \leq j \leq m$), we can rewrite the linear system (\ref{eq02}) as 
\begin{equation}
   (\KKO + \DeltaK)\svec \, = \, -\RR,
   \label{eq02aa}
\end{equation}
where $\DeltaK := \KK - \KKO$. After rearranging the terms in (\ref{eq02aa}), we obtain
\begin{equation}
    \svec \, = \, \stilde - \BB\svec,
    \label{eq02b}
\end{equation}
where 
\begin{equation}
\stilde := -\KKO^{-1}\RR
\qquad \mbox{and } \qquad
\BB := \KKO^{-1}\DeltaK. 
\label{eq:stilde&matB}
\end{equation}
So, using (\ref{eq02b}), we state the following sequence of steps:
\begin{equation}
    \skk\, = \, \stilde - \BB\sk, \qquad k = 0,\,1,\,2,\,\dots.
    \label{eq03}
\end{equation}

Choosing $\svec^{(0)} = \stilde$, the CA strategy consists in generating the first $q$ terms of the sequence (\ref{eq03}) to build an approximate solution $\shat$, given by
\begin{equation}
    \shat \, = \, y_{1}\svec^{(1)} + y_{2}\svec^{(2)} + \cdots + y_{q}\svec^{(q)} \, = \, \SB\yy,
		\label{eq03aa}
\end{equation}
where 
\begin{equation*}
    \SB \, = \, [\svec^{(1)} \,\,\, \svec^{(2)} \,\,\, \cdots \,\,\, \svec^{(q)}] \qquad \mbox{and} \qquad \yy = [y_{1} \,\,\, y_{2} \,\,\, \cdots \,\,\, y_{q}]^{T}. 
\end{equation*}
We remark that $q$ should be a small integer. Usually, $q \leq 10$.

The vector of coefficients $\yy$ is found by solving the small ($q \times q$) linear system 
\begin{equation*}
    \SB^{T}\KK\SB\yy \, = \, -\SB^{T}\RR,
\end{equation*}
obtained by substituting $\shat$ in the original linear system (\ref{eq02}) and then  multiplying both sides of the resulting equation on the left by $\SB^{T}$. Generally, the relative norm of the residual of the linear system (\ref{eq02}) is the measure adopted to decide if the approximate solution $\shat$ will be accepted.  

In this work, we propose a new iterative method to get an approximate solution to the linear system~(\ref{eq02}), named \textit{Iterative Combined Approximations} (ICA), based on the sequence $\left\{\sk\right\}^{+\infty}_{k=0}$ defined in (\ref{eq03}), and described below.  It is worth mentioning that a distinct iterative CA approach to address topological modifications for the static reanalysis of structures has been presented in \cite{ChenYang2004}.

Given  $\stilde$ as in \eqref{eq:stilde&matB}, and any initial point $\svecO$, the general term of the sequence~(\ref{eq03}), for any $k$, can be written as
\begin{equation}
    \skk \, = \, [\II - \BB + \BB^{2} - \BB^{3} + \cdots + (-1)^{k}\BB^{k}]\,\stilde + (-1)^{k+1}\BB^{k+1}\svecO.
    \label{eq03a}
\end{equation}

By considering any vector $p$-norm and the consistent norm of operators, if $\| \BB \| < 1$, the matrix  $\II + \BB$ is nonsingular (cf. \cite[\S 2.3.4]{Golub}), and its inverse can be expressed by means of the Newmann's series  \cite{JuanjuanHu2018,ZuoBaiYu2016,ZuoYuZhaoZhang2012}
\begin{equation*}
    (\II + \BB)^{-1} \, = \, \textstyle\sum^{+\infty}_{j=0} (-1)^{j}\BB^{j}.
\end{equation*}
Notice that, from \eqref{eq:stilde&matB}, since $  \BB=  \KKO^{-1} \KK  -  \II$, as long as $\KK $ does not differ too much from $\KKO$, the assumption $\| \BB \| < 1$ is reasonable.

Defining $\sstar := (\II + \BB)^{-1}\stilde$, from (\ref{eq03a}) we have
\begin{equation}
    \skk - \sstar \, = \, -\left[\textstyle \sum^{+\infty}_{j=k+1} (-1)^{j}\BB^{j}\right]\stilde + (-1)^{k+1}\BB^{k+1}\svecO.
    \label{eq04}
\end{equation}
Since
\begin{eqnarray}
        \textstyle\sum^{+\infty}_{j=k+1} (-1)^{j}\BB^{j} 
        & = & (-1)^{k+1}\BB^{k+1} + (-1)^{k+2}\BB^{k+2} + (-1)^{k+3}\BB^{k+3} + \cdots \nonumber \\ 
        & = & (-1)^{k+1}\BB^{k+1}(\II - \BB + \BB^{2} - \BB^{3} + \cdots) \nonumber \\ [3pt]
        & = & (-1)^{k+1}\BB^{k+1}(\II + \BB)^{-1}, \label{eq05}
\end{eqnarray} 
substituting (\ref{eq05}) into (\ref{eq04}), we obtain
\begin{eqnarray}
    \skk - \sstar &= & (-1)^{k+1}\BB^{k+1}[\svecO - (\II + \BB)^{-1}\stilde] \nonumber \\ [3pt]
    & = & (-1)^{k+1}\BB^{k+1}[\svecO - \sstar] 
    \label{eq06}
\end{eqnarray}
and hence
\begin{equation*}
   0 \, \leq \, \| \skk - \sstar \| \, \leq \, \| \BB \|^{k+1}\,\| \svecO -\sstar \|. 
\end{equation*}

Once we have $\lim_{k \rightarrow +\infty} \| \BB \|^{k+1} = 0$ (provided that $\| \BB \| < 1$), and assuming that $\| \svecO -\sstar \|$ remains bounded, by the Squeeze Theorem we conclude that 
\begin{equation*}
    \lim_{k \rightarrow +\infty} \| \skk - \sstar \| \, = \, 0,
\end{equation*}
and, consequently,
    $\lim_{k \rightarrow +\infty} \skk \, = \, \sstar$.

Now, we will show that $\sstar$ is the (unique) solution of the linear system $\KK\svec \, = \, -\RR$ if, and only if, $\sstar$ is the (unique) solution of $(\II + \BB)\svec = \stilde$. In fact, 
\begin{equation*} 
    \begin{array}{lll}
        \KK\sstar \, = \, -\RR & \Longleftrightarrow & (\KKO + \DeltaK)\sstar \, = \, -\RR \\ [3pt] 
                               & \Longleftrightarrow & \sstar \, = \, -\KKO^{-1}\RR - \KKO^{-1}(\DeltaK)\sstar \\  [3pt]
                               & \Longleftrightarrow & \sstar \, = \, \stilde - \BB\sstar \\ [3pt] 
                               & \Longleftrightarrow & (\II + \BB)\sstar = \stilde.
    \end{array}
\end{equation*}

Using (\ref{eq06}), we observe that
\begin{equation*}
    \begin{array}{lll}
        \skk - \sstar & = & (-1)^{k+1}\BB^{k+1}[\svecO - (\II + \BB)^{-1}\stilde] \\ [3pt]
                      & = & (-1)\BB(-1)^{k}\BB^{k}[\svecO - (\II + \BB)^{-1}\stilde] 
                      = -\BB(\sk - \sstar)
    \end{array}
\end{equation*}
and, consequently,
    \begin{equation}
        \| \skk - \sstar \| \, \leq \, \| \BB \|\,\| \sk - \sstar \|.
        \label{eq07}
    \end{equation} 
Thus, by (\ref{eq07}) and under the assumption that $\| \BB \| < 1$, we conclude that the sequence $\left\{\sk\right\}^{+\infty}_{k=0}$ generated by the ICA method has a linear rate of convergence.

Let $\uhat$ be an approximate solution of the nonlinear system (\ref{eq01b}) found by means of Newton's method, such that $\| \RR(\uhat,\,\bsrho)\| < \varepsilon$. It can be shown that the partial derivatives of the objective function $F(\bsrho) = \lbold^{T}\uu(\bsrho)$ are given by $\partial F/\partial \rho_{e} \, = \, \bslambda_{e}^{T}\,\partial \RR/\partial \rho_{e}$, where $\bslambda$ is the solution of the adjoint linear system 
    \begin{equation}
        \KKhat\bslambda = -\lbold,
        \label{eq09}
    \end{equation}
$\KKhat$ is the global tangent stiffness matrix evaluated at $\uhat$, and $\bslambda_{e}$  encompasses the components of the vector $\bslambda$ associated with the $e$-th element. 
		
		In order to avoid the factorization of $\KKhat$ (and, thus, to accelerate the solution of our topology optimization problems), the ICA method can also be applied to find an approximate solution of the linear system~(\ref{eq09}). As done in~(\ref{eq03}), we define the sequence
		$$
		    \bslambdakk = \bslambdatilde - \BBhat\bslambdak, \qquad k = 0,\,1,\,2,\,\dots
		$$
		where $\bslambdatilde = -\KKO^{-1}\lbold$, $\BBhat = \KKO^{-1}\DeltaKhat$, $\DeltaKhat = \KKhat-\KKO$, and we take $\bslambdaO = \bslambdatilde$. Following the same steps shown above, we can prove that, if $\| \BBhat \| < 1$, the sequence $\{\bslambdak\}^{+\infty}_{k=0}$ converges to the unique solution $\bslambda^{*} = -\KKhat^{-1}\lbold$ of the linear system (\ref{eq09}) with a linear rate of convergence. 
		
%===========================================
\section{Test problems, proposed strategies and implementation details}\label{sec:numtests}

In this section we present the computational settings and choices considered in this work. 

\subsection{Test problems} \label{sub:testproblems}

In our numerical investigation, we have addressed the two structures and the two compliant mechanisms described below.\\

\begin{figure}[htbp!]
   \centering 
      \begin{tabular}{cc} 
          \includegraphics[width=0.47\textwidth]{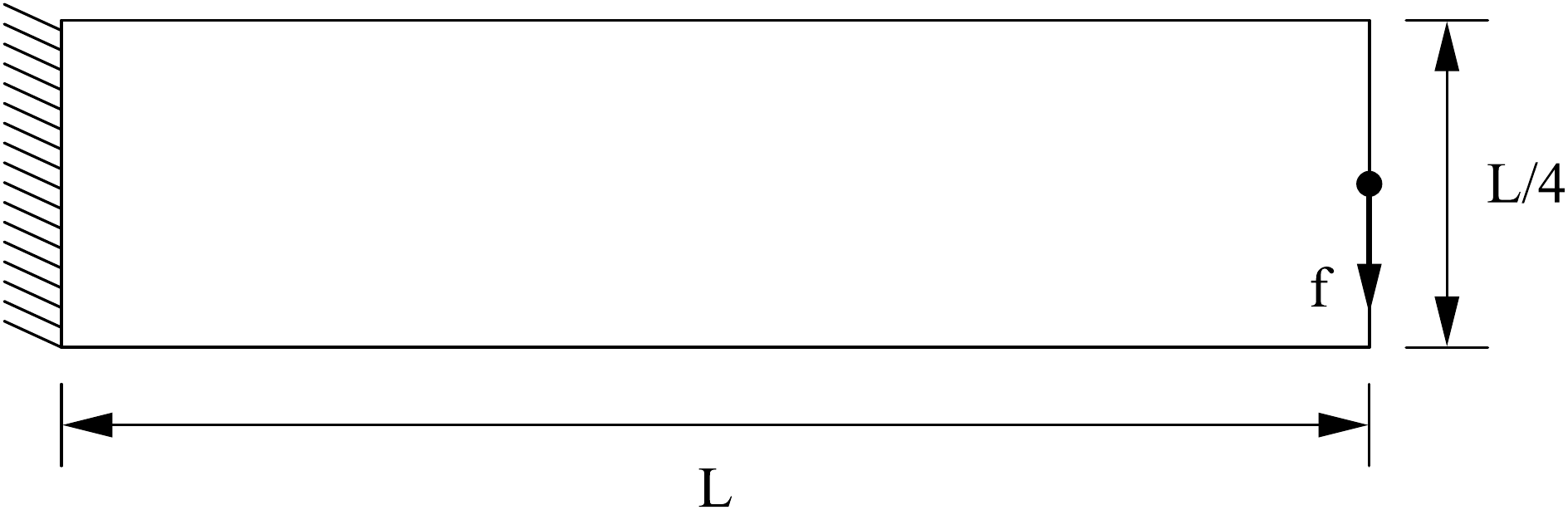}
        & \includegraphics[width=0.47\textwidth]{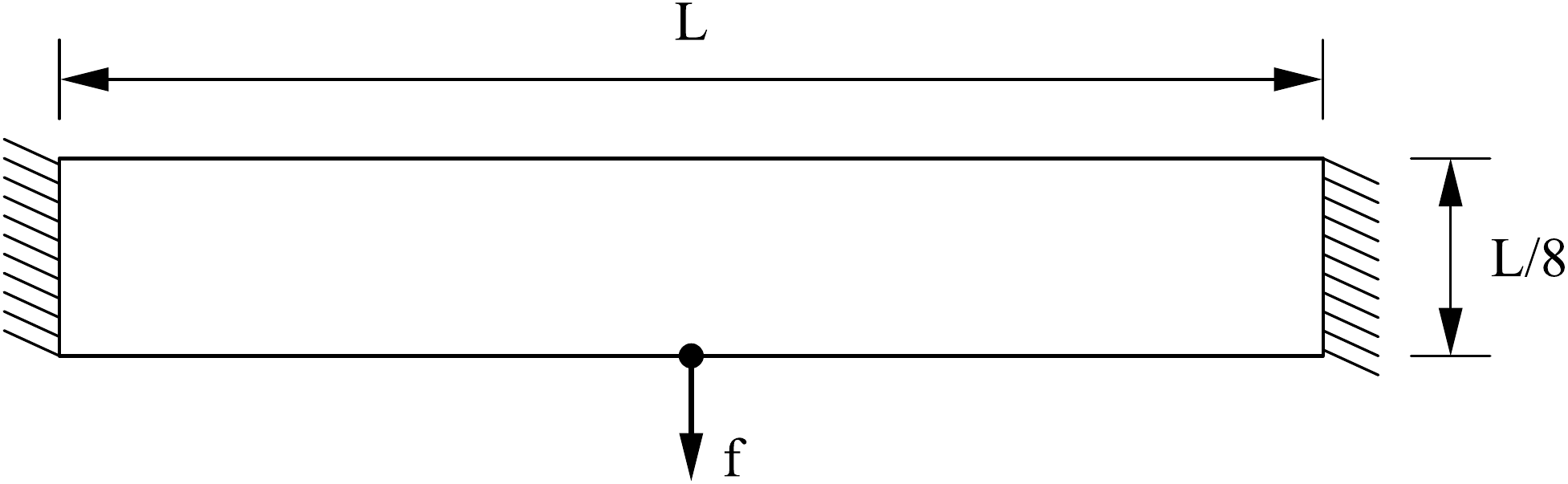} \\ 
      \end{tabular}
      \caption{Domains of the structures: cantilever beam (left) and slender beam (right).} 
      \label{fig:domainStruct}
\end{figure}

\noindent \textbf{Cantilever beam.} Fig.~\ref{fig:domainStruct}~(left) shows the design domain for this problem. Following the suggestions of \cite{Chen2019}, we define $L = 120\,mm$\, and suppose that the material used to construct the structure has a Young's modulus of $3000\,N/mm^{2}$ and a Poisson's ratio of $0.4$. The thickness of the structure is equal to $1\,mm$. An external load of $120\,N$ is applied downwards at the midpoint of the right side. The structure is supported at the left edge, where all of the nodes are fixed. The domain is discretized into $40000  \ (400 \times 100)$ finite elements, and the optimal structure must contain $50\,\%$ of the domain's volume. \\

\noindent \textbf{Slender beam.} The design domain of this structure is presented in Fig.~\ref{fig:domainStruct} (right), in which $L = 400\,mm$. An external load of $40\,N$ is applied downwards at the center of the domain's basis, and the structure is supported at the left and at the right edges, where all of the nodes are fixed. The Poisson's ratio and the Young's modulus of the material are, respectively, $0.3$ and $3000\,N/mm^{2}$. The structure has a thickness of $1\,mm$. The domain is discretized into $45000 \ (600 \times 75) $ finite elements, and the optimal structure must contain $20\,\%$ of the domain's volume. \\

\begin{figure}[htbp!]
   \centering 
      \begin{tabular}{cc} 
          \includegraphics[width=0.47\textwidth]{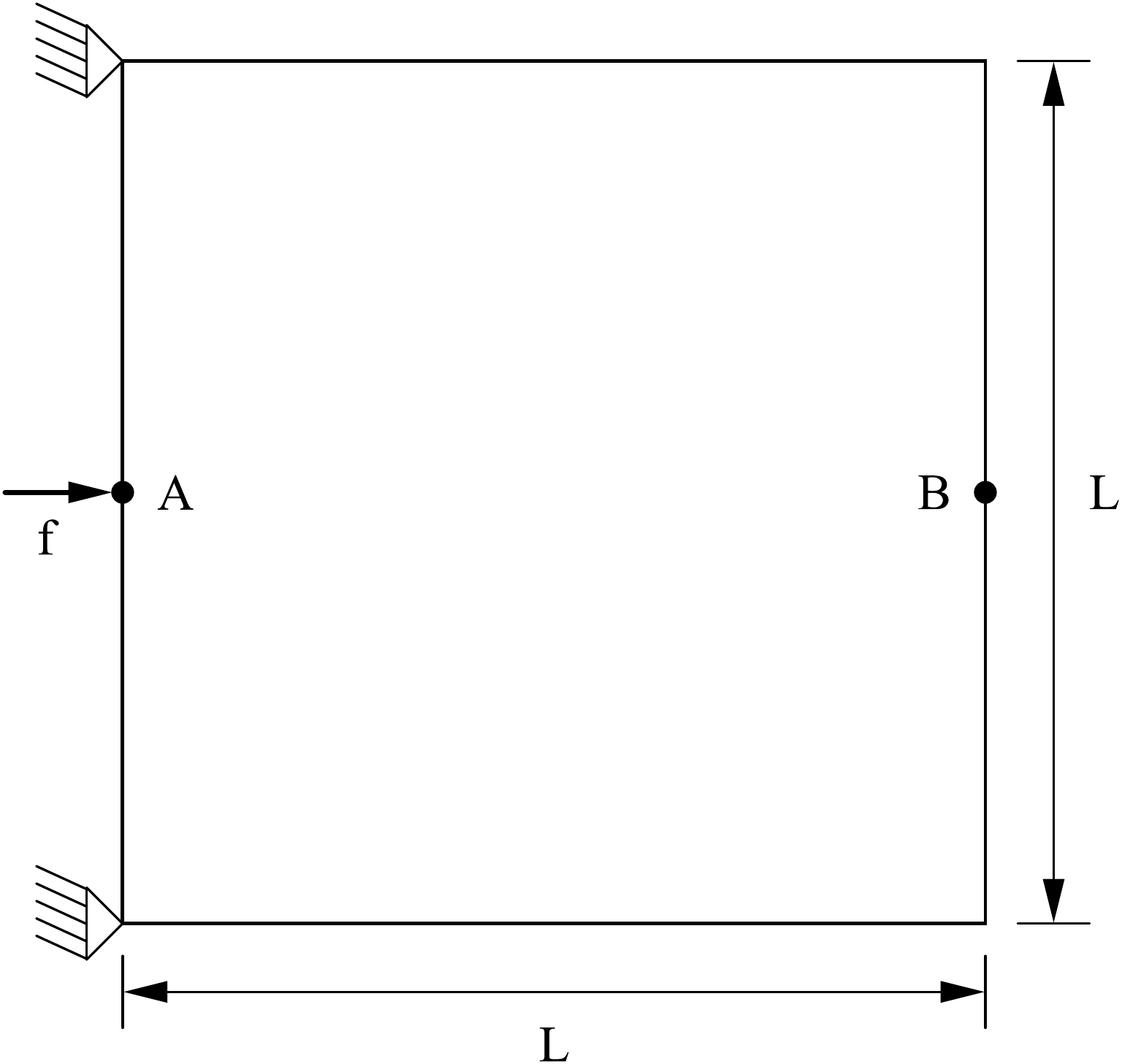}
        & \includegraphics[width=0.47\textwidth]{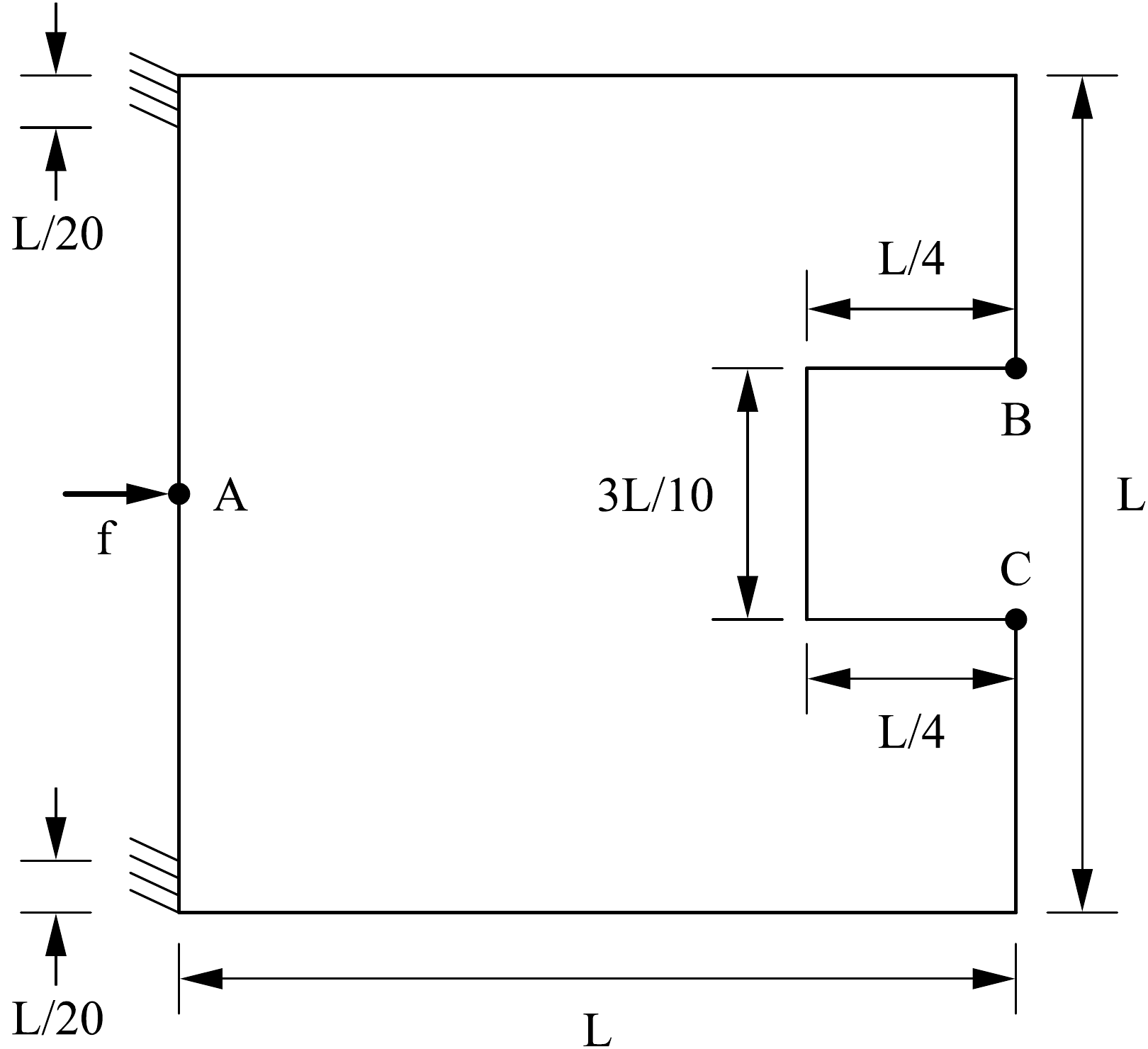} \\ 
      \end{tabular}
      \caption{Domains of the mechanisms: inverter (left) and gripper (right).} 
      \label{fig:domainMech}
\end{figure}

\noindent \textbf{Inverter.} Fig.~\ref{fig:domainMech}~(left) exhibits the square domain for the inverter, with $L = 300\,\mu m$ and thickness equal to $7\,\mu m$.
At the center of the left edge (point $A$), there is an external horizontal force of $50\,mN$ pointing to the right. Supports at the lower and the upper left corners prevent horizontal and vertical displacements. The Poisson's ratio and the Young's modulus of the material are, respectively, $0.3$ and $180\,mN/\mu m^{2}$. The stiffness of the springs of the model at points $A$ (input port) and $B$ (output port) are, respectively, $k_{in} = 4.0\,mN/\mu m$ and $k_{out} = 1.0\,mN/\mu m$. The volume of the optimal compliant mechanism is limited to $20\,\%$ of the domain's volume. Due to symmetry, only the upper half of the domain is discretized into $45000 \ (300 \times 150)$ finite elements. \\

\noindent \textbf{Gripper.} The design domain for the gripper is shown in Fig.~\ref{fig:domainMech}~(right), where $L = 320\,\mu m$. The domain has a thickness of $7\,\mu m$. In this compliant mechanism, an external force of $4\,mN$ is applied horizontally at point $A$. Two sets of supports (each one with $L/20\,\mu m$ of length) are located at the upper and lower left corners, preventing displacements. The Poisson's ratio and the Young's modulus of the material are, respectively, $0.3$ and $180\,mN/\mu m^{2}$. The stiffness of the spring at point $A$ (input port) is $k_{in} = 0.2\,mN/\mu m$. At points $B$ and $C$ (output ports), the stiffness is equal to $k_{out} = 1.0\,mN/\mu m$. Due to symmetry, only the upper half of the domain is discretized into $51200 \ (320 \times 160)$ finite elements. The optimal mechanism must contain $20\,\%$ of the domain's volume.

\subsection{Proposed strategies}

Whenever the objective function of problem (\ref{prob01}) needs to be computed, it is necessary to apply Newton's method to solve the nonlinear system (\ref{eq01b}), which in turn requires the solution of several linear systems in the form (\ref{eq02}). Therefore, the frequency of factorizations of $\KK$ is the key for speeding up the overall algorithm while keeping the solutions of the system (\ref{eq01b}) sufficiently accurate.

In the context of the nonlinear analysis of frame structures, Amir, Kirsch and Sheinman \cite{AmirKirschSheinman2008} presented some strategies based in the pure CA method reviewed in Section \ref{sec:inca}, in contrast to Newton's and the Modified Newton methods. They considered the application of the CA method using $2$ to $6$ vectors of the sequence (\ref{eq03}) to obtain the approximate solution (\ref{eq03aa}) of the linear system~(\ref{eq02}). 

In this work, we propose  the five \emph{Inexact Combined Approximations} (ICA) strategies described below. All of them rely on the use of the iterative scheme (\ref{eq03}) to approximately solve the linear system (\ref{eq02}) at each iteration of Newton's method. The strategies differ on the frequency of $\KK$ factorizations, on the updating scheme for matrix $\DeltaK$ and on the use of (\ref{eq03}) for the solution of the adjoint linear system (\ref{eq09}) related to the objective function gradient computation. To help future reference, each strategy is identified by its abbreviation.

\begin{description}

\item[\tt upK1:] $\KK$ is factored only at the first iteration of Newton's method; $\DeltaK$ is updated at each iteration of Newton's method. 

\item[\tt upK1g:] same as {\tt upK1}, with the addition that the linear system (\ref{eq09}) is also solved inexactly by means of the iterative scheme (\ref{eq03}).

\item[\tt upK100:] $\KK$ is factored only at the first iteration of Newton's method; $\DeltaK$ is updated every 100 iterations of Newton's method. 

\item[\tt upK100g:] same as {\tt upK100}, with the addition that the linear system (\ref{eq09}) is also solved inexactly by means of the iterative scheme (\ref{eq03}).

\item[\tt upK03K100g:] $\KK$ is factored only at every 3 iterations of the SPLP method; $\DeltaK$ is updated every 100 iterations of Newton's method, and the linear system (\ref{eq09}) is solved inexactly by means of the iterative scheme (\ref{eq03}). 

\end{description}

In order to evaluate the performance of the proposed strategies, they are compared with two factorization schemes commonly found in literature:

\begin{description}
\item[\tt N] (Newton): $\KK$ is factored at every iteration of Newton's method. 

\item[\tt MN] (Modified Newton): $\KK$ is factored only at the first iteration of Newton's method.
\end{description}

It is important to mention that, independently of the adopted strategy, in the first five iterations of the SPLP method, Newton's method is always applied to obtain the approximate solutions of the nonlinear system (\ref{eq01b}) to reduce the impact of the starting solution on the overall optimization process.

\subsection{Implementation Details}\label{sec:impdt}

The performance of the proposed strategies will be investigated for the four benchmark structures and compliant mechanisms presented in Subsection \ref{sub:testproblems}. 

When we apply an optimization algorithm for solving problem (\ref{prob01}), at every iteration we need to solve the nonlinear system (\ref{eq01b}) in order to compute the objective function. No matter the strategy adopted here, our stopping criterion for Newton's method is 
$\| \RR(\uu^{(\ell)},\,\bsrho) \|_{\infty} \leq 10^{-5}$.

For strategies {\tt upK1}, {\tt upK1g}, {\tt upK100}, {\tt upK100g} and {\tt upK03K100g}, in a fixed iteration $\ell$ of Newton's method, the linear system (\ref{eq02}) is solved inexactly by means of the ICA method described in Section \ref{sec:inca}. In the solution of this system, we used a stopping criterion based on the value of the relative magnitude of the residual of the linear system (\ref{eq02}), given by
\begin{equation*}
    \widehat{R}_{k}^{(\ell)} = \frac{\| \KK\sk + \RR \|_{\infty}}{\| \RR \|_{\infty}}, \qquad k = 0,\,1,\,2,\,\dots.
		\label{stopica1}
\end{equation*}
where $\KK \equiv \KK(\uu^{(\ell)},\,\bsrho)$ and $\RR = \RR(\uu^{(\ell)},\,\bsrho)$. If we get $\widehat{R}_{k^{*}}^{(\ell)} < \varepsilon_{R}$ (where $\varepsilon_{R}$ is a small positive constant) for some $k^{*} \leq 10$, we consider that $\svec^{(\ell)} = \svec^{(k^{*})}$ is a good approximation for the solution of (\ref{eq02}). Otherwise, we update the factorization of $\KK$, and solve the linear system (\ref{eq02}) exactly. 

When strategies {\tt upK1g}, {\tt upK100g} and {\tt upK03K100g} are adopted, the ICA method is also applied for inexactly solving the linear system (\ref{eq09}). In this case, the relative magnitude of the residual is defined as
\begin{equation*}
    T_{k} = \frac{\| \KKhat\bslambdak + \lbold \|_{\infty}}{\| \lbold \|_{\infty}}, \qquad k = 0,\,1,\,2,\,\dots,
		\label{stopica2}
\end{equation*}
and we require that $T_{\overline{k}} < \varepsilon_{T}$ for some $\overline{k} \leq 10$. If this criterion is not satisfied, we update the factorization of $\KKhat$, and solve the linear system (\ref{eq09}) exactly. In our numerical tests, we took $\varepsilon_{R} = 10^{-2}$ and $\varepsilon_{T} = 10^{-8}$. A more stringent value is adopted for $\varepsilon_{T}$ in order to guarantee a high level of accuracy in the computation of the derivatives of the objective function.

With the aim of ensuring the global convergence of Newton's method and avoiding small decreases of $\| \RR(\uu,\,\bsrho) \|^{2}_{2}$ even with large steps, for all of the strategies we adopt the Armijo's line search~\cite{Nocedal2006}.

When dealing with topology optimization problems under geometric nonlinearity, the value of the penalty parameter $p$ of the SIMP method plays a crucial role in the stabilization of the overall optimization process. Based on our previous experience~\cite{GS2014}, the initial value of $p$ is set to $1.0$ and it is increased by $\Delta p = 0.1$ after every 10 iterations of the SPLP method until it reaches $3.0$. The minimum value $\rho_{\min}$ allowed for the densities in problem (\ref{prob01}) was set to $0.001$, to avoid the singularity of the matrix $\KK$. The density filter proposed by Bruns and Tortorelli \cite{BrunsTortorelli2003} was applied to prevent the checkerboard-like pattern in the optimal distribution of material \cite{DiazSigmund1995}. In our tests, the number of elements used to define the filter radius were 10, 5, 7.5 and 5, respectively, for the cantilever beam, the slender beam, the inverter and the gripper. 

The algorithm presented was coded in C$++$. The CPLEX 12.1 software library was used for solving the subproblems of the SPLP algorithm, and the linear systems were solved using the CHOLMOD 1.7 library. All of the tests were performed on a personal computer with an Intel Core i7-6500U processor, under the Ubuntu Linux 16.4 operating system.

%===========================================
\section{Performance of the strategies considering a fixed budget}\label{sec:perfstr}

With the aim of providing a fair comparison of the performance of the several strategies proposed here, we have first solved all of the problems with a fixed budget of 300 iterations of the SPLP method. The tables presented in this and in the following sections contain the value reached for the objective function $F(\bsrho)$, the accumulated number of iterations of Newton's method and the total time spent to achieve the optimal solution and to evaluate $F(\bsrho)$, that is composed by the sum of the time consumed for assembling $\KK$ and its factorization, in the assembling of the right-hand side (RHS) vectors, and in the solution of the linear systems (\ref{eq02}) and (\ref{eq09}). Moreover, the tables show the time spent for solving the SPLP subproblems, for applying the filter to the densities, and for computing the gradient vector $\nabla F(\bsrho)$. We remark that the CPU times are always expressed in seconds. 

\begin{table}[htbp!]
    \centering
    \caption{Results obtained for 300 iterations of SPLP (cantilever beam and slender beam).} \vspace{0.3cm}
    \label{table_100_iter_cb_sb}
        \begin{tabular}{lc@{}c@{}c@{}c@{}c@{}c@{}c@{}}
        \hline\noalign{\smallskip}
        \textbf{Cantilever beam} & \multirow{2}{*}{\tt N} & \multirow{2}{*}{\tt MN} & \multirow{2}{*}{\tt upK1} & \multirow{2}{*}{\tt upK1g} & \multirow{2}{*}{\tt upK100} & \multirow{2}{*}{\tt upK100g} & {\tt upK03} \\
        $\mathbf{400 \times 100}$ &  &  &  &  &  &  & {\tt K100g} \\ 
        \noalign{\smallskip}\hline\noalign{\smallskip} 
        $F(\bsrho)$ - value           & $2036.42$ \,  & $2036.62$ \, & $2036.42$ \, & $2036.57$ \,  & $2036.62$ \,& $2036.57$ \, & $2036.49$ \\
        \# Newton iter.          & $679$       & $1147$    & $826$     & $814$       & $839$     & $814$     & $904$     \\ 
         \hline
        \noalign{\smallskip}
        CPU time (s)  &  &  &  &  &  &  &  \\ 
        Total                & $1579.93$   & $1210.52$ & $1208.69$ & $1029.80$   & $1150.91$ & $986.30$  & $732.38$  \\
	    $F(\bsrho)$             & $1218.75$   & $868.40$  & $847.89$  & $688.44$    & $801.32$  & $637.12$  & $373.62$  \\     
        $\KK$                 & $181.18$    & $184.92$  & $197.34$  & $199.15$    & $145.51$  & $148.82$  & $95.27$   \\
        RHS                  & $25.37$     & $37.25$   & $32.55$   & $39.17$     & $32.69$   & $39.34$   & $56.28$   \\
	    Factorizations         & $1000.31$   & $628.47$  & $597.50$  & $415.95$    & $602.27$  & $414.50$  & $161.29$  \\
	    Linear systems         & $11.89$     & $17.76$   & $20.50$   & $34.17$     & $20.85$   & $34.46$   & $60.79$   \\
	    $\nabla F(\bsrho)$     & $104.95$    & $102.12$  & $105.70$  & $102.45$    & $106.19$  & $107.59$  & $104.99$  \\
	    SPLP solving         & $159.17$    & $145.15$  & $158.65$  & $143.24$    & $146.34$  & $143.54$  & $155.87$  \\
	    Filtering             & $90.76$     & $88.78$   & $90.23$   & $89.47$     & $90.87$   & $91.75$   & $91.82$   \\
	    Other                 & $6.30$      & $6.07$    & $6.22$    & $6.20$      & $6.20$    & $6.30$    & $6.08$    \\
	    \hline\hline\noalign{\smallskip}      
                \textbf{Slender beam} & \multirow{2}{*}{\tt N} & \multirow{2}{*}{\tt MN} & \multirow{2}{*}{\tt upK1} & \multirow{2}{*}{\tt upK1g} & \multirow{2}{*}{\tt upK100} & \multirow{2}{*}{\tt upK100g} & {\tt upK03} \\
        $\mathbf{600 \times 75}$ &  &  &  &  &  &  & {\tt K100g} \\ 
        \noalign{\smallskip}\hline\noalign{\smallskip} 
        $F(\bsrho)$ - value           & $138.62$   & $138.52$   & $138.68$     & $138.77$   & $138.68$   & $138.77$  & $138.62$ \\
        \# Newton iter.          & $648$      & $758$      & $658$        & $667$      & $658$      & $670$     & $838$    \\
         \hline
        \noalign{\smallskip}
        CPU time (s)  &  &  &  &  &  &  &  \\         
        Total                 & $1622.93$  & $1262.76$  & $1293.72$    & $1201.68$  & $1263.94$  & $1170.56$ & $876.84$ \\
	    $F(\bsrho)$            & $1208.74$  & $853.88$   & $878.07$     & $770.80$   & $846.39$   & $741.24$  & $463.44$ \\       
        $\KK$                & $189.38$   & $210.10$   & $191.42$     & $199.88$   & $157.05$   & $162.87$  & $106.44$ \\
        RHS                  & $26.97$    & $30.74$    & $29.07$      & $43.71$    & $29.01$    & $46.34$   & $59.35$  \\
	    Factorizations           & $979.93$   & $599.03$   & $640.57$     & $480.06$   & $643.33$   & $479.84$  & $232.96$ \\
	    Linear systems         & $12.46$    & $14.01$    & $17.01$      & $47.15$    & $17.00$    & $52.19$   & $64.69$  \\
	    $\nabla F(\bsrho)$    & $34.74$    & $35.33$    & $35.09$      & $35.12$    & $35.06$    & $35.37$   & $35.35$  \\
	    SPLP solving        & $351.85$   & $345.64$   & $352.93$     & $367.64$   & $354.74$   & $365.88$  & $350.45$ \\
	    Filtering           & $25.71$    & $25.97$    & $25.73$      & $26.23$    & $25.86$    & $26.16$   & $25.69$  \\
	    Other                 & $1.89$     & $1.94$     & $1.90$       & $1.89$     & $1.89$     & $1.91$    & $1.91$   \\	    
        \noalign{\smallskip}\hline
    \end{tabular}
\end{table}

From Table \ref{table_100_iter_cb_sb}, one may notice that the optimal values of the objective function among the strategies are very similar for the cantilever beam and for the slender beam. With Newton's method, the evaluation of the objective function takes $77.1\,\%$ of the total time demanded for solving the problem for the cantilever beam, and $74.5\,\%$ for the slender beam. This percentage is substantially reduced if the {\tt upK03K100g} strategy is adopted. 
From the total time spent to obtain the optimal structure, the computation of $F(\bsrho)$ consumes 
 $51.0\,\%$ for the cantilever beam, and $52.9\,\%$ for the slender beam.

A slight reduction on the time consumed to assemble $\KK$ is noted if {\tt upK100} is adopted, in comparison with the {\tt upK1} strategy. For these two strategies, and for both structures, the time spent on the factorization of $\KK$ remains very similar. On the other hand, with the adoption of the  strategies {\tt upK1g} and {\tt upK100g}, we observe a significant reduction of the time spent in the factorization of $\KK$ as compared with the counterparts {\tt upK1} and  {\tt upK100} -- about $30\,\%$ for the cantilever beam, and $25\,\%$ for the slender beam.

As we relax the proposed strategies, it is noticeable the decreasing of the  demanded time for computing the objective function, and consequently, for solving the problem. Indeed, all of the five ICA strategies outperformed the original and the modified Newton methods.  For the cantilever beam, the reduction attained on the time spent to evaluate $F(\bsrho)$ by strategies {\tt upK1}, {\tt upK1g}, {\tt upK100} and {\tt upK100g}, in comparison with Newton's method, is of $30.4\,\%$, $43.5\,\%$, $34.3\,\%$ and $47.7\,\%$, respectively.  The reduction is even more  pronounced if we consider the {\tt upK03K100g} strategy , reaching $69.3\,\%$ for this structure.  When it comes to the slender beam, similar percentage reductions have been achieved, namely $27.4\,\%$, $36.2\,\%$, $30.0\,\%$, $38.7\,\%$, and  $61.7\,\%$ for strategies {\tt upK1}, {\tt upK1g}, {\tt upK100}, {\tt upK100g}, and  {\tt upK03K100g}, respectively.

\begin{figure}[htbp!]
   \centering 
      \begin{tabular}{cc} 
        \includegraphics[width=0.47\textwidth]{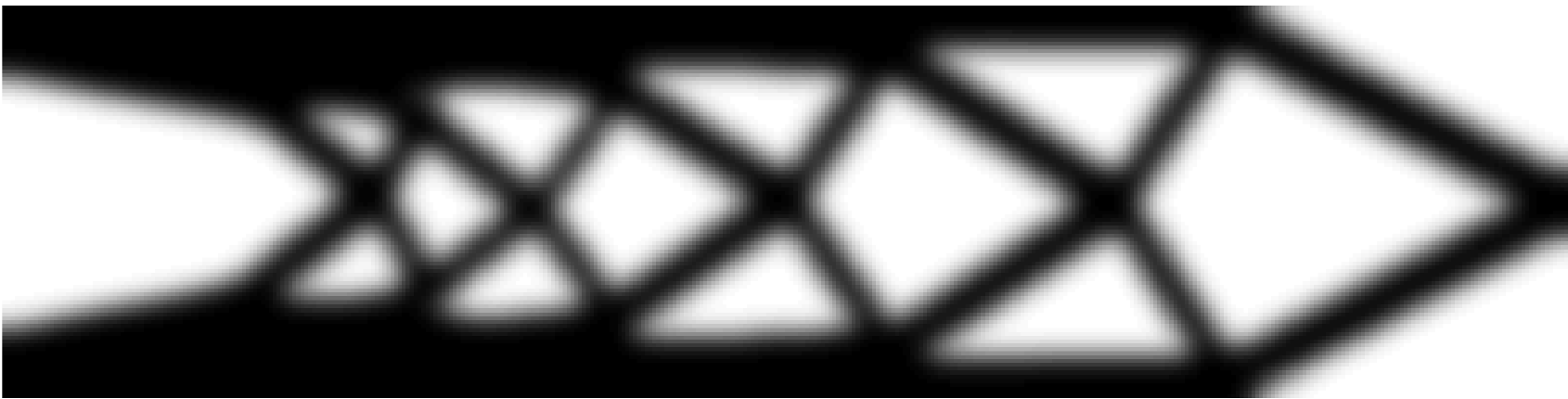}
        &\includegraphics[width=0.47\textwidth]{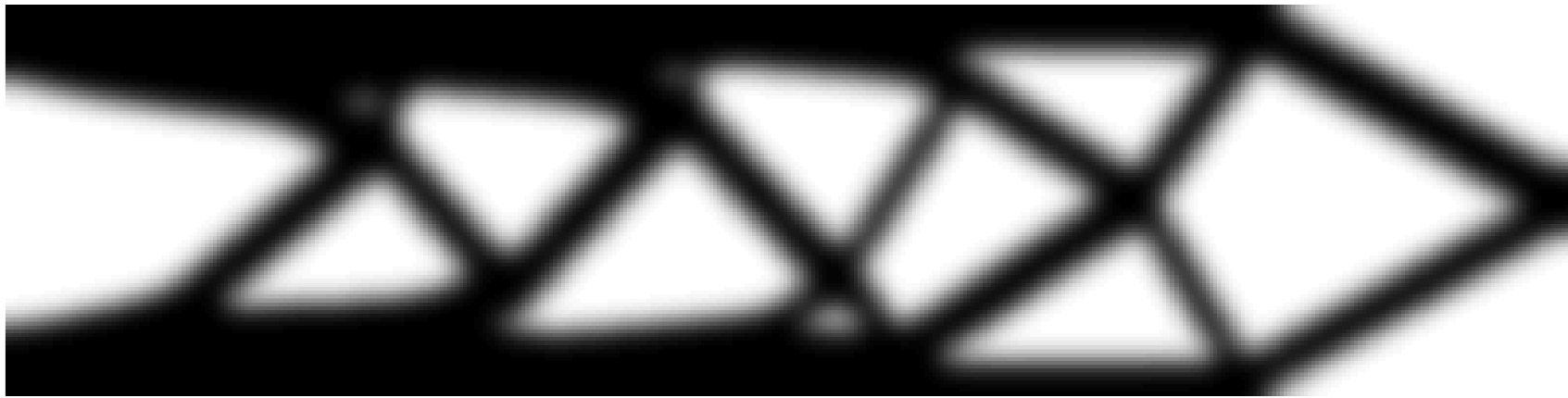} \\[15 pt]
          \includegraphics[width=0.47\textwidth]{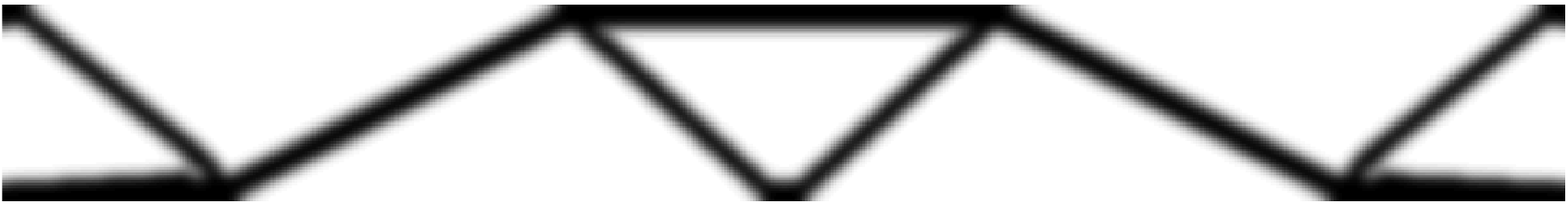}
        & \includegraphics[width=0.47\textwidth]{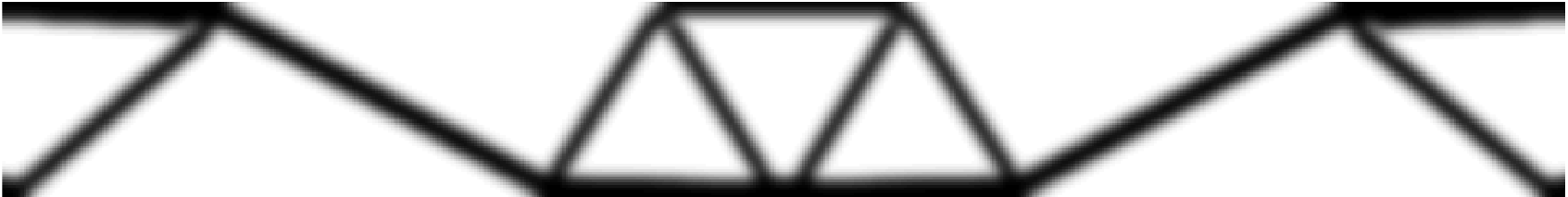} \\[6 pt] 
          Small displacements         & Large displacements 
      \end{tabular}
      \caption{Optimal topologies for the cantilever beam $400 \times 100$ (top) and the slender beam $600 \times 75$ (bottom).} 
      \label{figure_cb_sb}
\end{figure}

To put the effect of the geometric nonlinearity in perspective, we have also considered the two structures as described in Section~\ref{sub:testproblems}, but using the \emph{linear elastic} assumption to define the stiffness matrix, so that it varies just with the densities, and do not depend on the displacements any longer. 
Adopting the policy of factoring such a matrix at every call of the objective function, i.e. the \emph{original} strategy analyzed in~\cite{SGS2019}, we have solved the corresponding topology optimization problem, reaching the so-called \emph{small displacement} configuration.  
Fig.~\ref{figure_cb_sb} depicts the optimal topologies, with small and large displacements, for the cantilever beam and the slender beam, respectively.  For the cantilever beam, it is evident the lack of symmetry about the medium horizontal line presented by the structure under large displacements, to cope with the nonlinear requirement, as compared with the structure under small displacements.  When it comes to the slender beam, it is also significant the effect of the geometric nonlinearity upon the optimal configuration. It must be stressed that the structures obtained for the seven strategies presented on Table \ref {table_100_iter_cb_sb} are quite similar, so only the results for the {\tt upK100g} strategy are shown.

\begin{figure}[htbp!]
   \centering 
      \begin{tabular}{cc} 
           \includegraphics[width=0.36\textwidth]{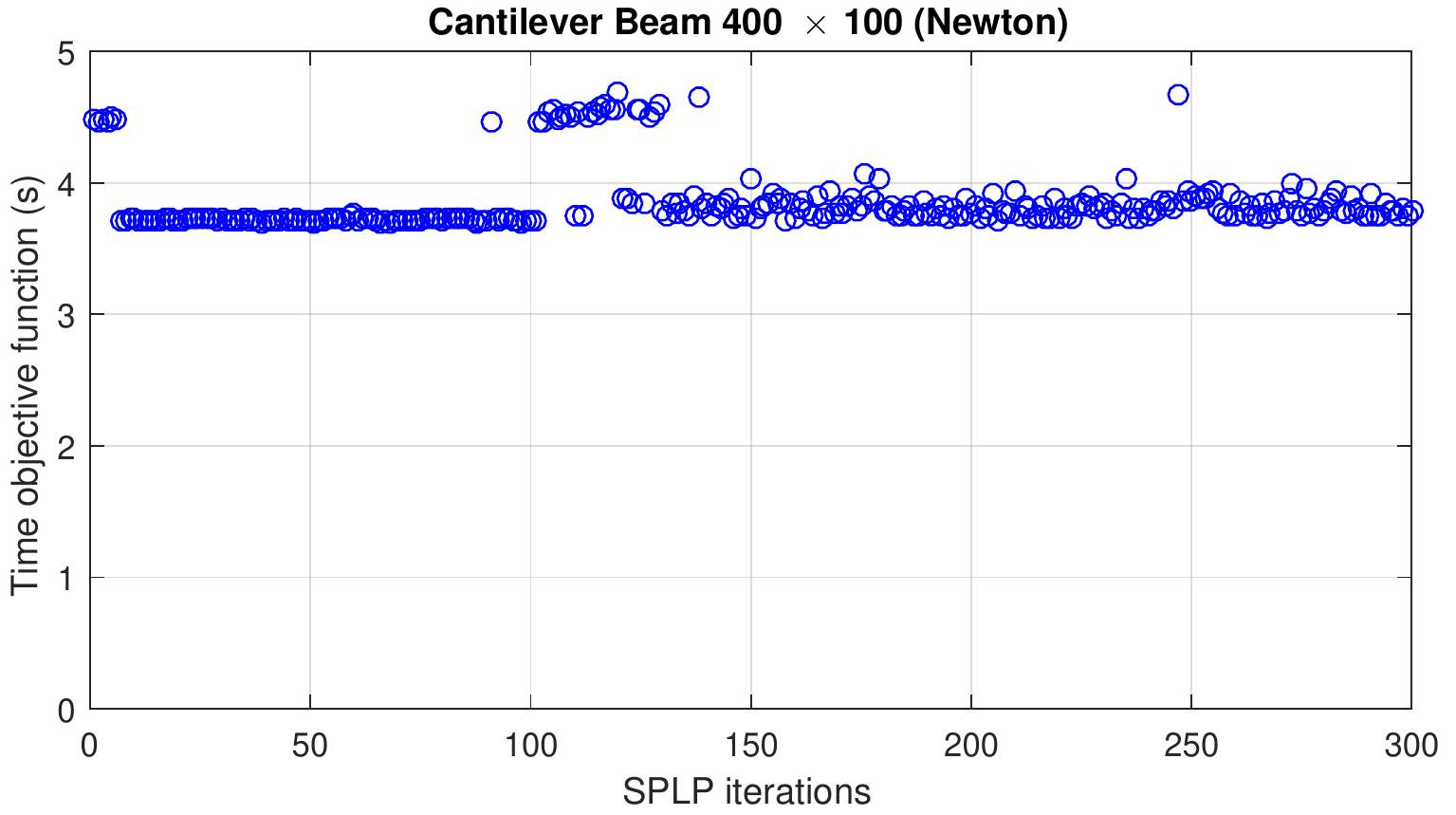}
        & \includegraphics[width=0.36\textwidth]{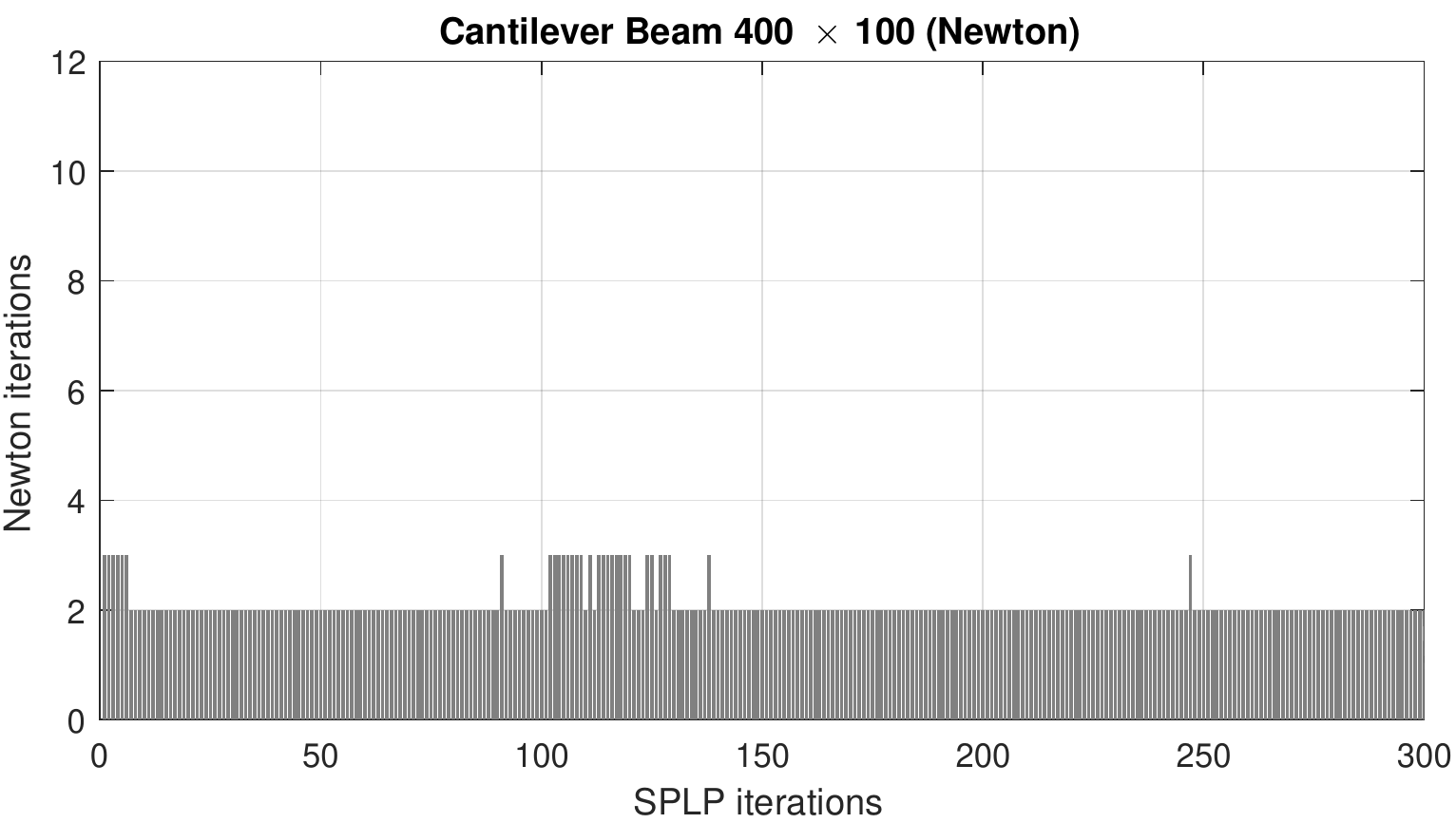} \\
	    \includegraphics[width=0.36\textwidth]{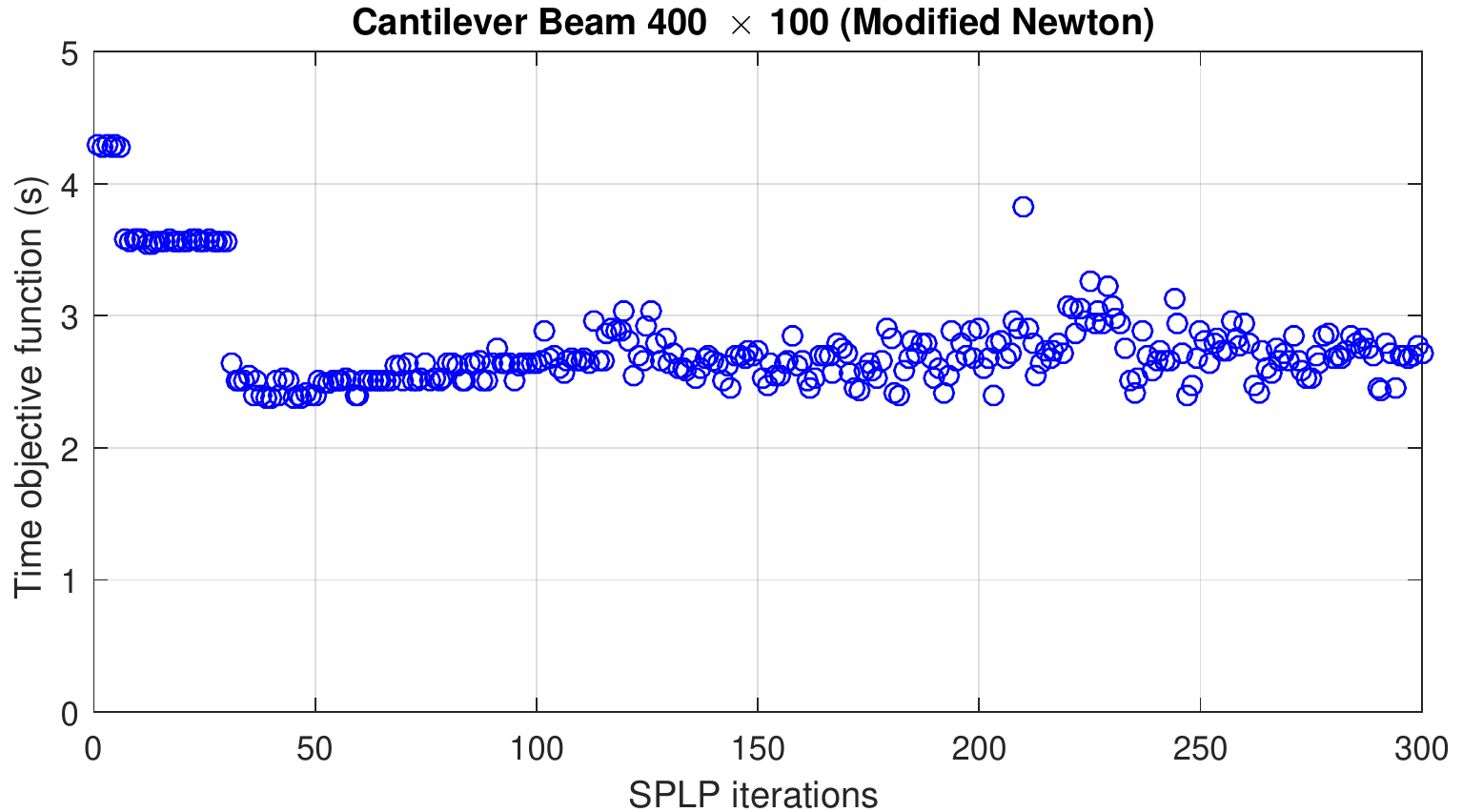}
        & \includegraphics[width=0.36\textwidth]{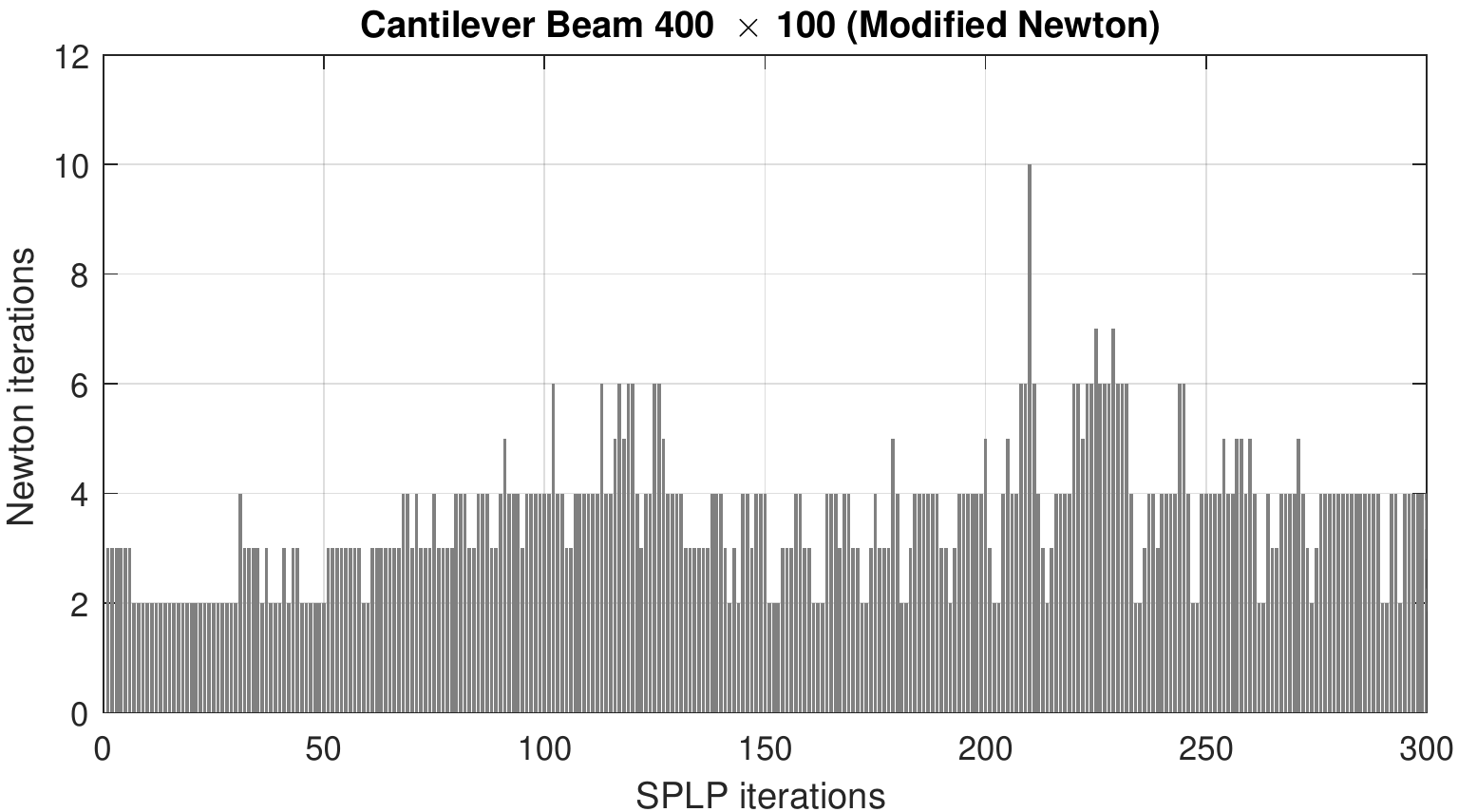} \\
	     \includegraphics[width=0.36\textwidth]{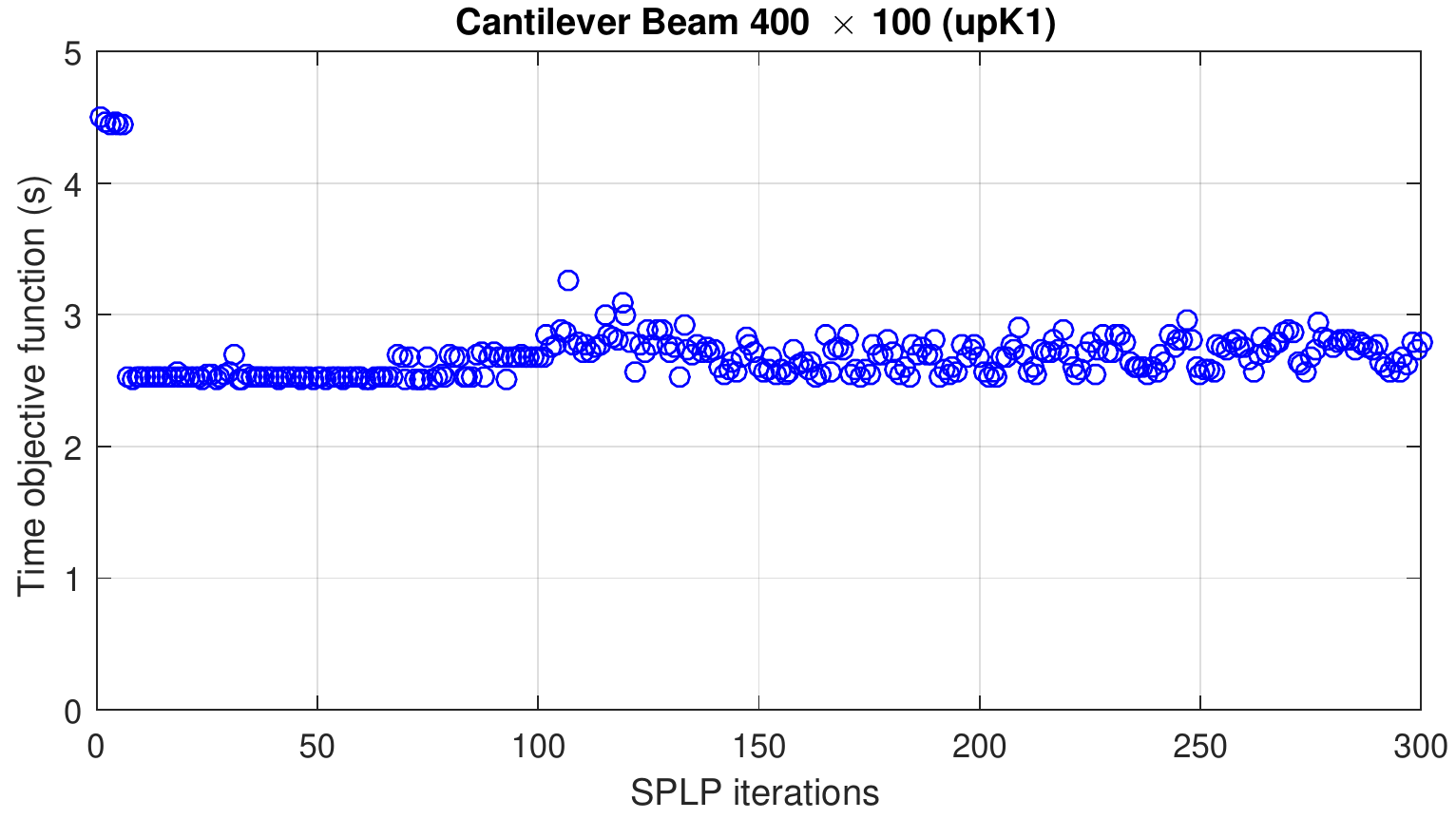}
        & \includegraphics[width=0.36\textwidth]{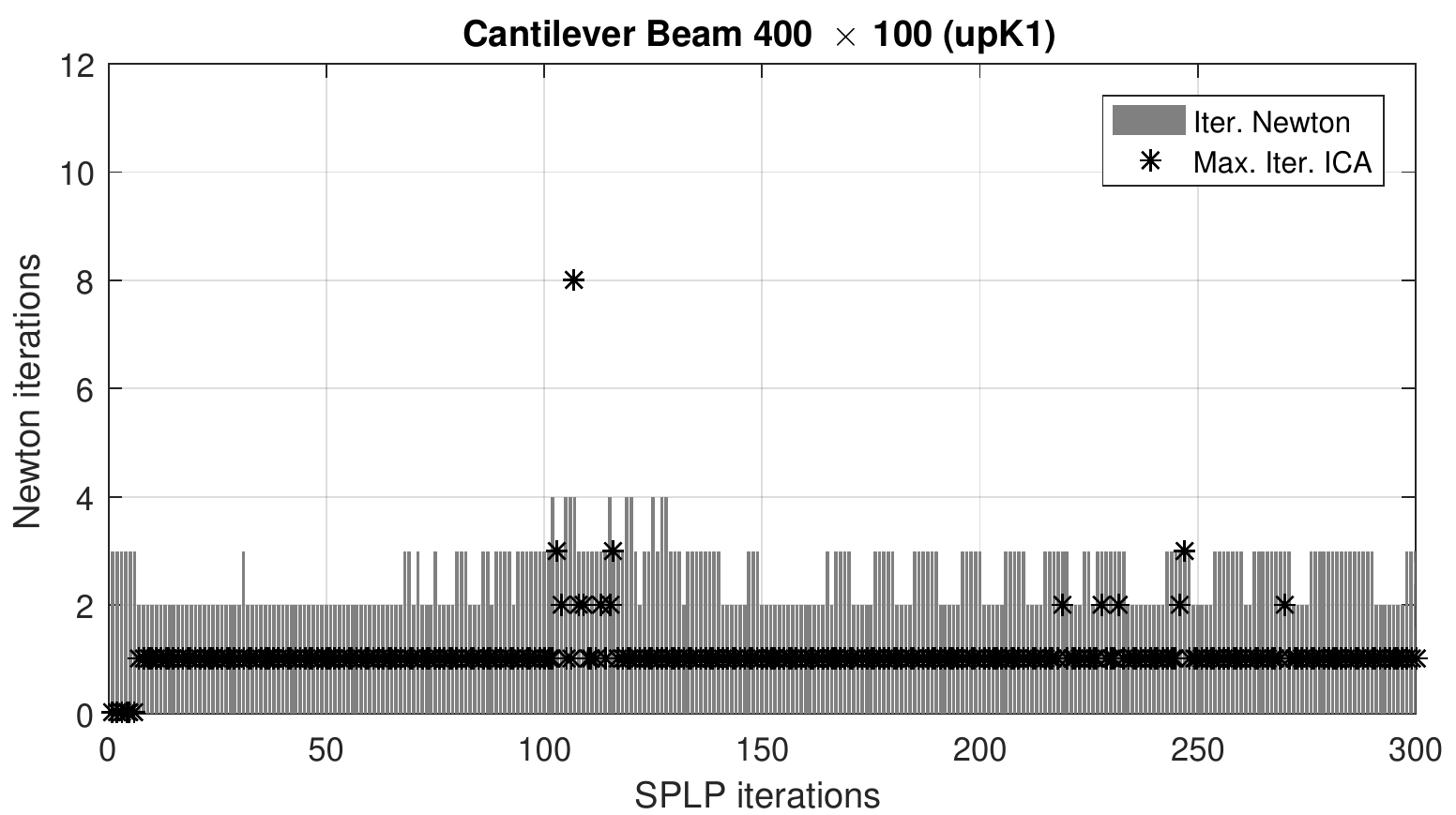} \\
	     \includegraphics[width=0.36\textwidth]{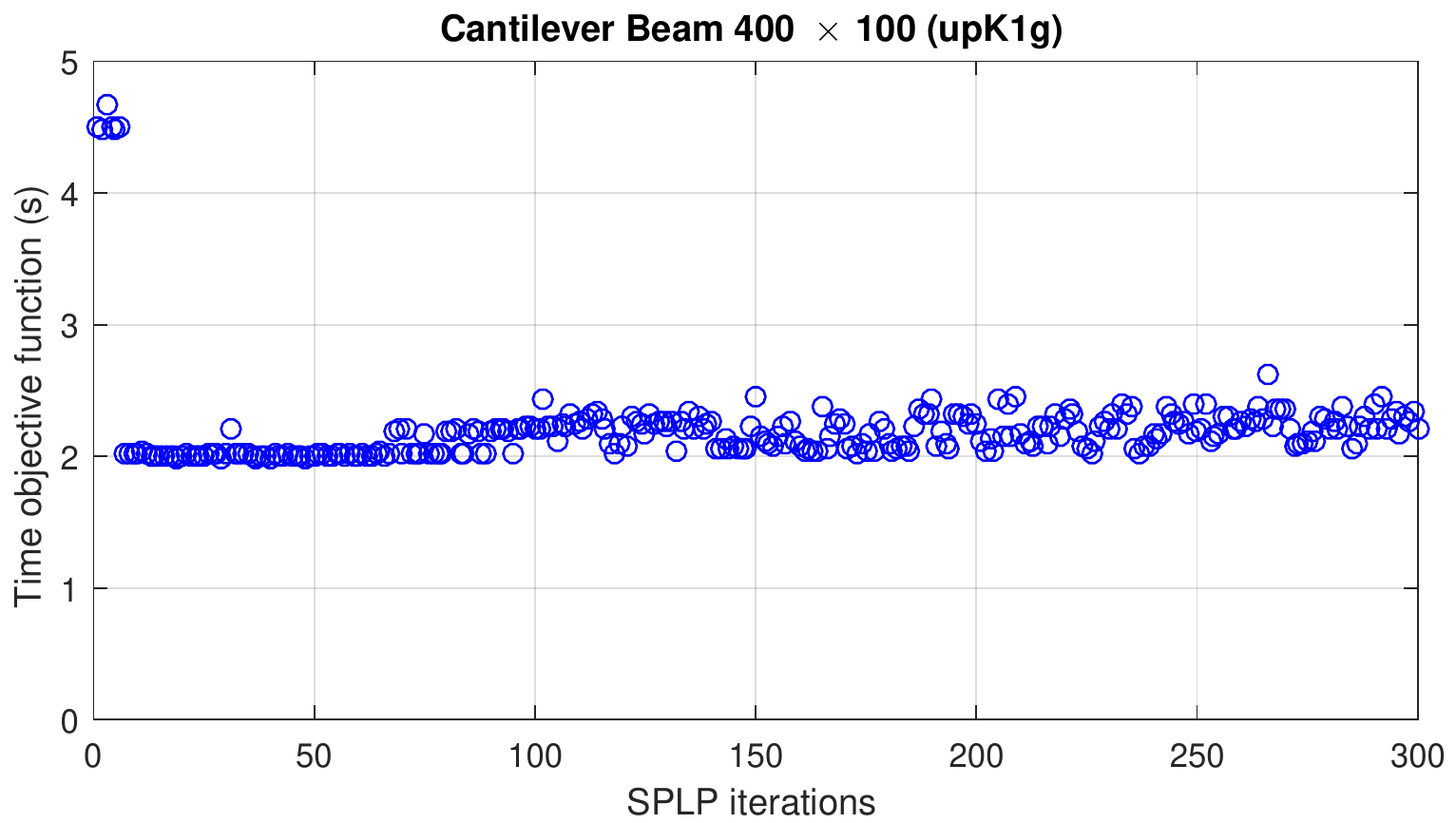}
        & \includegraphics[width=0.36\textwidth]{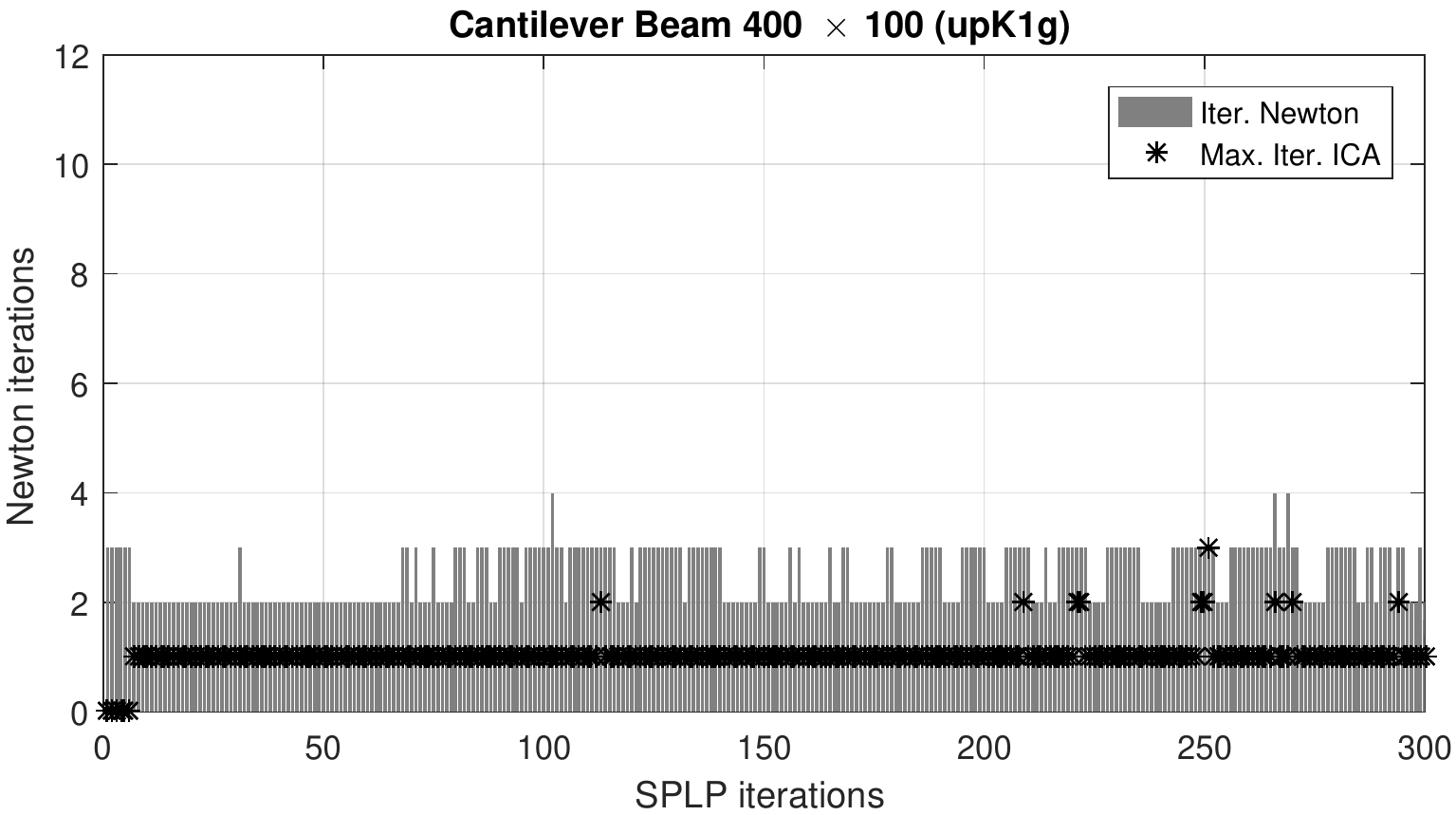} \\
	     \includegraphics[width=0.36\textwidth]{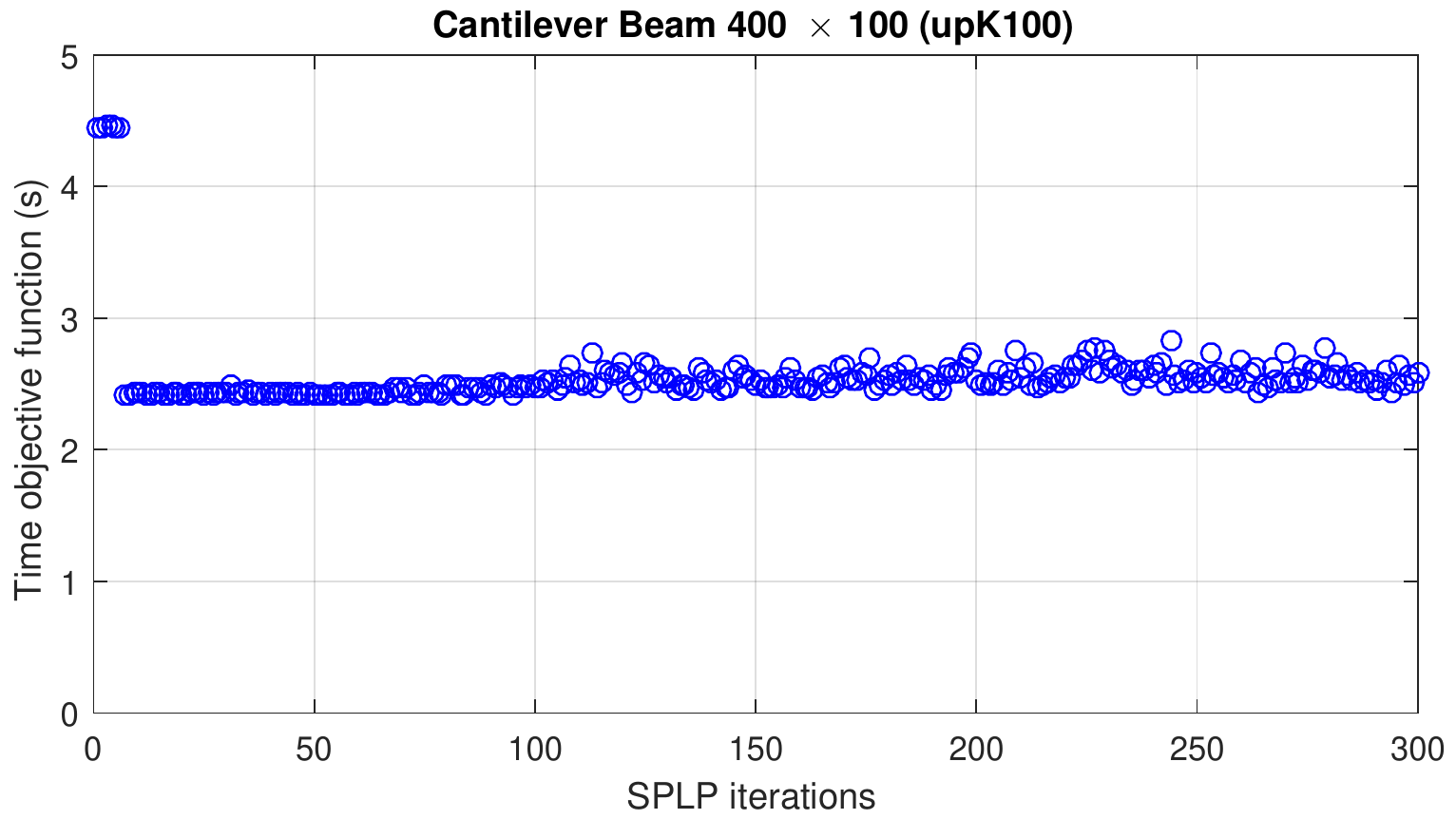}
        & \includegraphics[width=0.36\textwidth]{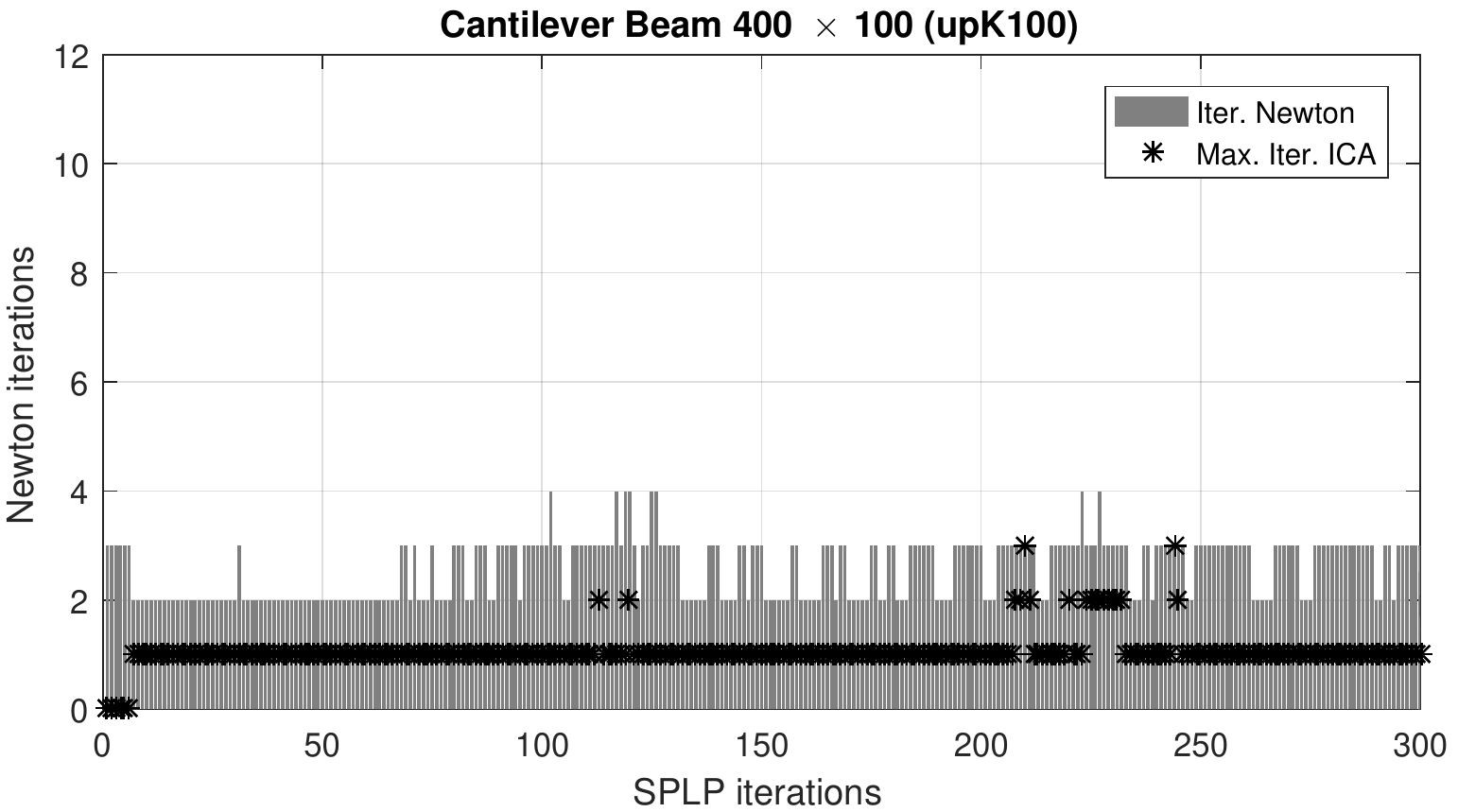} \\
	    \includegraphics[width=0.36\textwidth]{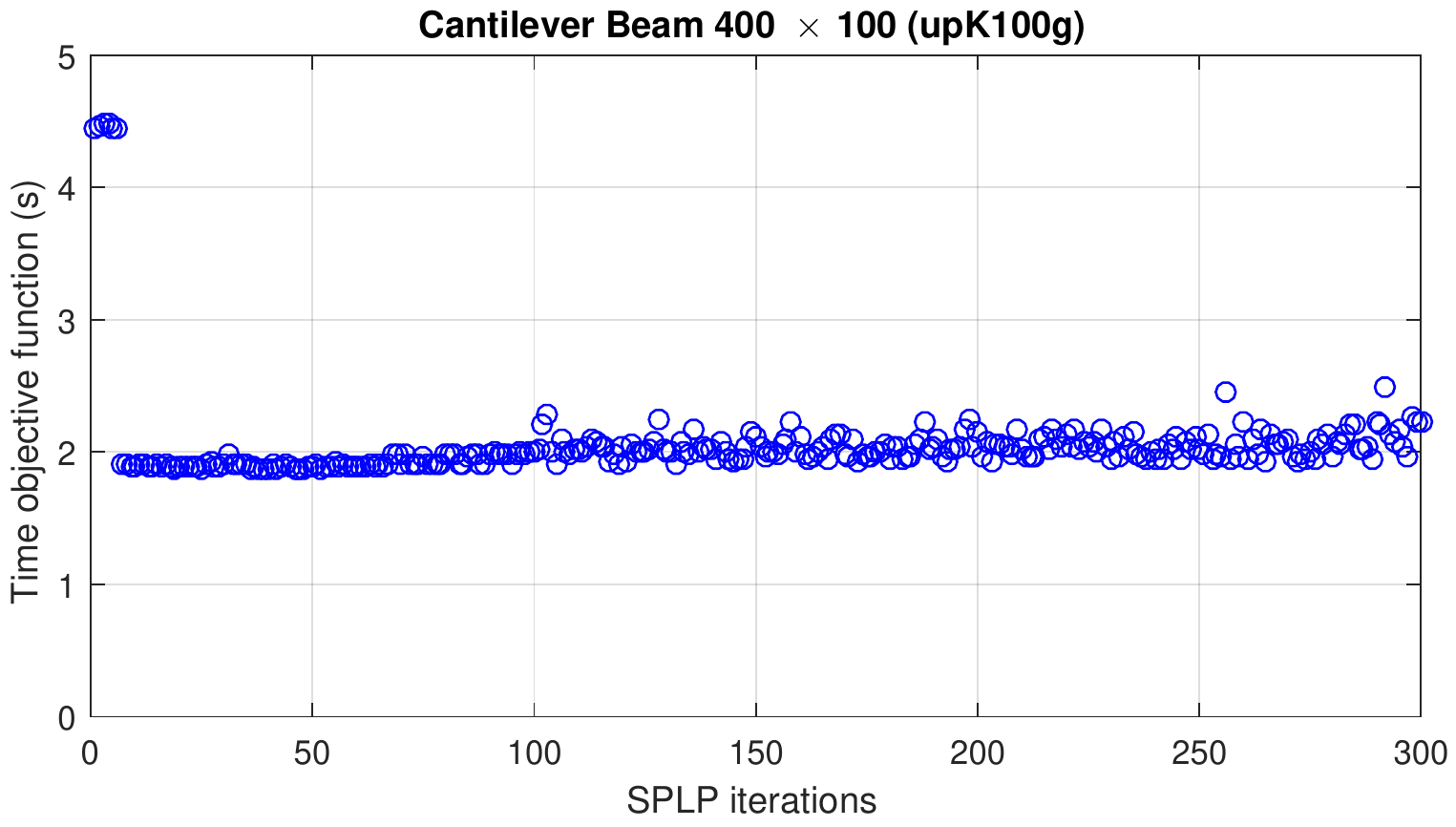}
        & \includegraphics[width=0.36\textwidth]{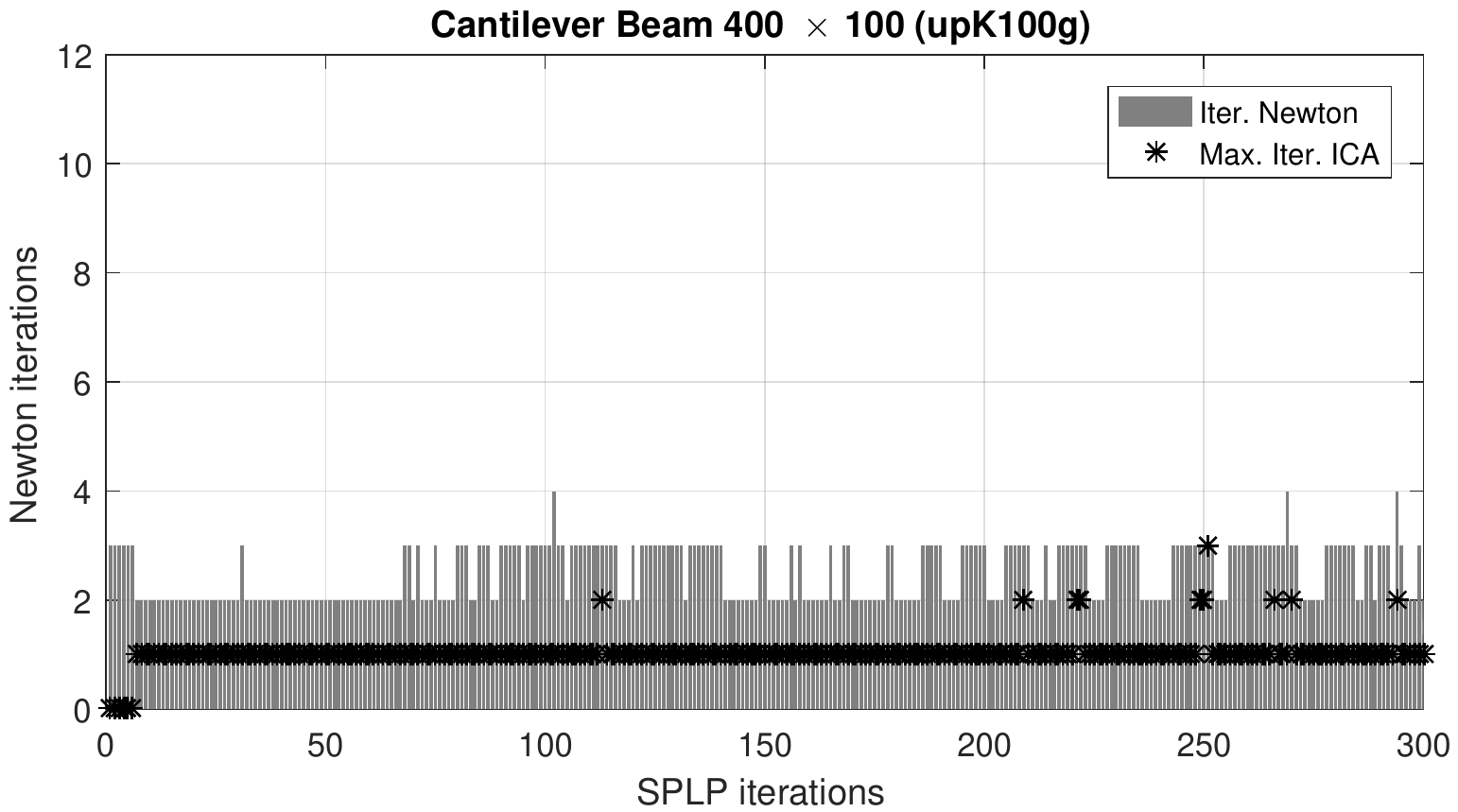} \\
	     \includegraphics[width=0.36\textwidth]{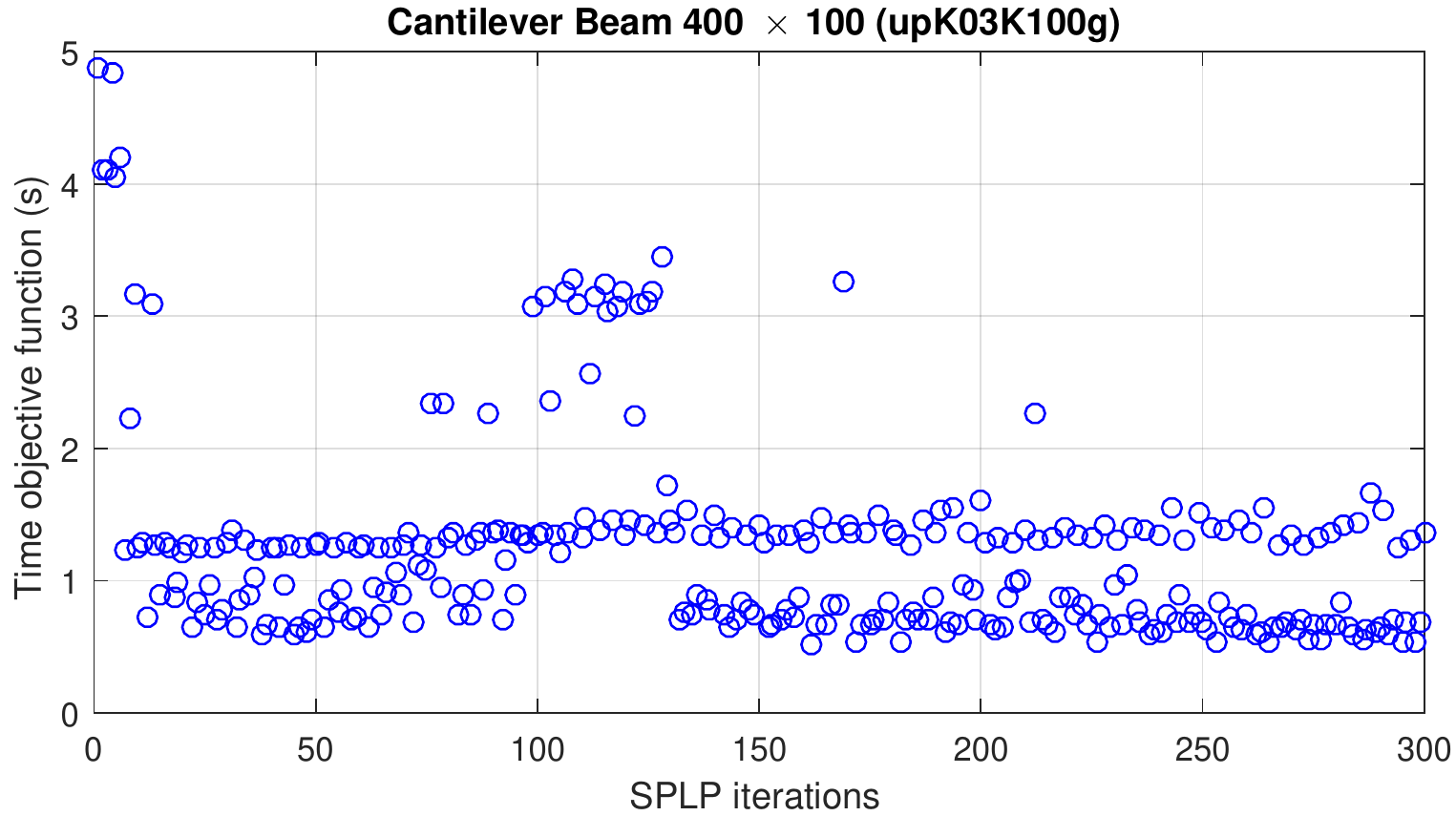}
        & \includegraphics[width=0.36\textwidth]{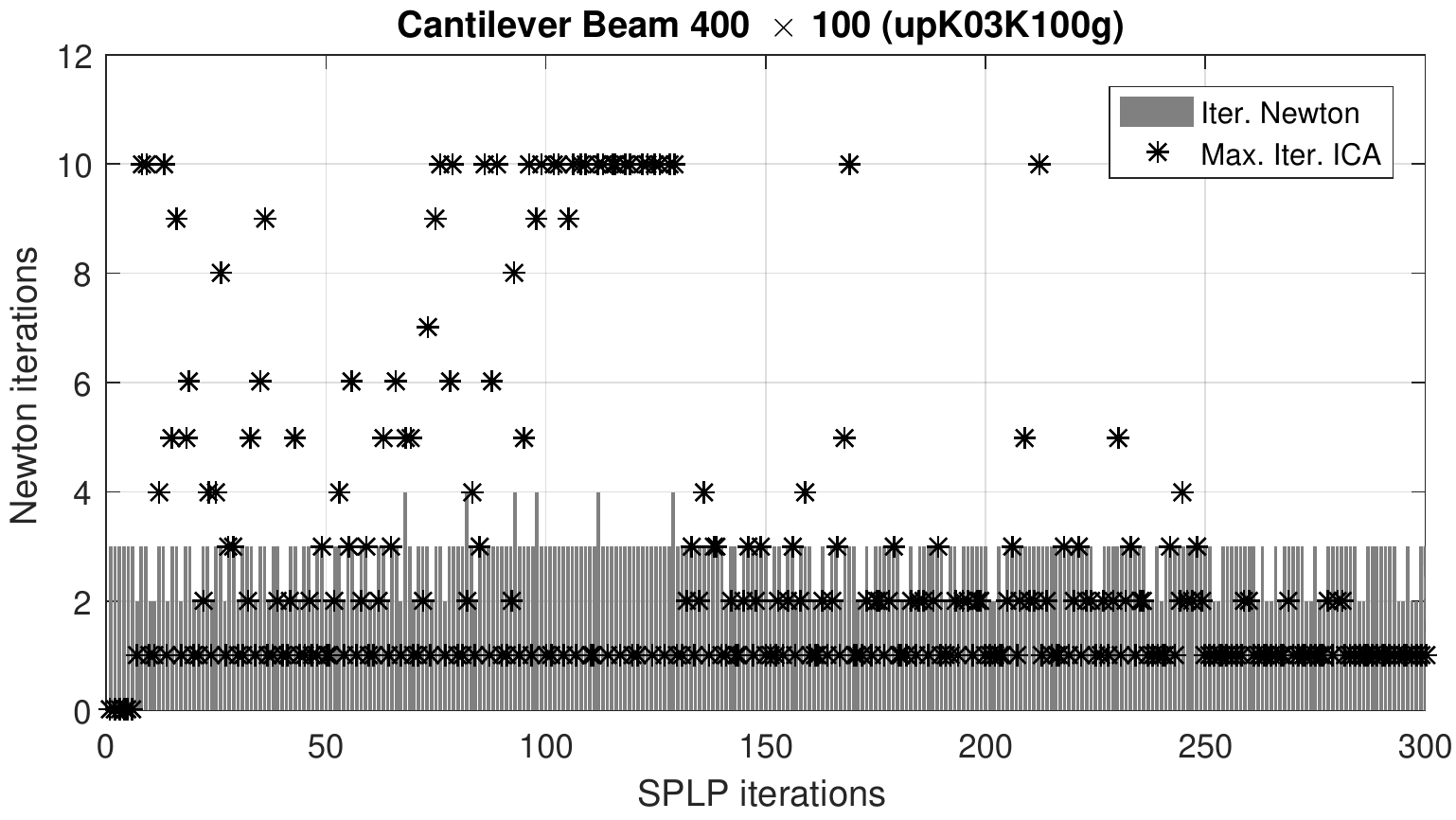} \\
      \end{tabular}
      \caption{Time spent to compute the objective function (left) and number of iterations of Newton's Method (right) in the solution of the cantilever beam ($400 \times 100$).} 
      \label{time_fobj_iter_newton_cb400x100}
\end{figure}

Figure \ref{time_fobj_iter_newton_cb400x100} shows the time consumed by each strategy to evaluate the objective function, and the number of iterations performed by Newton's method at each iteration of the SPLP method, together with the maximum number of iterations performed by the iterative scheme (\ref{eq03}) to approximately solve the linear system (\ref{eq02}), for the cantilever beam. 

We can observe from Fig. \ref{time_fobj_iter_newton_cb400x100} that, if Newton's method is employed, it takes at most 3 iterations to reach the approximate solution of the nonlinear system~(\ref{eq01b}), but, due to the expensive iterations, the average time spent to compute the objective function is about $4.0$ seconds. With the adoption of the Modified Newton method, we can note a good decrease in this average time as compared with Newton's method (it oscillates between $2.5$ and $3.0$ seconds), but a slower convergence is experienced in some occasions. For example, at the iteration $210$ of the SPLP method,  $10$ iterations were necessary to obtain the approximate solution of (\ref{eq01b}). This behavior can be explained by the fact that, for the Modified Newton method, the matrix $\KK$ is factored only at the initial point, without carrying any information about the changes of $\KK$ along the iterative process. 

The {\tt upK1} strategy represents an enhancement of the Modified Newton method, since the iterative scheme (\ref{eq03}) takes into account the variations on $\KK$ by means of the matrix $\DeltaK$. This is reflected in the expressive acceleration in the convergence of such inexact Newton method, as we can see in Fig.~\ref{time_fobj_iter_newton_cb400x100}. Moreover, in almost all of the callings of this inexact Newton method, the approximate solution of the linear system (\ref{eq02}) with the desired precision was obtained already with the first term of the sequence (\ref{eq03}). Despite of these facts, the time consumed to evaluate $F(\bsrho)$ by the Modified Newton method and by strategy {\tt upK1} are quite similar. On the other hand, strategy {\tt upK1g} showed a noticeable reduction in the time spent to compute $F(\bsrho)$ in comparison with {\tt upK1}, since it does not require the factorization of $\KK$ for the solution of the adjoint linear system (\ref{eq09}). Moreover, we observed the fast convergence of this inexact Newton method. According to Fig. \ref{time_fobj_iter_newton_cb400x100}, the average time consumed to compute $F(\bsrho)$ is about $2.2$ seconds when strategy {\tt upK1g} is applied, a reduction of 43,5\,\% when compared to Newton's method.

When the {\tt upK100} strategy is employed, a $26.2\,\%$  reduction is attained in the time spent to assemble $\KK$, as compared with {\tt upK1}. The decrease reached in this step by {\tt upK100g} is of $25.3\,\%$  in comparison with {\tt upK1g}. Strategies {\tt upK1g} and {\tt upK100g} required much less time to evaluate $F(\bsrho)$ compared to Newton's method, showing a reduction of $43.5\,\%$ and $47.8\,\%$, respectively.

The adoption of strategy {\tt upK03K100g} produced the most remarkable improvement in the time spent to evaluate the objective function, in relation to Newton's method. As we can see in Fig. \ref{time_fobj_iter_newton_cb400x100}, in most of the iterations, the evaluation of $F(\bsrho)$ demanded less than $1.5$ seconds, requiring less than $1$ second in many cases, while Newton's method spent about $4$ seconds to compute the same function. It should be mentioned that the maximum number of iterations of the iterative scheme (\ref{eq03}) can be reached when we are far from the solution, due to the lack of information on the changes of the factorization of $\KK$. Whenever this happens, the factorization of $\KK$ is updated, and we solve the full linear system (\ref{eq02}). Fortunately, this rarely occurs, so the speed of convergence of the inexact Newton method is almost not affected. According to Table \ref{table_100_iter_cb_sb}, strategy {\tt upK03K100g} requires $41.4\,\%$ less time than {\tt upK100g} to compute $F(\bsrho)$. This reduction increases to 69.3\,\% when the comparison is made with Newton's method. Taking also into account the decreasing in the time spent with the factorization of $\KK$, the reduction are even more noticeable, reaching $61.1\,\%$ and $83.9\,\%$ when compared to the {\tt upK100g} strategy and Newton's method, respectively.

Now, we will discuss the results obtained for the inverter and for the gripper mechanisms, considering the fixed budget of 300 iterations of the SPLP method. These results are summarized on Table \ref{table_100_iter_inv_grp}. One may notice that the maximum variation between the optimal values of the objective function was of $0.09\,\%$ for the inverter and of $0.26\,\%$ for the gripper, which is very reasonable from the practical point of view. As expected, when Newton's method is employed, the factorization of $\KK$ takes most of the total time for solving the problems. In fact, the factorization of $\KK$ consumes $90\,\%$ and $86.5\,\%$ of the total time for the inverter and for the gripper, respectively. When the {\tt upK03K100g} strategy is adopted, these values are reduced to $41.1\,\%$ and $35.9\,\%$ for the inverter and  the gripper, respectively.

\begin{table}[htbp!]
    \centering
    \caption{Results obtained for 300 iterations of SPLP (Inverter and Gripper).} \vspace{0.3cm}
    \label{table_100_iter_inv_grp}
        \begin{tabular}{lc@{}c@{}c@{}c@{}c@{}c@{}c@{}}
        \hline\noalign{\smallskip}
                        \textbf{Inverter} & \multirow{2}{*}{\tt N} & \multirow{2}{*}{\tt MN} & \multirow{2}{*}{\tt upK1} & \multirow{2}{*}{\tt upK1g} & \multirow{2}{*}{\tt upK100} & \multirow{2}{*}{\tt upK100g} & {\tt upK03} \\
        $\mathbf{300 \times 150}$ &  &  &  &  &  &  & {\tt K100g} \\
        \noalign{\smallskip}\hline\noalign{\smallskip} 
        $F(\bsrho)$ - value              & $-8.6969$ \,  & $-8.6989$ \,  & $-8.7020$ \, & $-8.7035$ \, & $-8.6997$  \,& $-8.7028$ \,& $-8.7046$ \\
        \# Newton iter.            & $696$       & $1065$      & $859$     & $849$      & $848$     & $841$     & $929$     \\
       \hline
        \noalign{\smallskip}
        CPU time (s)  &  &  &  &  &  &  &  \\ 
        Total                 & $1883.89$   & $1334.16$   & $1344.36$ & $1097.37$  & $1293.52$ & $1066.04$ & $756.60$  \\
	    $F(\bsrho)$           & $1691.98$   & $1144.86$   & $1149.91$ & $905.20$   & $1102.38$ & $860.02$  & $564.05$  \\       
        $\KK$                  & $161.20$    & $194.80$    & $178.60$  & $179.35$   & $128.35$  & $132.06$  & $106.75$  \\
        RHS                    & $28.87$     & $38.90$     & $42.83$   & $54.16$    & $42.22$   & $54.37$   & $66.69$   \\
	    Factorizations        & $1487.78$   & $891.76$    & $895.51$  & $614.71$   & $899.80$  & $616.35$  & $310.82$  \\
	    Linear systems          & $14.13$     & $19.40$     & $32.97$   & $56.98$    & $31.98$   & $57.24$   & $79.79$   \\
	    $\nabla F(\bsrho)$   & $35.44$     & $36.17$     & $37.69$   & $36.37$    & $36.27$   & $36.61$   & $36.57$   \\
	    SPLP solving       & $128.50$    & $124.56$    & $126.68$  & $126.87$   & $126.22$  & $140.25$  & $127.22$  \\
	    Filtering              & $26.23$     & $26.77$     & $28.17$   & $27.09$    & $26.88$   & $27.29$   & $26.90$   \\
	    Other                & $1.74$      & $1.80$      & $1.91$    & $1.84$     & $1.77$    & $1.87$    & $1.86$    \\	      
        \hline\hline\noalign{\smallskip}
        \textbf{Gripper} & \multirow{2}{*}{\tt N} & \multirow{2}{*}{\tt MN} & \multirow{2}{*}{\tt upK1} & \multirow{2}{*}{\tt upK1g} & \multirow{2}{*}{\tt upK100} & \multirow{2}{*}{\tt upK100g} & {\tt upK03} \\
        $\mathbf{320 \times 160}$ &  &  &  &  &  &  & {\tt K100g} \\
        \noalign{\smallskip}\hline\noalign{\smallskip} 
        $F(\bsrho)$ - value            & $-3.8481$ & $-3.8478$   & $-3.8556$  & $-3.8493$ & $-3.8577$ & $-3.8487$ & $-3.8460$  \\
        \# Newton iter.      & $659$     & $888$       & $732$      & $725$     & $732$     & $744$     & $864$      \\
         \hline
        \noalign{\smallskip}
        CPU time (s)  &  &  &  &  &  &  &  \\ 
        Total                   & $2530.57$ & $1828.94$   & $1817.90$  & $1488.72$ & $1766.46$ & $1437.66$ & $999.59$   \\
	    $F(\bsrho)$              & $2188.84$ & $1485.33$   & $1477.07$  & $1143.88$ & $1425.90$ & $1090.35$ & $650.17$   \\       
        $\KK$                      & $208.06$  & $237.01$    & $218.70$   & $222.47$  & $169.91$  & $174.48$  & $128.89$   \\
        RHS                       & $32.40$   & $39.31$     & $41.90$    & $52.61$   & $42.26$   & $53.62$   & $72.89$    \\
	    Factorizations            & $1932.43$ & $1189.14$   & $1186.75$  & $816.90$  & $1183.57$ & $809.60$  & $358.68$   \\
	    Linear systems            & $15.95$   & $19.87$     & $29.72$    & $51.48$   & $30.16$   & $52.65$   & $89.71$    \\
	    $\nabla F(\bsrho)$       & $102.56$  & $104.83$    & $103.80$   & $106.30$  & $104.30$  & $106.54$  & $108.18$   \\
	    SPLP solving              & $147.57$  & $145.99$    & $145.25$   & $145.73$  & $144.11$  & $147.71$  & $148.28$   \\
	    Filtering                 & $85.64$   & $86.76$     & $85.83$    & $86.43$   & $86.19$   & $86.64$   & $86.44$    \\
	    Other                    & $5.96$    & $6.03$      & $5.95$     & $6.38$    & $5.96$    & $6.42$    & $6.52$     \\	
        \noalign{\smallskip}\hline 
    \end{tabular}
\end{table} 

Analogously to the results achieved for the rigid structures, the reductions attained in the time spent to compute $F(\bsrho)$\, are remarkable. Besides, {\tt upK03K100g} is much cheaper than the other strategies. As an example, for the inverter, strategies {\tt upK1}, {\tt upK1g}, {\tt upK100}, {\tt upK100g} and {\tt upK03K100g} spent, respectively, $32.0\,\%$, $46.5\,\%$, $34.8\,\%$, $49.2\,\%$ and $66.7\,\%$ less time to evaluate $F(\bsrho) $ in comparison with Newton's method. 

\begin{figure}[htbp!]
   \hspace*{-0.2cm}
      \begin{tabular}{cccc} 
          \includegraphics[width=0.22\textwidth]{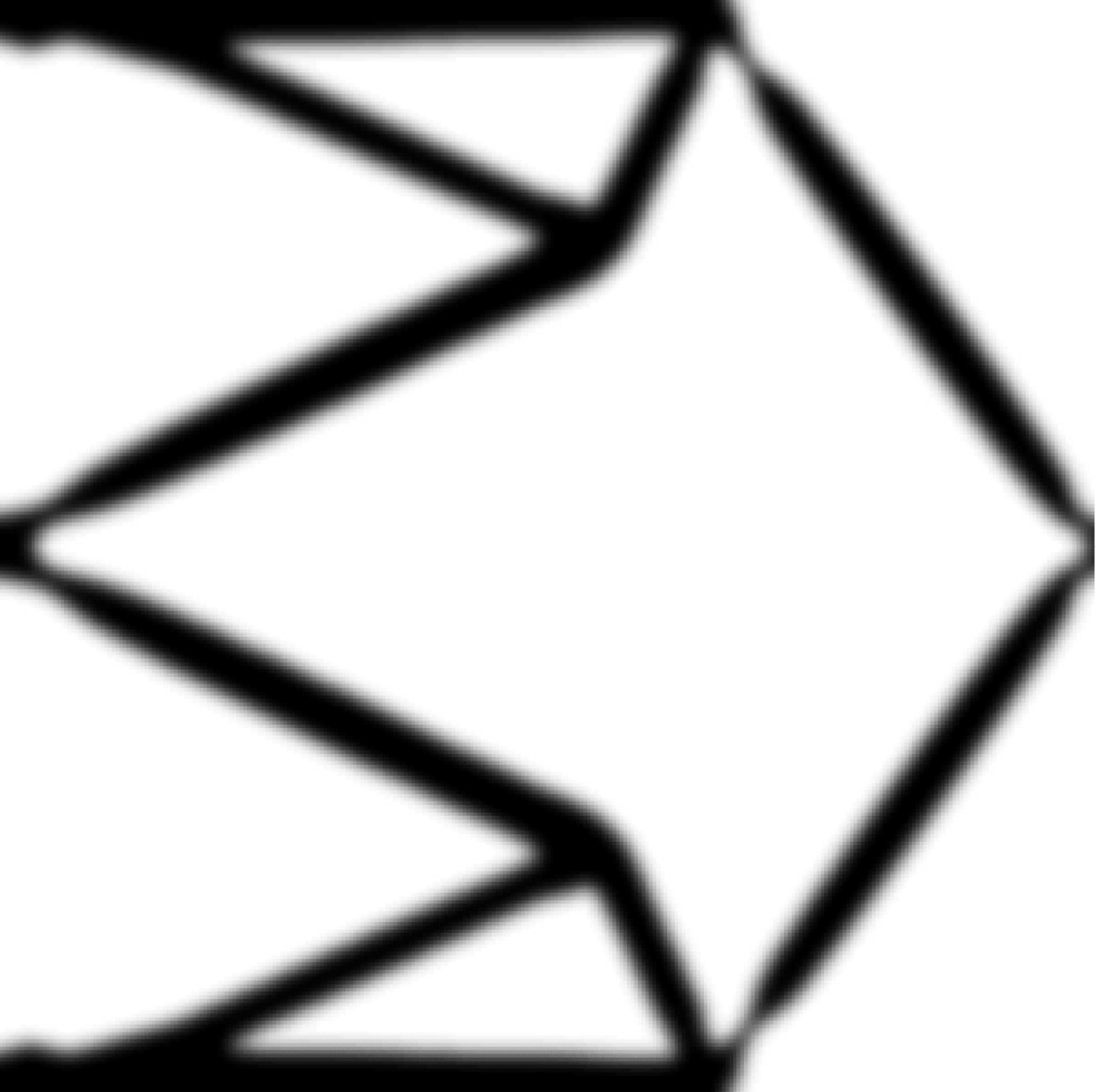}
        & \includegraphics[width=0.22\textwidth]{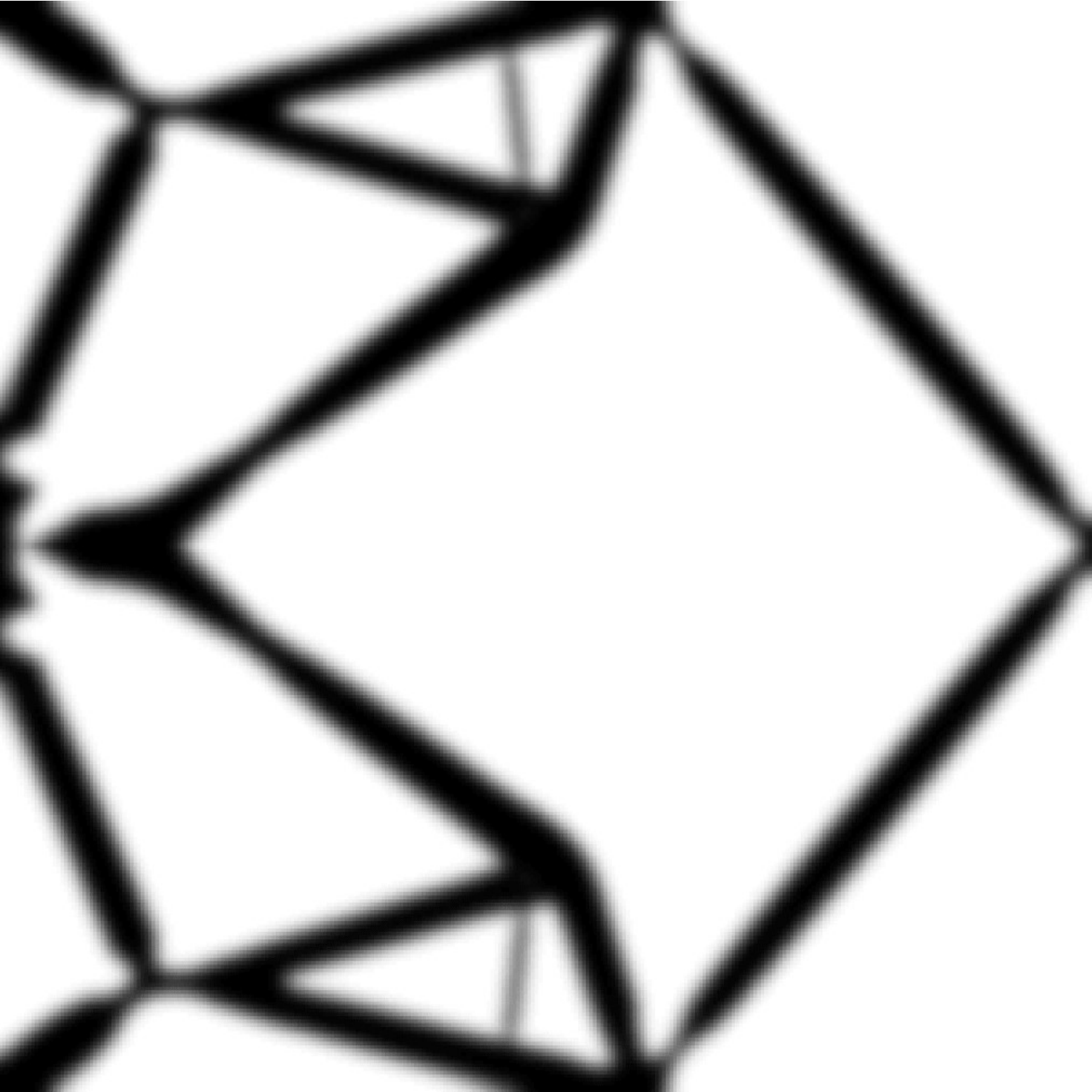}
        & \includegraphics[width=0.22\textwidth]{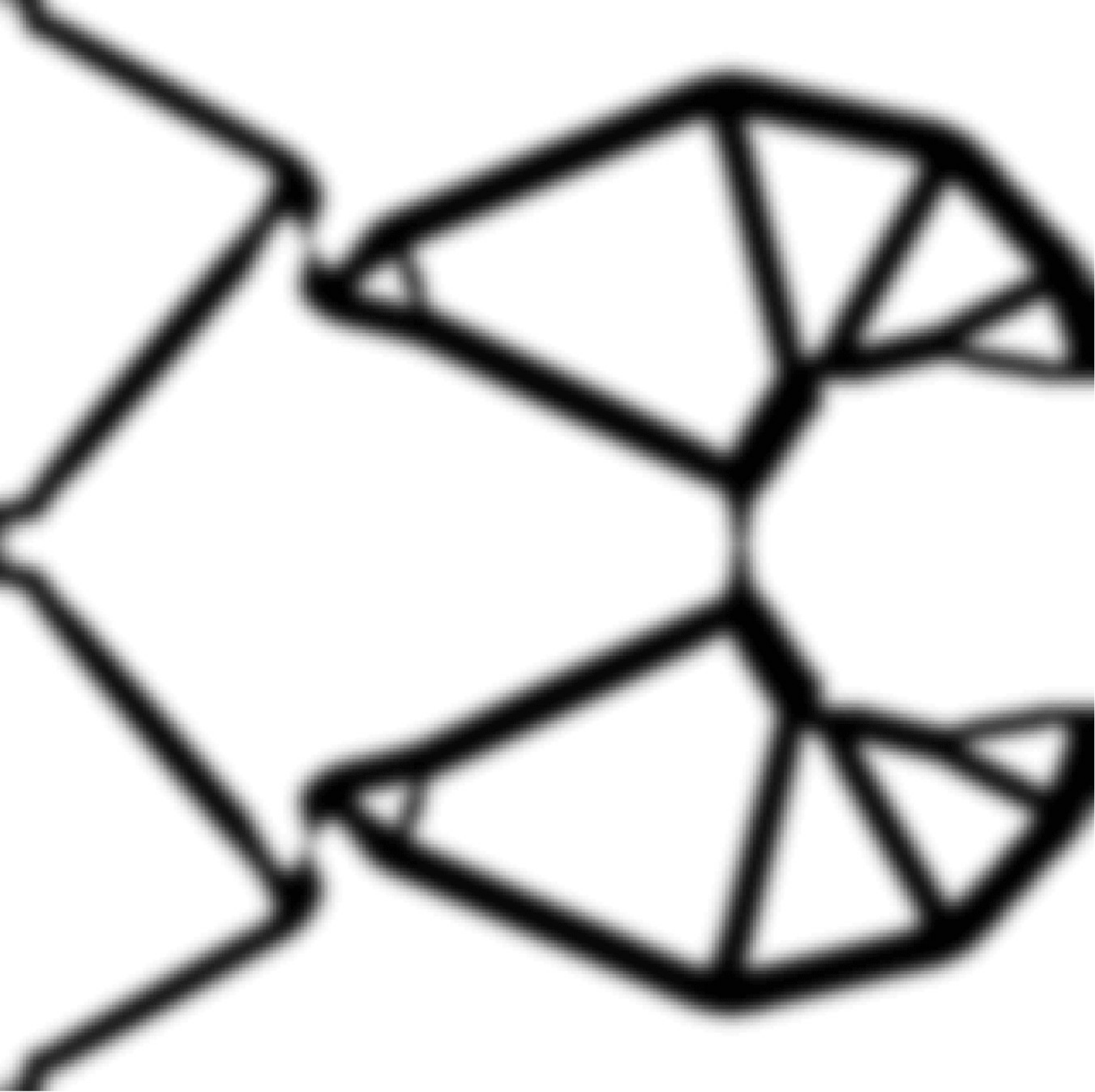}
        & \includegraphics[width=0.22\textwidth]{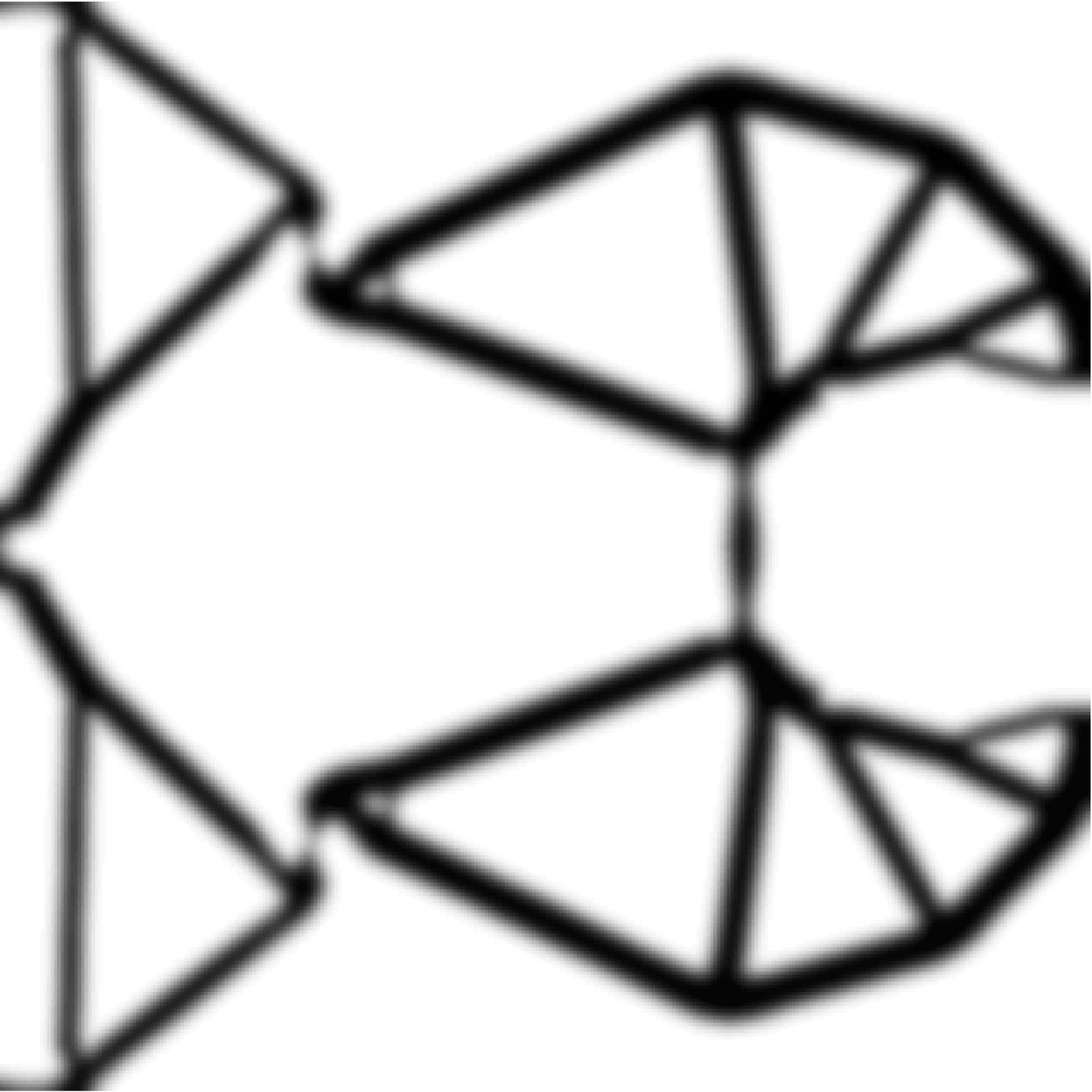} \\[6 pt]
          Small displacements & Large displacements & Small displacements & Large displacements 
      \end{tabular}
      \caption{Optimal topologies for the inverter $300 \times 150$ (left) and the gripper $320 \times 160$ (right).} 
      \label{figure_inverter_gripper}
\end{figure}

Aiming to assess the effect of the geometric nonlinearity for the mechanisms, following the same  ideas applied to the structures, we have also generated the optimal configurations under small displacements for the inverter and the gripper, so they can be compared with the optimal configurations obtained when large displacements are considered. The varied topologies obtained are presented in Fig.~\ref{figure_inverter_gripper}. Since there was no significant difference between the structures obtained with the strategies presented on Table \ref{table_100_iter_inv_grp}, only the results for the {\tt upK100g} strategy are shown.

\begin{figure}[htbp!]
   \centering 
      \begin{tabular}{cc} 
           \includegraphics[width=0.36\textwidth]{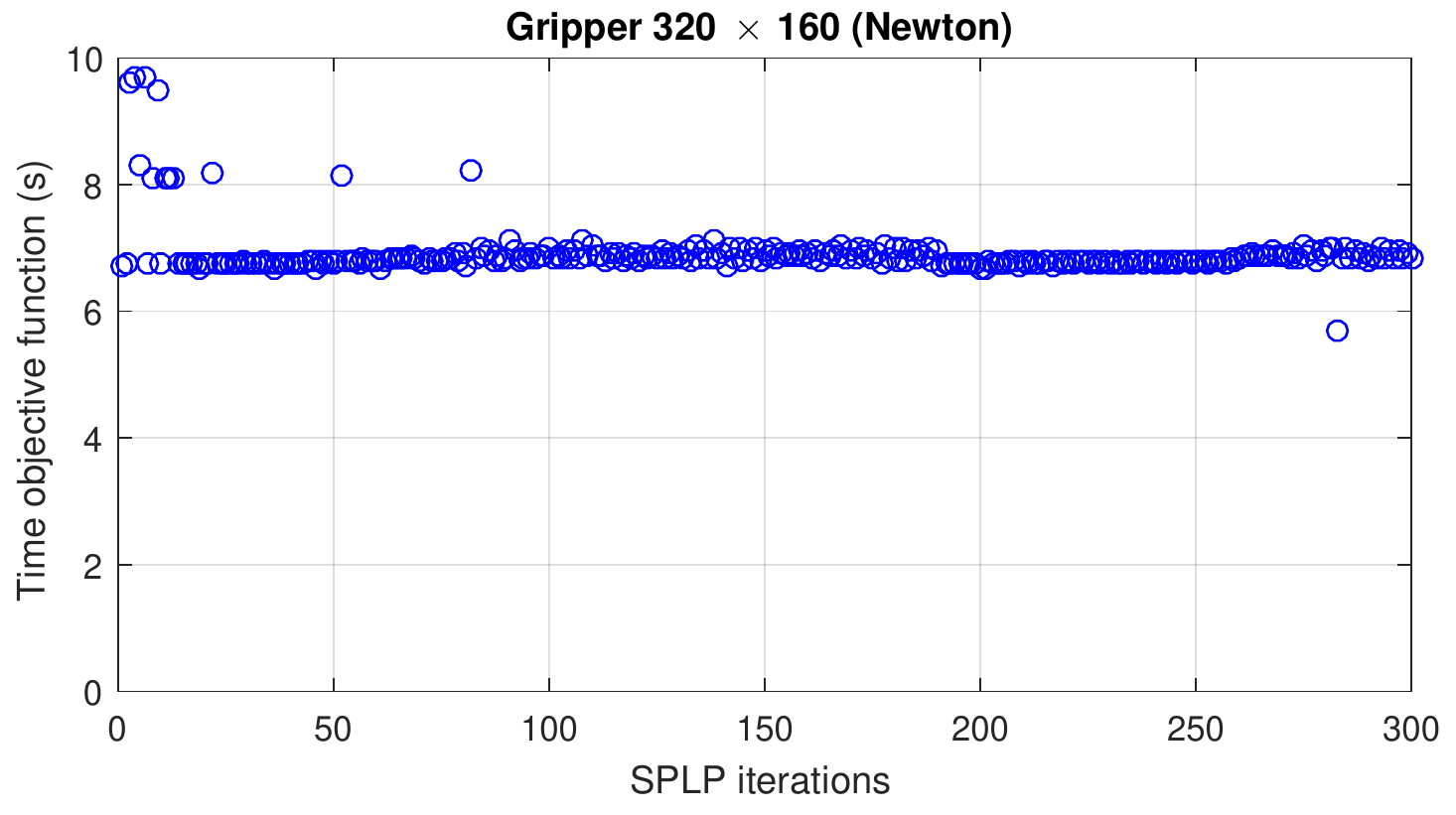}
        & \includegraphics[width=0.36\textwidth]{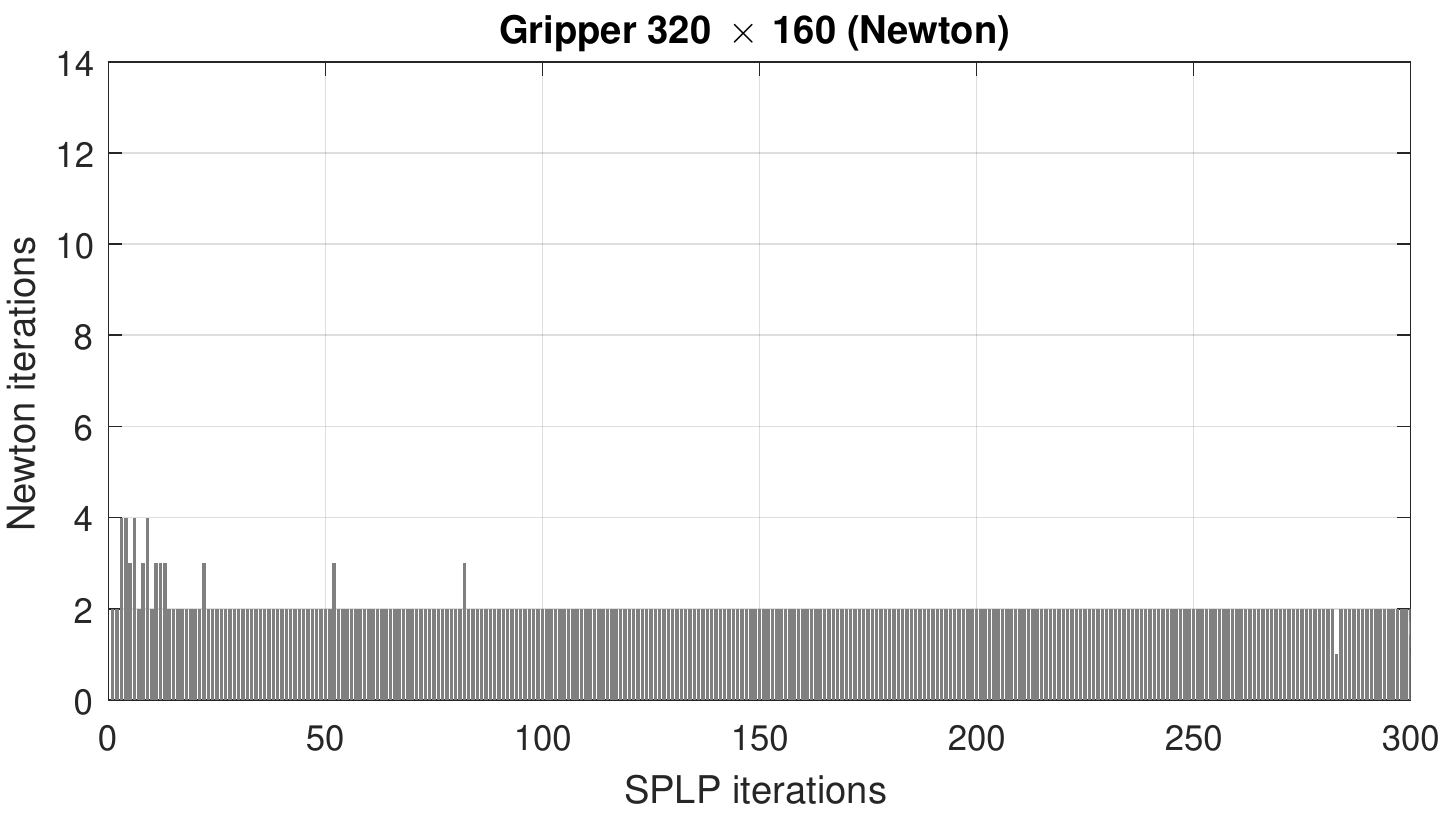} \\
	     \includegraphics[width=0.36\textwidth]{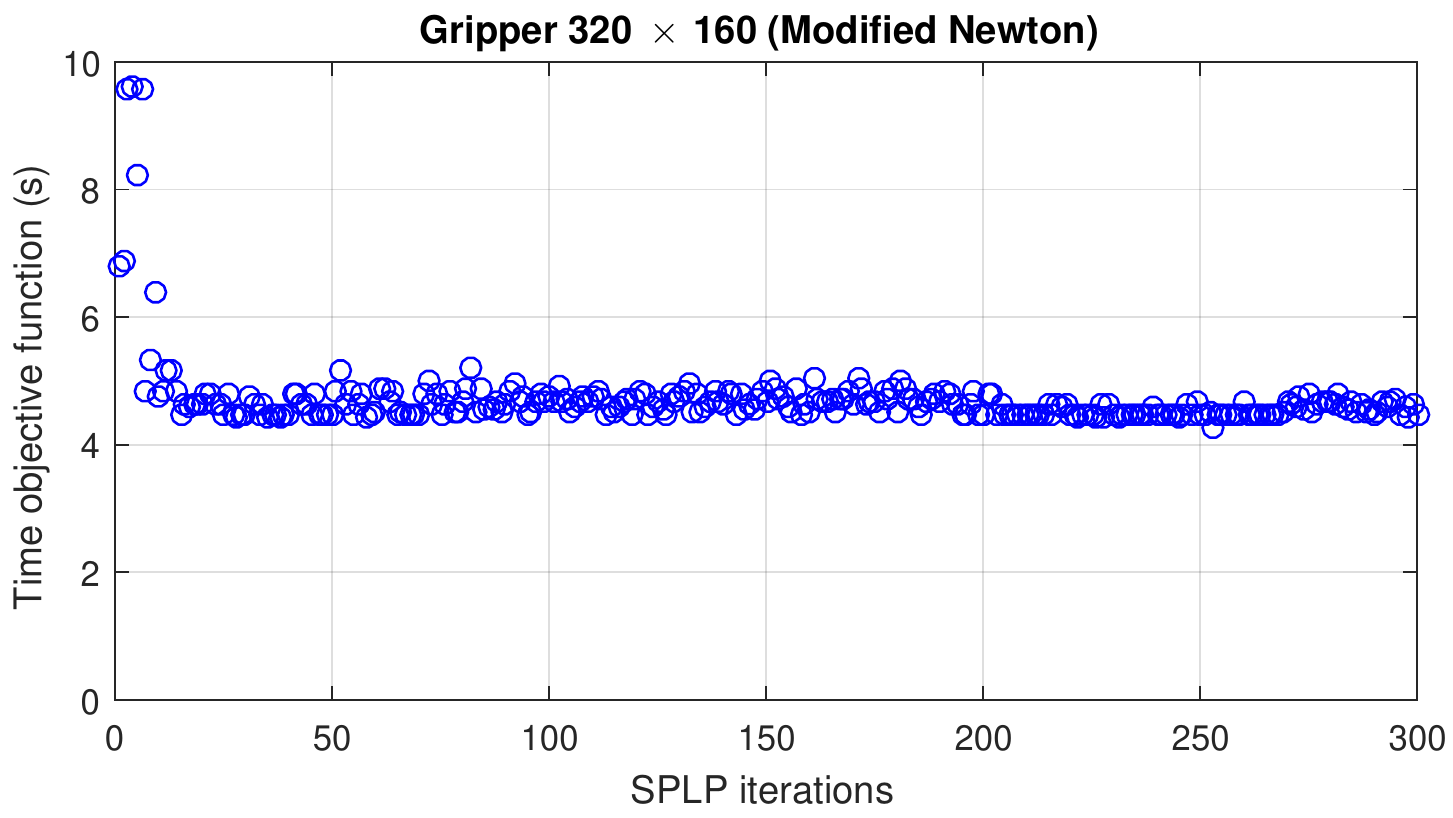}
        & \includegraphics[width=0.36\textwidth]{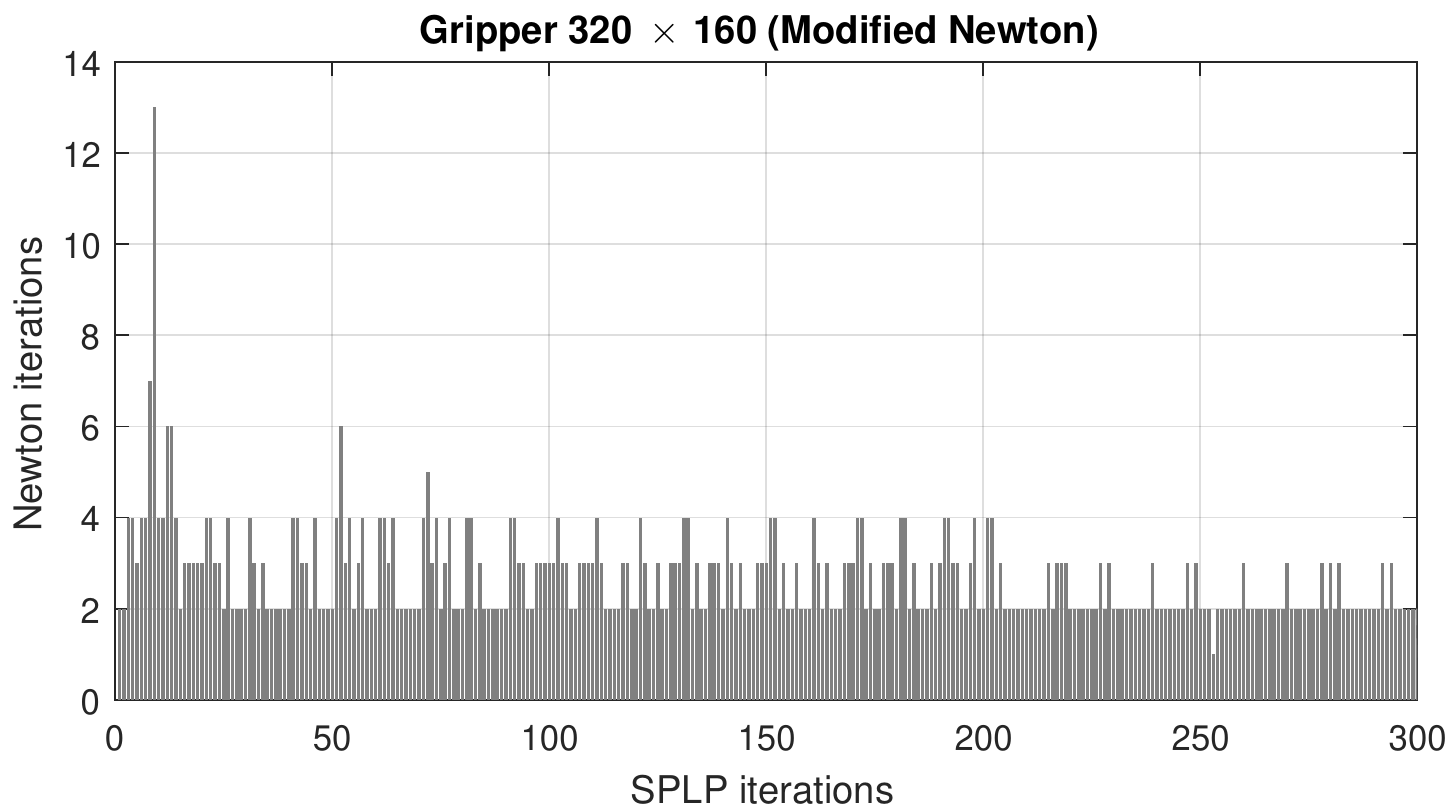} \\
	     \includegraphics[width=0.36\textwidth]{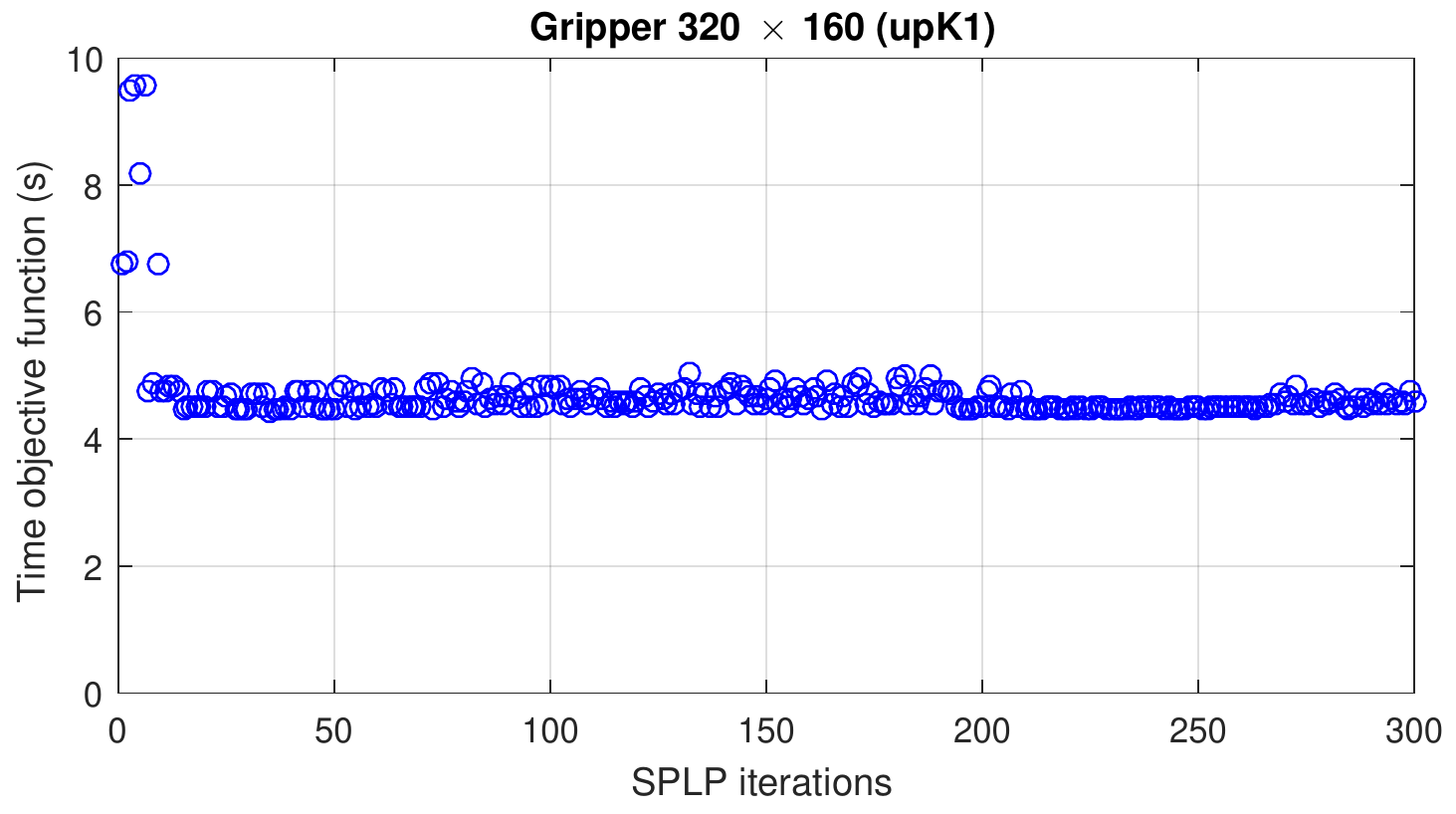}
        & \includegraphics[width=0.36\textwidth]{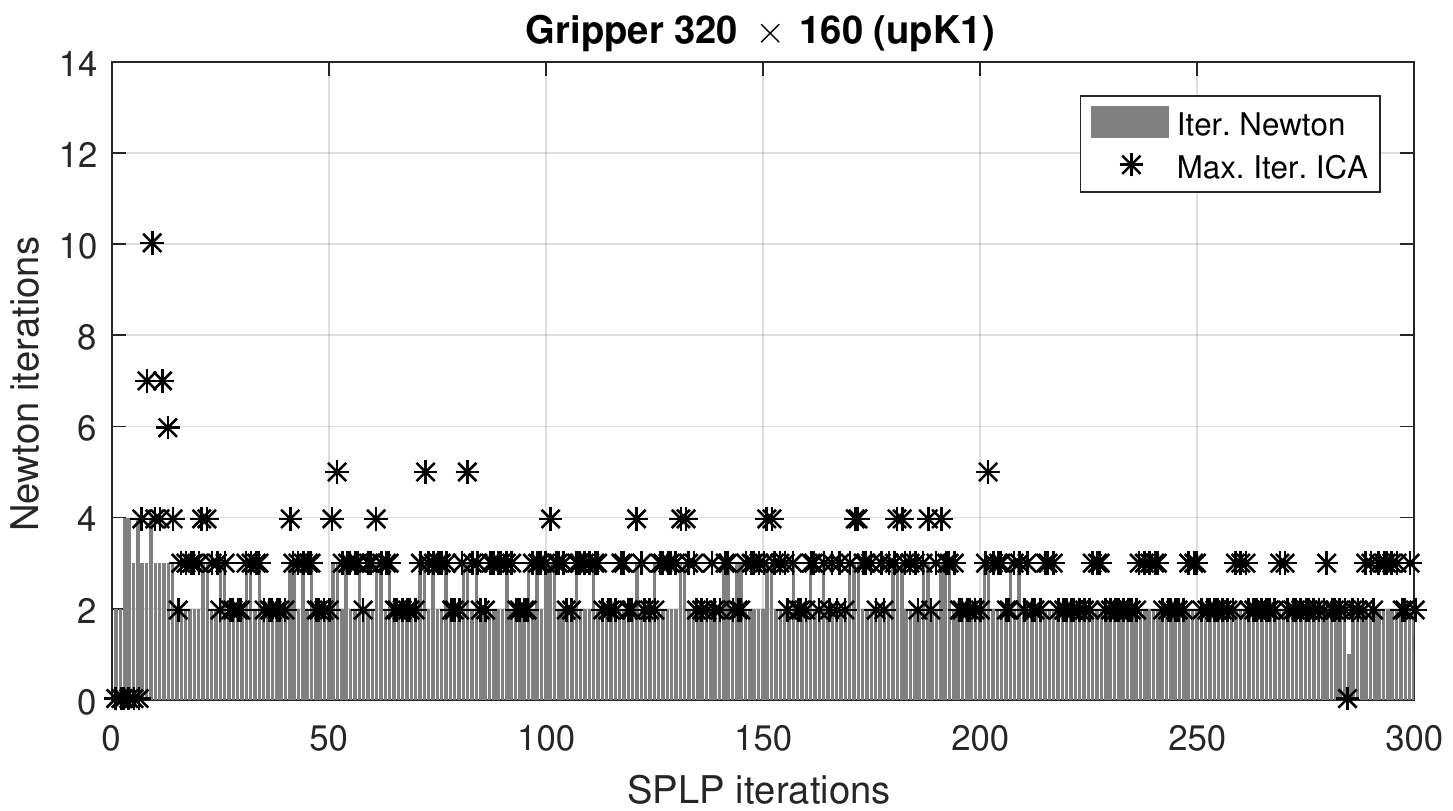} \\
	     \includegraphics[width=0.36\textwidth]{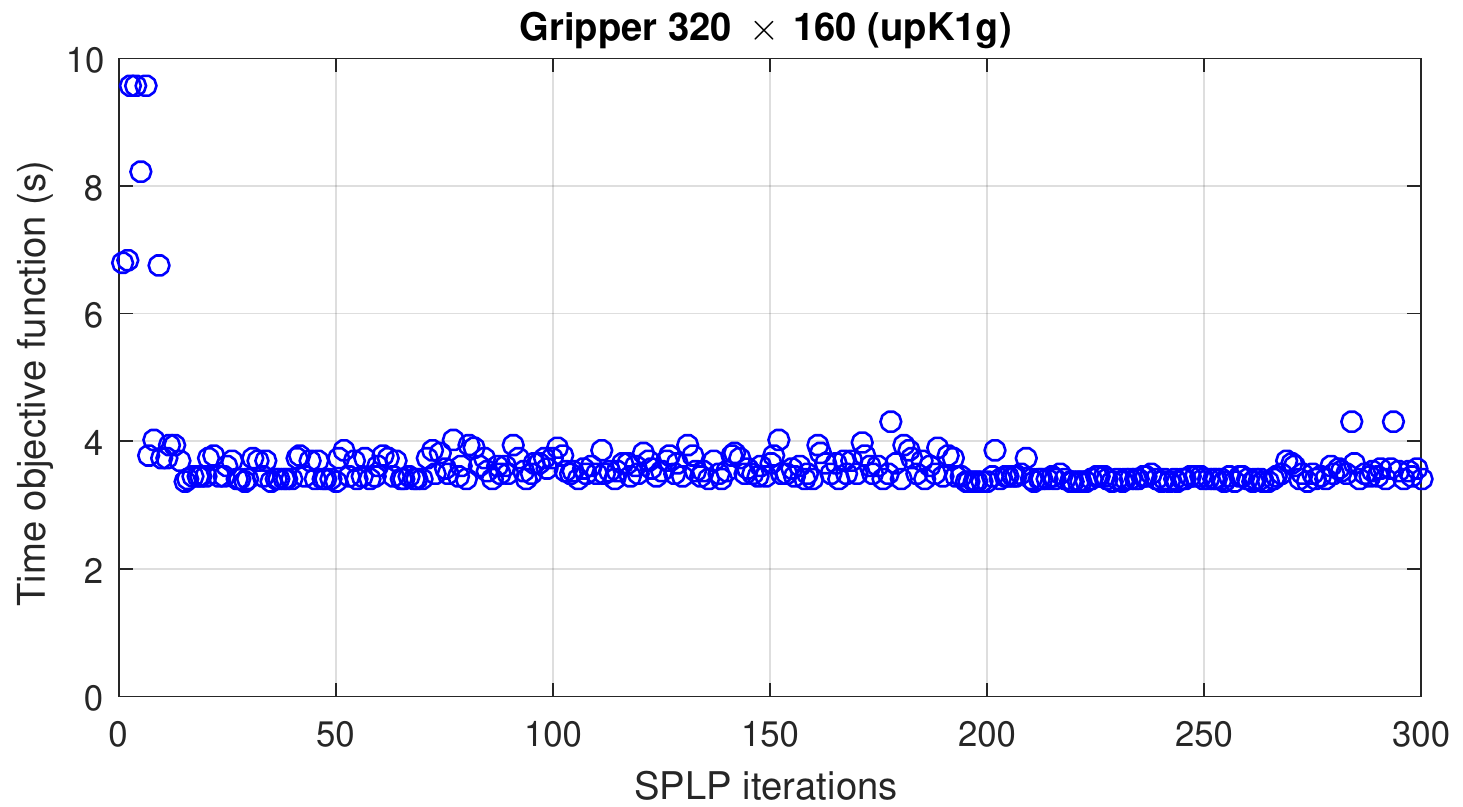}
        & \includegraphics[width=0.36\textwidth]{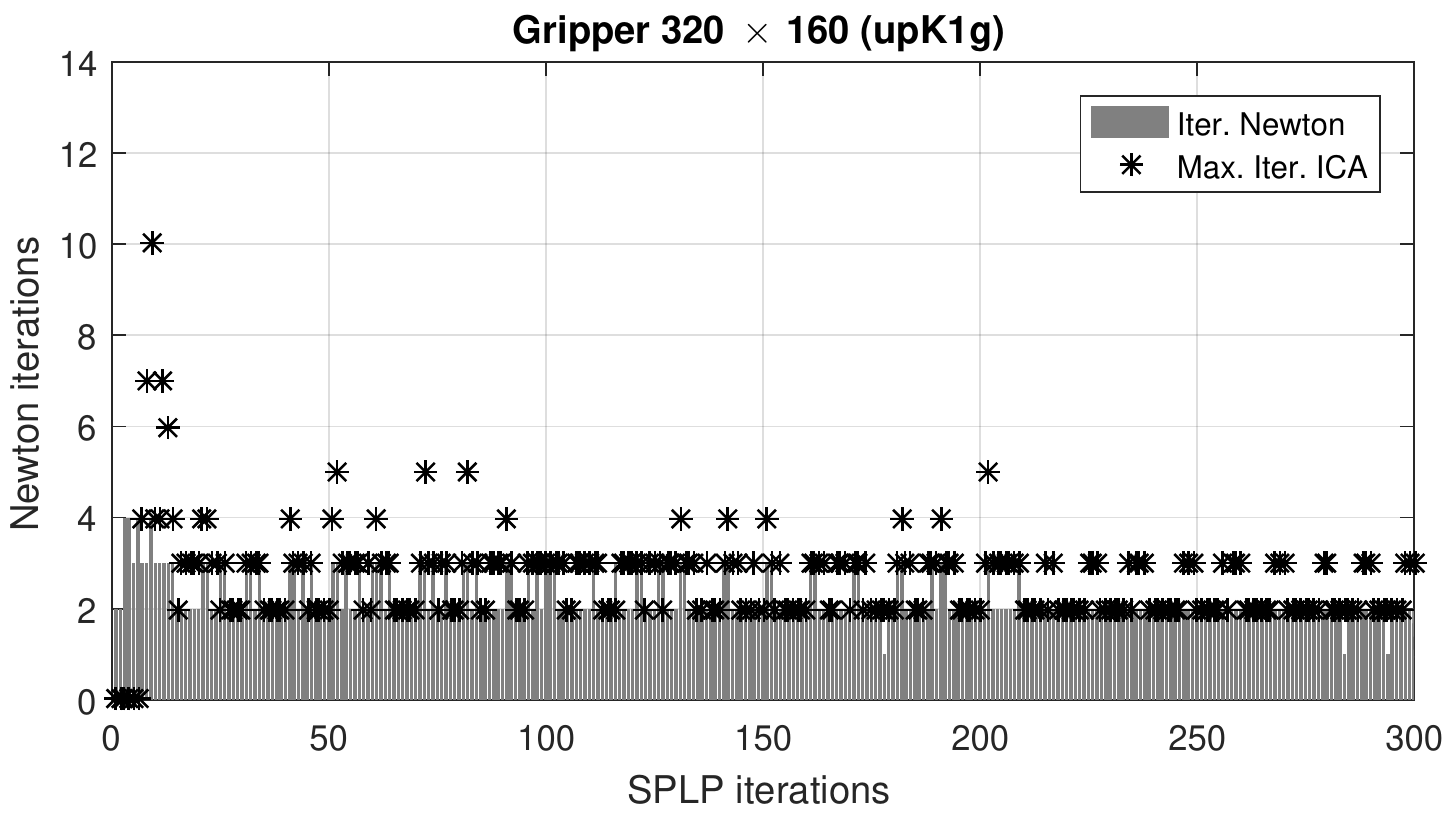} \\
	     \includegraphics[width=0.36\textwidth]{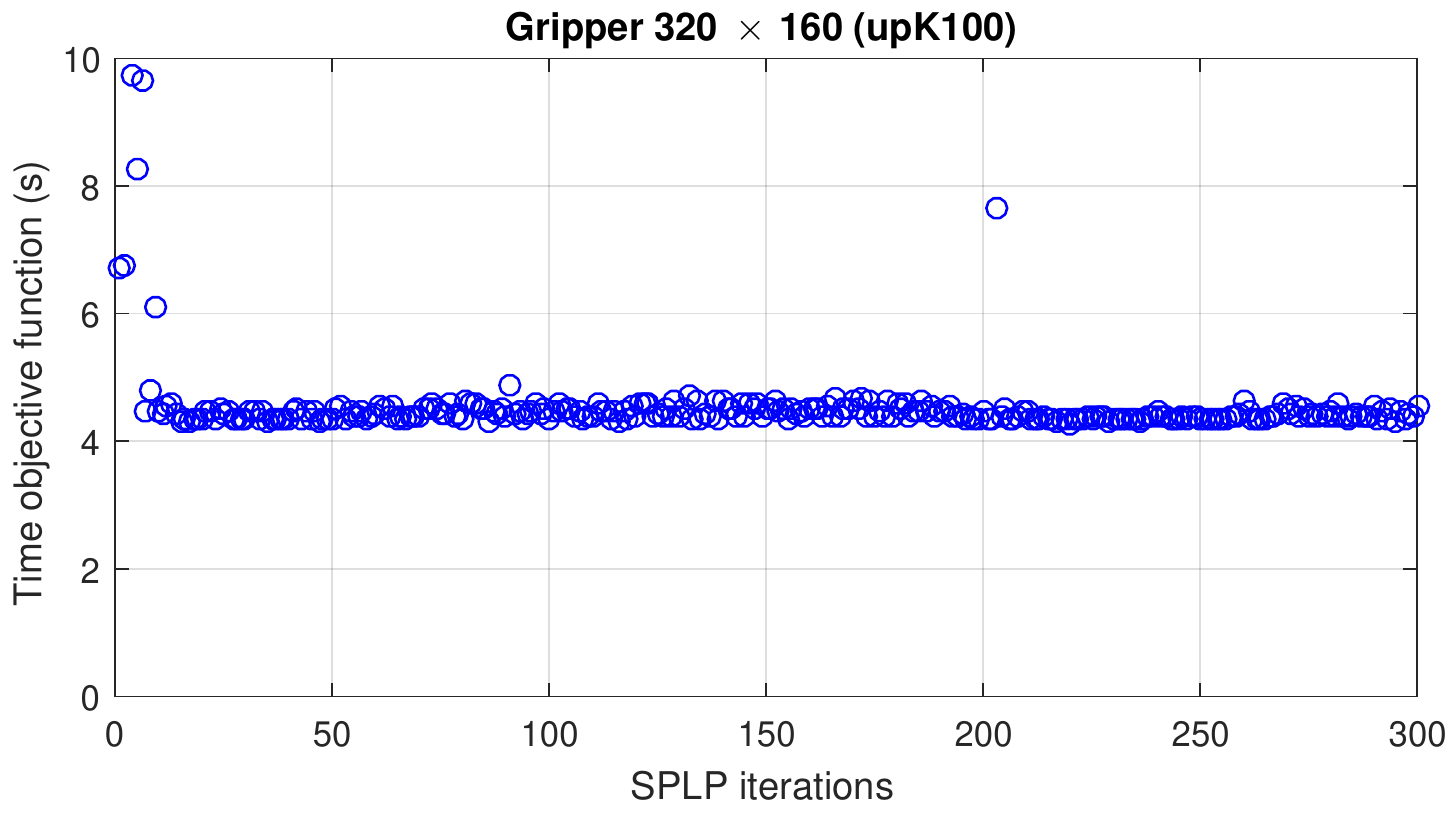}
        & \includegraphics[width=0.36\textwidth]{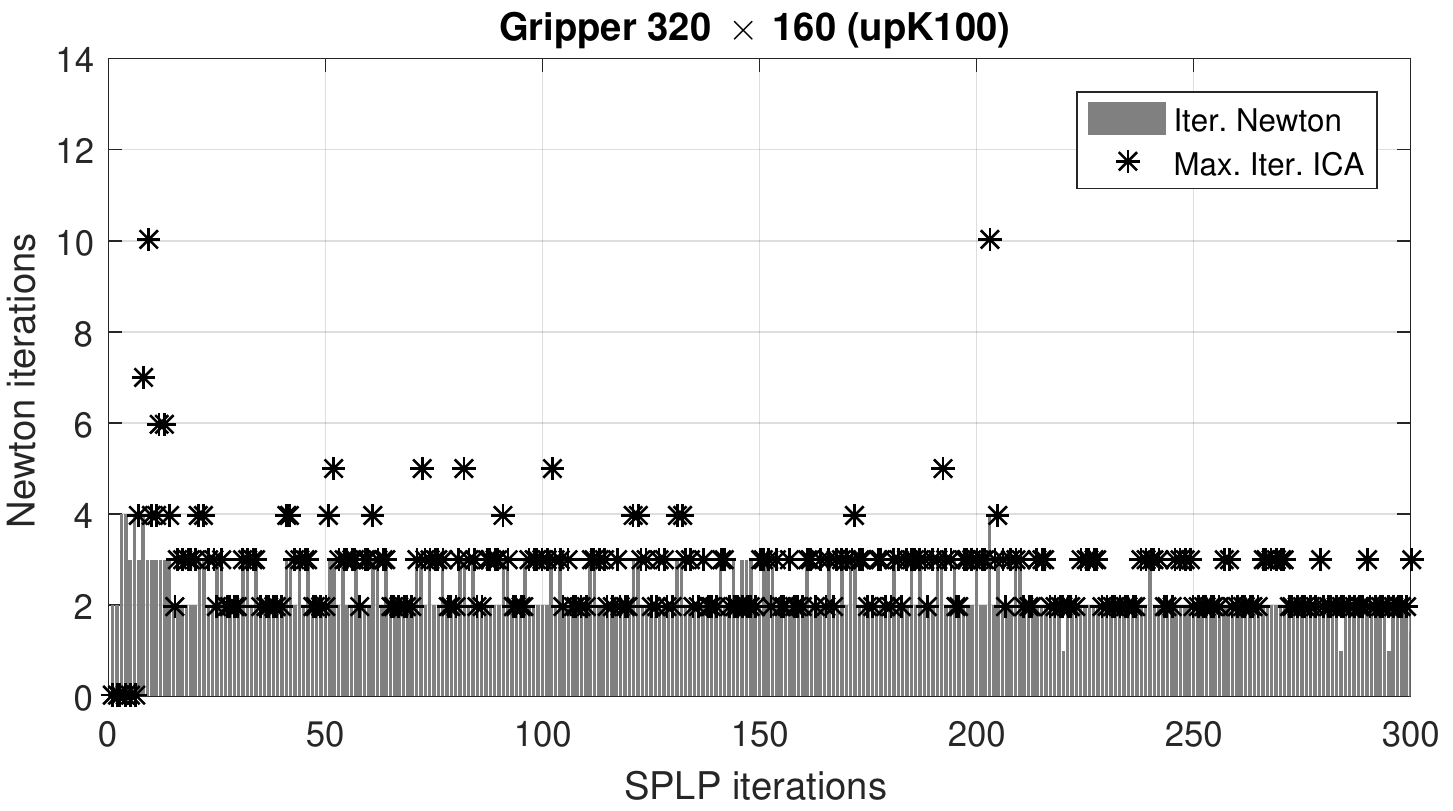} \\
	     \includegraphics[width=0.36\textwidth]{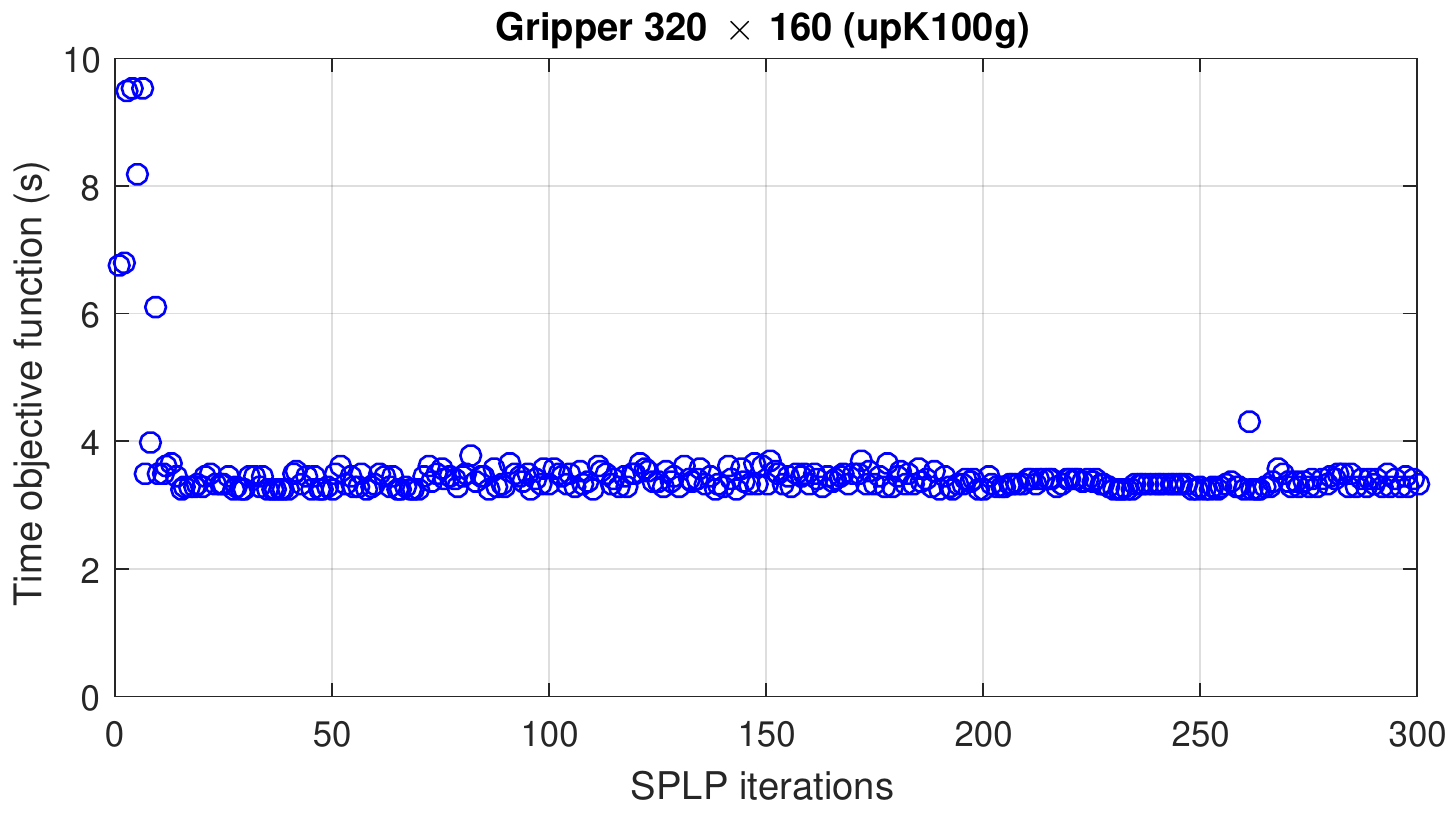}
        & \includegraphics[width=0.36\textwidth]{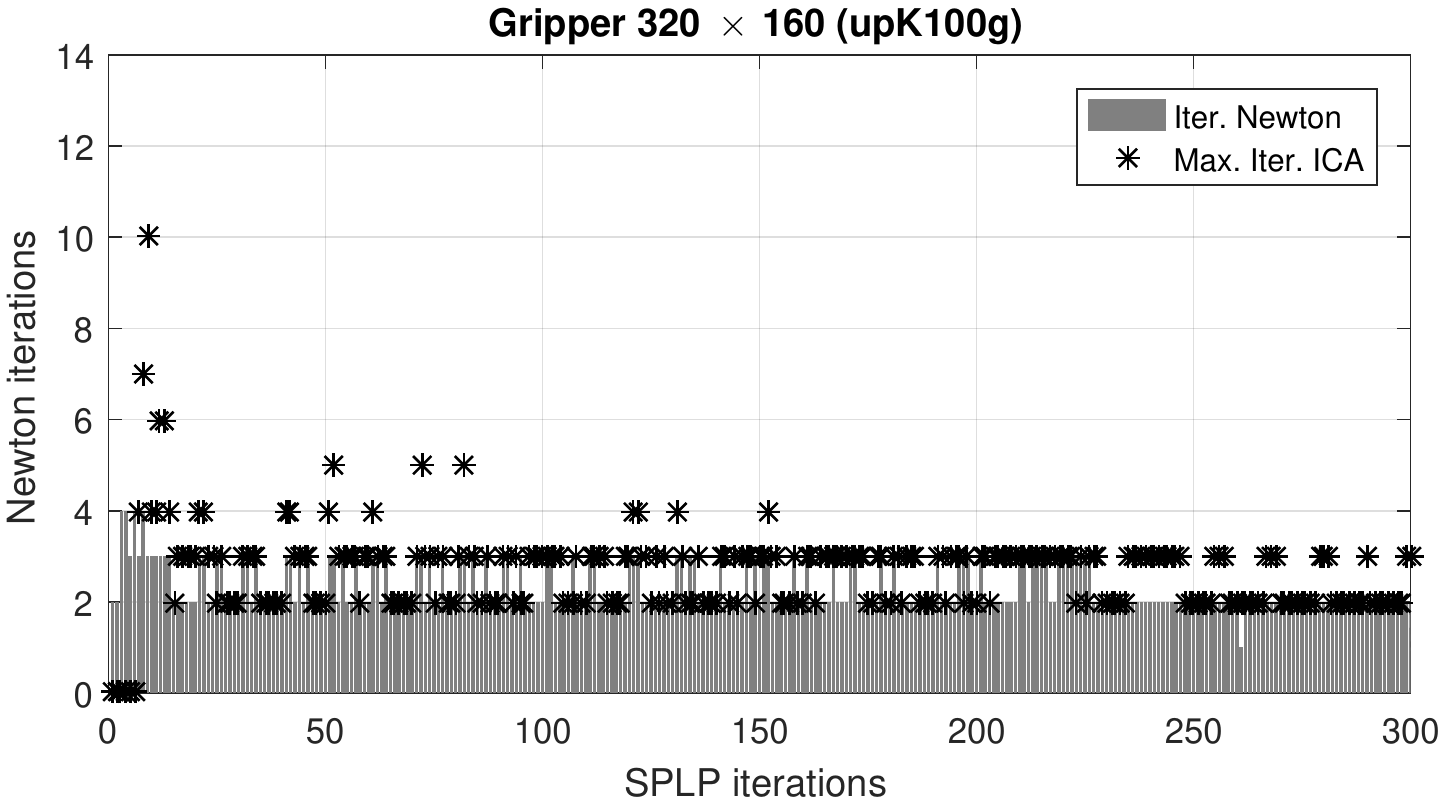} \\
	     \includegraphics[width=0.36\textwidth]{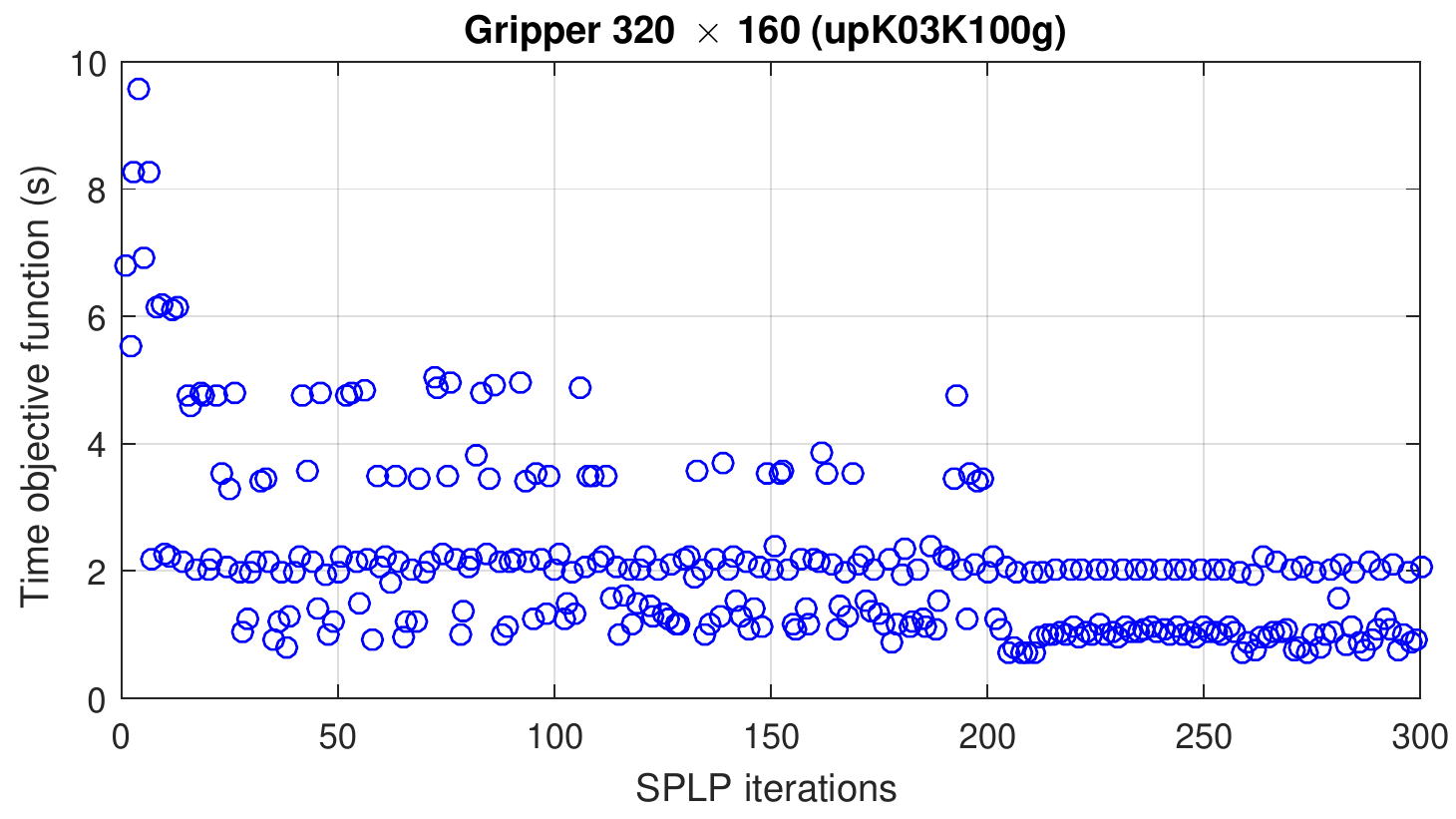}
        & \includegraphics[width=0.36\textwidth]{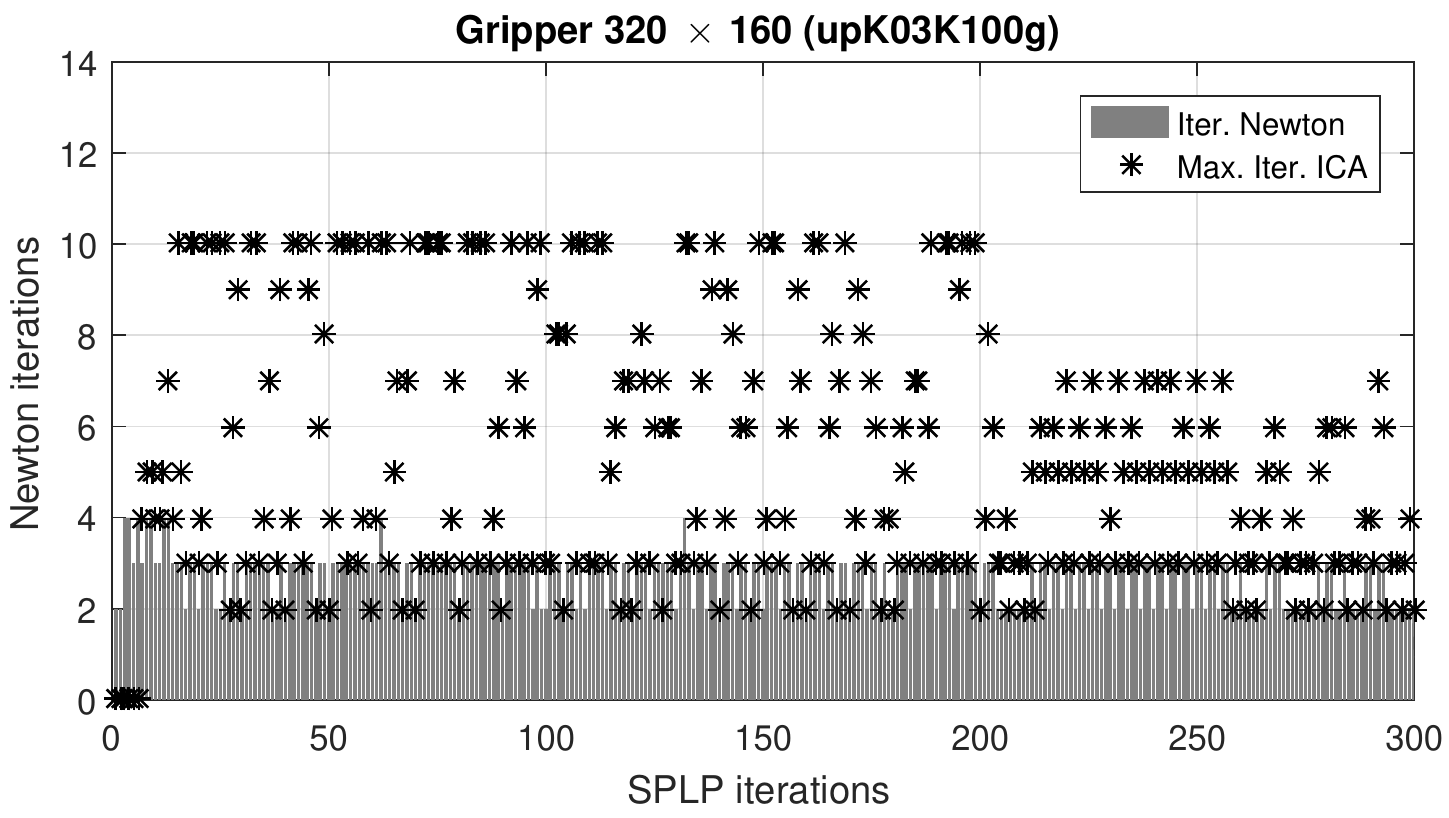} \\
      \end{tabular}
      \caption{Time spent to compute the objective function (left) and number of iterations for Newton's Method (right) in the solution of the gripper ($320 \times 160$).} 
      \label{time_fobj_iter_newton_gripper320x160}
\end{figure}

Figure \ref{time_fobj_iter_newton_gripper320x160} shows the time spent by each strategy to compute $F(\bsrho)$ and the number of iterations performed by Newton's method at each iteration of the SPLP method, along with the maximum number of iterations required to approximately solve the linear system (\ref{eq02}) for the gripper mechanism. As we can see, the convergence of Newton's method is very fast. In fact, in almost all of the iterations of the SPLP method, the approximate solution of the nonlinear system (\ref{eq01b}) was found in just 2 iterations. On the other hand, the average time spent to compute $F(\bsrho)$ is equal to $7.0$ seconds. With the adoption of the Modified Newton method, we can observe that, at the tenth iteration of the SPLP method, 13 iterations were necessary to find the approximate solution of the nonlinear system (\ref{eq01b}), but, in general, $F(\bsrho)$ is computed in no more than $5.5$ seconds. 

The overall performance of the Modified Newton method and of strategies {\tt upK1} and {\tt upK100} were very similar, but the last two strategies required less iterations for solving the nonlinear system (\ref{eq01b}), due to the frequent updating of the $\DeltaK$ matrix. Moreover,  in most cases, the approximate solution of the linear system (\ref{eq02}) was found in no more than 3 iterations of the scheme described in~(\ref{eq03}). Despite of similarities between strategies {\tt upK1} and {\tt upK100}, the latter spent $22.3\,\%$ less time to assemble $\KK$ than the former. 

Strategies {\tt upK1g} and {\tt upK100g} required much less time to factorize $\KK$ than its counterparts {\tt upK1} and {\tt upK100}. In fact, for the two compliant mechanisms, {\tt upK1g} spent $31.3\,\%$ less time than {\tt upK1}. The same amount of reduction was attained by the {\tt upK100g} strategy in comparison with {\tt upK100}. We also observe in Fig. \ref{time_fobj_iter_newton_gripper320x160} that the convergence of the inexact Newton variants have similar overall behavior.

At the iterations of the SPLP method in which the factorization of $\KK$ was reused,  the upK03K100g strategy took no more than $2$ seconds to compute $F(\bsrho)$. In some occasions, the time spent was less than $1$ second. Because of the smaller frequency of factorizations of $\KK$, the average number of iterations to obtain the solution of the linear system~(\ref{eq02}) was larger than those required by the other inexact Newton strategies. Despite this fact, the accuracy of the approximate solutions of the linear system~(\ref{eq02}) and of the nonlinear system~(\ref{eq01b}) were not spoiled. For the gripper mechanism, this strategy took $70.3\,\%$ less time than Newton's method to compute $F(\bsrho)$, and the reduction of the time consumed in the factorization of $\KK$ reached $81.4\,\%$.

%===========================================
\section{Effect of mesh refinement}\label{sec:meshref}

The time consumed to evaluate the objective function (and, consequently, the total time spent for solving the topology optimization problem) is directly influenced by the mesh refinement of the design domain. To investigate the impact of the number of finite elements on the time spent to compute $F(\bsrho)$, we chose four different meshes for the slender beam and for the inverter, considering a fixed budget of 300 iterations of the SPLP method. For the slender beam, the design domain was discretized into $200 \times 25 = 5000$, $400 \times 50 = 20000$, $600 \times 75 = 45000$ and $800 \times 100 = 80000$ elements. In each case, the filter radius was defined accordingly, corresponding to the length of $2.5$, $5.0$, $7.5$ and $10.0$ elements, respectively. The upper half of the design domain of the inverter was discretized into $200 \times 100 = 20000$, $300 \times 150 = 45000$, $400 \times 200 = 80000$ and $500 \times 250 = 125000$ elements, and the filter radius was set to the length of $5.0$, $7.5$, $10.0$ and $12.5$ elements, respectively. The tables that show how the performance of each strategy was affected by the mesh refinement are available in Appendix \ref{sec:tables_300it}.

\subsection{Slender beam}

\begin{figure}[htbp!]
   \centering
 \begin{tabular}{cc} 
         \includegraphics[width=0.45\textwidth]{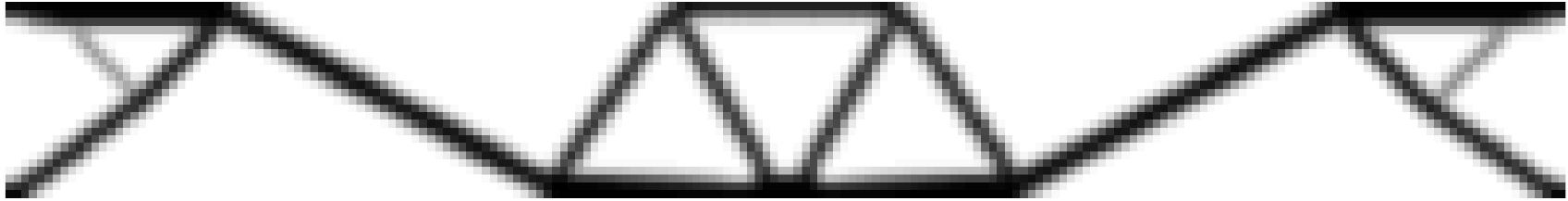} 
       & \includegraphics[width=0.45\textwidth]{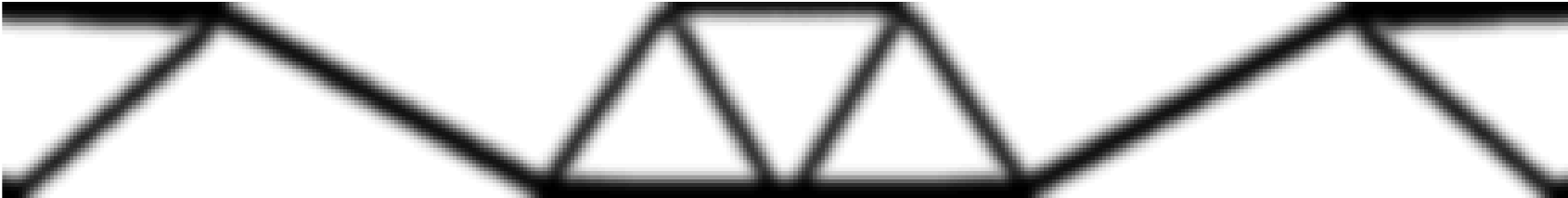} \\[3pt]
         $200 \times 25$ & $400 \times 50$ \\[10pt]
         \includegraphics[width=0.45\textwidth]{sb600x75_upK100g_crop.pdf} 
       & \includegraphics[width=0.45\textwidth]{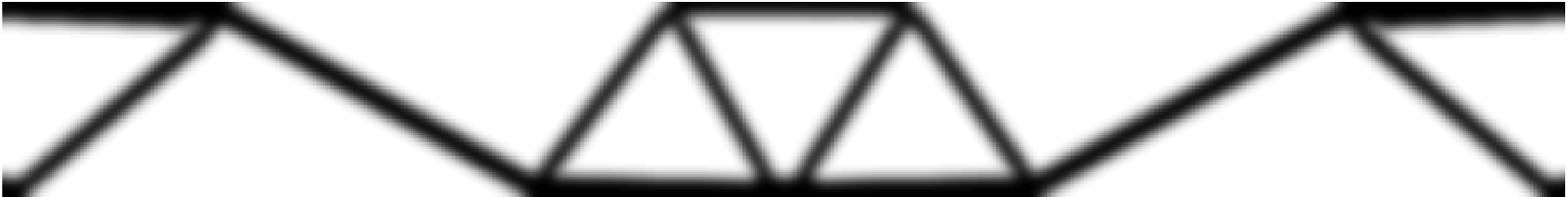} \\[3pt] 
         $600 \times 75$ & $800 \times 100$
  \end{tabular}
      \caption{Optimal topologies for the slender beam.} 
      \label{figure_sb2}
\end{figure}

First of all, it is worth mentioning that, for all of the meshes considered, the value of $F(\bsrho)$ attained after 300 iterations of the SPLP method differ very little regardless of the strategy adopted, corroborating the robustness of the reanalysis strategies proposed here. Since the optimal topologies are also very similar among strategies, only the results related to the {\tt upK100g} strategy are shown in Fig. \ref{figure_sb2}.

\begin{figure}[htbp!]
   \centering    
      \includegraphics[width=0.85\textwidth]{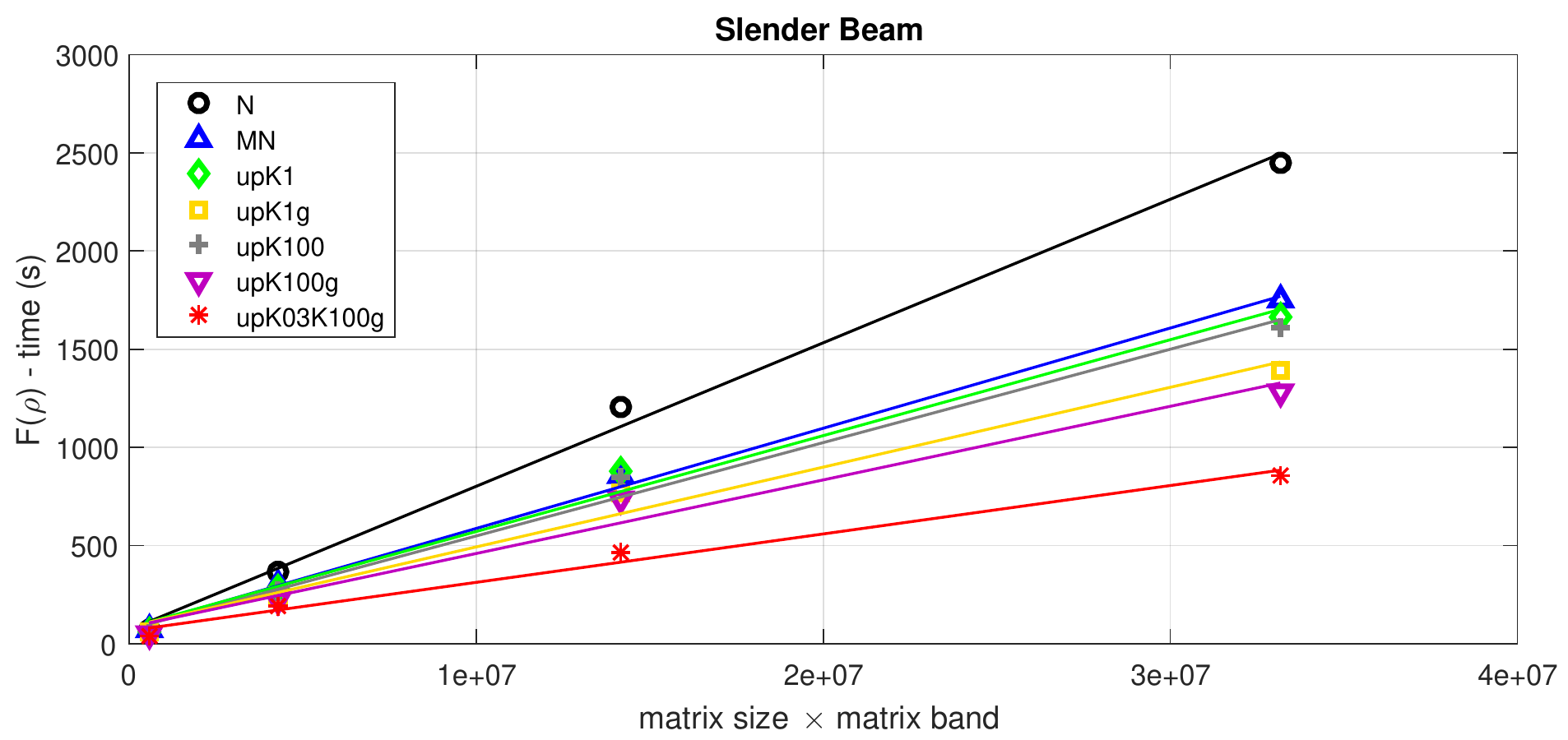}
      \caption{$F(\bsrho)$-time for each strategy according to the dimension of $\KK$ (slender beam).} 
      \label{time_fobj_discretization_slender_beam}
\end{figure}

Figure \ref{time_fobj_discretization_slender_beam}, however, shows that the time required to evaluate the objective function may vary significantly among strategies. This figure suggest that, for all of the methods, there is an almost linear relation between the time spent to compute $F(\bsrho)$ and the product of the size of matrix $\KK$ and its bandwidth. Moreover, all of our strategies not only provide a considerable decrease on the overall time spent to compute $F(\bsrho)$\, but also this reduction become more evident as we refine the design domain. In fact, the slopes of the regression lines of the Modified Newton method and strategies {\tt upK1}, {\tt upK1g}, {\tt upK100}, {\tt upK100g} and {\tt upK03K100g} are, respectively, $30.2\,\%$, $33.2\,\%$, $44.4\,\%$, $35.0\,\%$, $48.7\,\%$ and $66.3\,\%$ smaller than the slope of the line associated with Newton's method.

\begin{figure}[htbp!]
   \centering    
       \begin{tabular}{cc} 
          \includegraphics[width=0.45\textwidth]{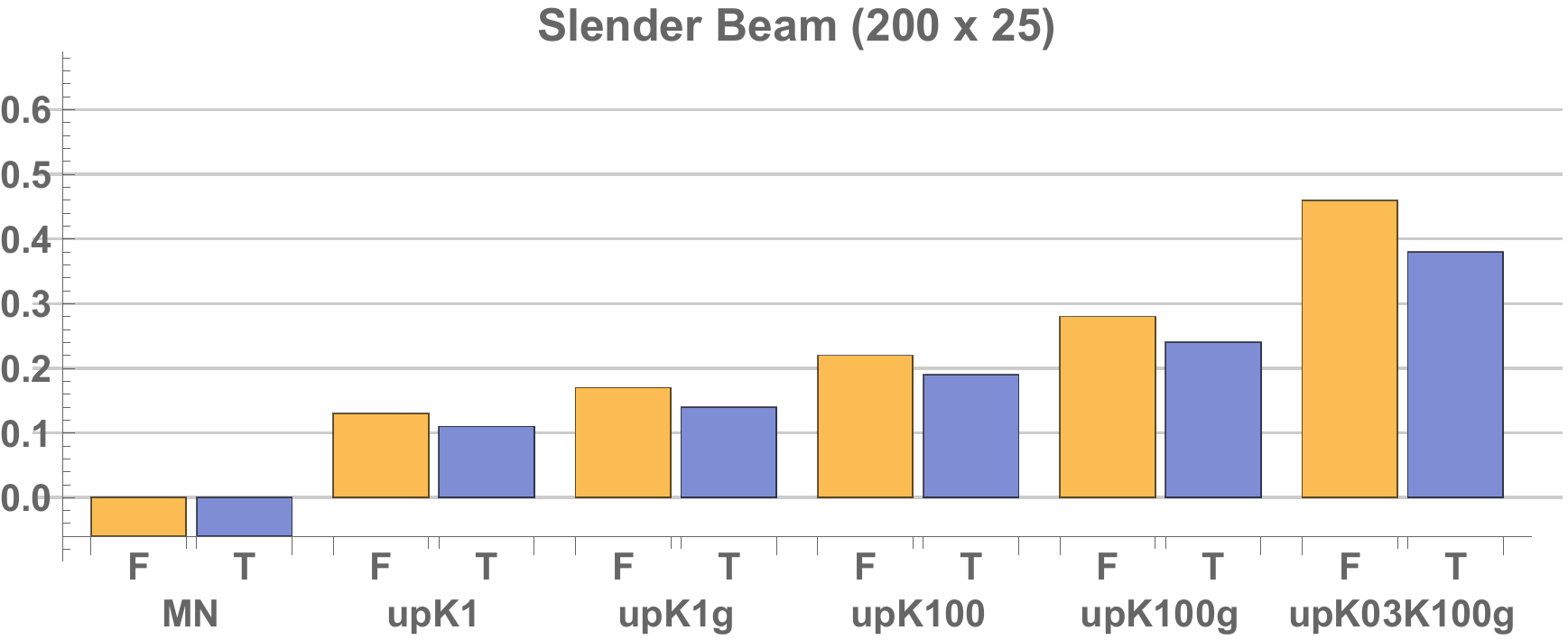}
        & \hspace{0.3cm}\includegraphics[width=0.45\textwidth]{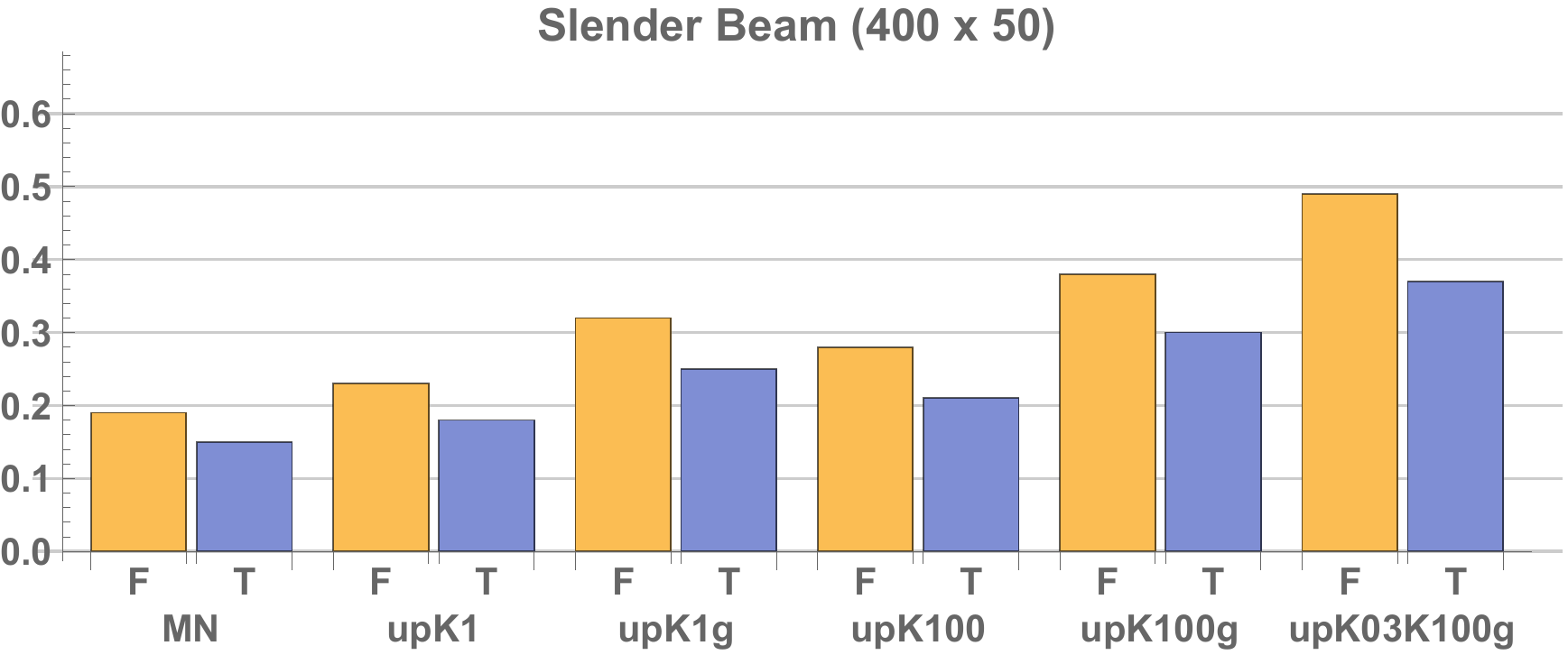} \\[10pt]
          \includegraphics[width=0.45\textwidth]{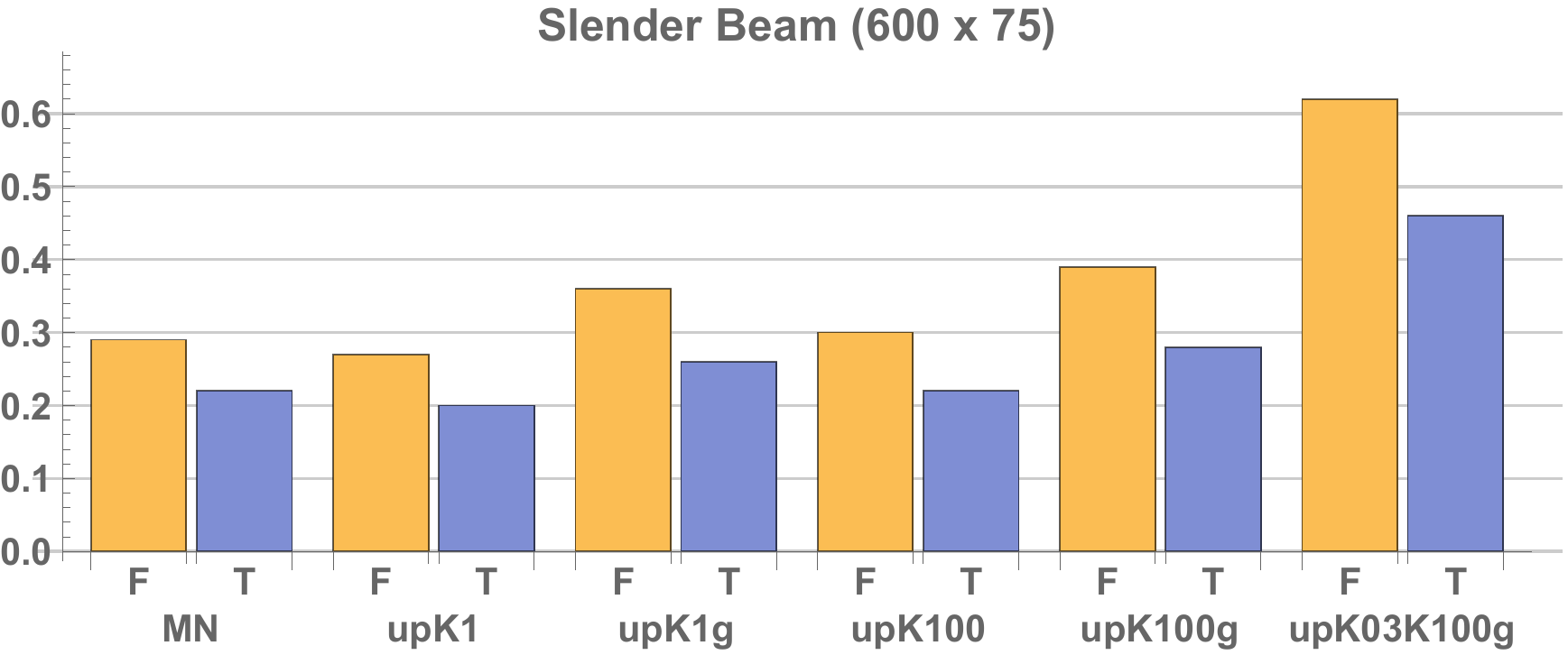}
        & \hspace{0.3cm}\includegraphics[width=0.45\textwidth]{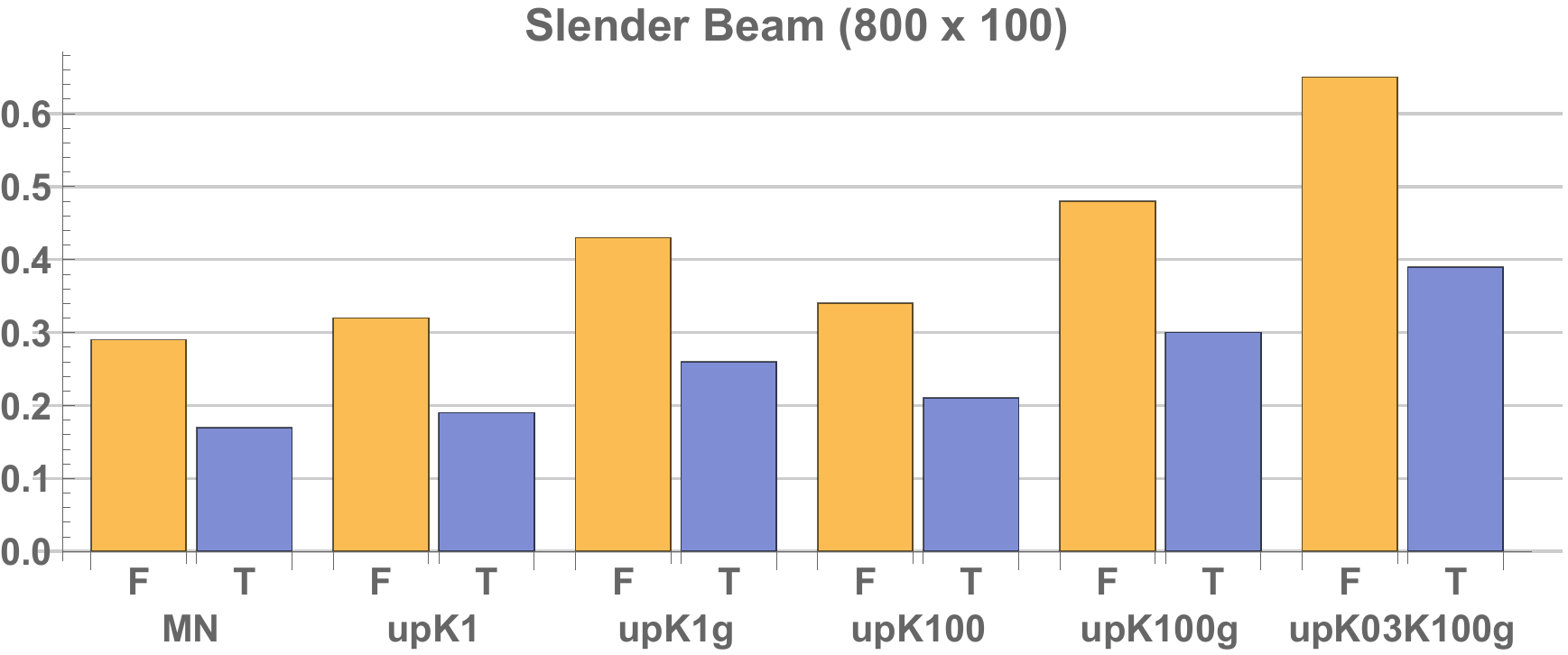} \\ 
      \end{tabular}
      \caption{Percentage of decrease in the time demanded to evaluate the objective function (\textsf{F}) and in the total time (\textsf{T}), as compared with Newton's method, for the slender beam.}
      \label{percentage_reduction_fobj_sb}
\end{figure}

Figure \ref{percentage_reduction_fobj_sb} shows the percentage of decrease in the time spent to evaluate $F(\bsrho)$ with respect to Newton's method, for the four discretizations considered for the slender beam. For the coarsest mesh, the Modified Newton method took almost twice the number of iterations of Newton's method, resulting in an increase of $5.7\,\%$ in the time spent to compute the objective function. For the remaining discretizations, the performance of the Modified Newton method and strategy {\tt upK1} were very similar, although the latter took fewer iterations due to its frequent updating of $\DeltaK$. It is clear from the figure that the inexact solution of the system related to the objective function gradient becomes more relevant as the mesh gets finer, so the decrease in time is more pronounced for the {\tt upK1g}, {\tt upK100g} and {\tt upK03K100g} strategies.

Although strategies {\tt upK1g} and {\tt upK100g} attained an average reduction of $34.2\,\%$\, and $38.5\,\%$, respectively, for the last two discretizations, {\tt upK03K100g} showed the best figures for all of the meshes, as expected. For this strategy, the reduction of the time consumed in evaluating $F(\bsrho)$ varied from $46.0\,\%$ to $65.0\,\%$. When it comes to  the time for assembling $\KK$, the reduction ranged between $38.4\,\%$ and $45.0\,\%$, while the factorization of $\KK$ showed an even more significant result, with a reduction in the range $71.6\,\%$ to $77.4\,\%$.

Figure \ref{percentage_reduction_fobj_sb} also shows the percentage of decrease in the total time spent by the inexact Newton strategies, with respect to Newton's method. The difference between the two percentages presented for each strategy reveals the relevance of the remained steps performed at each iteration of the SPLP method. In short, for the Modified Newton method and strategies {\tt upK1} and {\tt upK100}, the decrease in the time spent to obtain the optimal topologies was, in average, $18.1\,\%$, $19.3\,\%$ and $21.4\,\%$, respectively. When strategies {\tt upK1g} and {\tt upK100g} were adopted, an average decrease of $25.6\,\%$ and $29.5\,\%$ was obtained. Finally, the {\tt upK03K100g} strategy showed an average decrease of $40.0\,\%$.

\subsection{Inverter}

\begin{figure}[htbp!]
   \hspace*{-0.2cm}
      \begin{tabular}{cccc} 
          \includegraphics[width=0.22\textwidth]{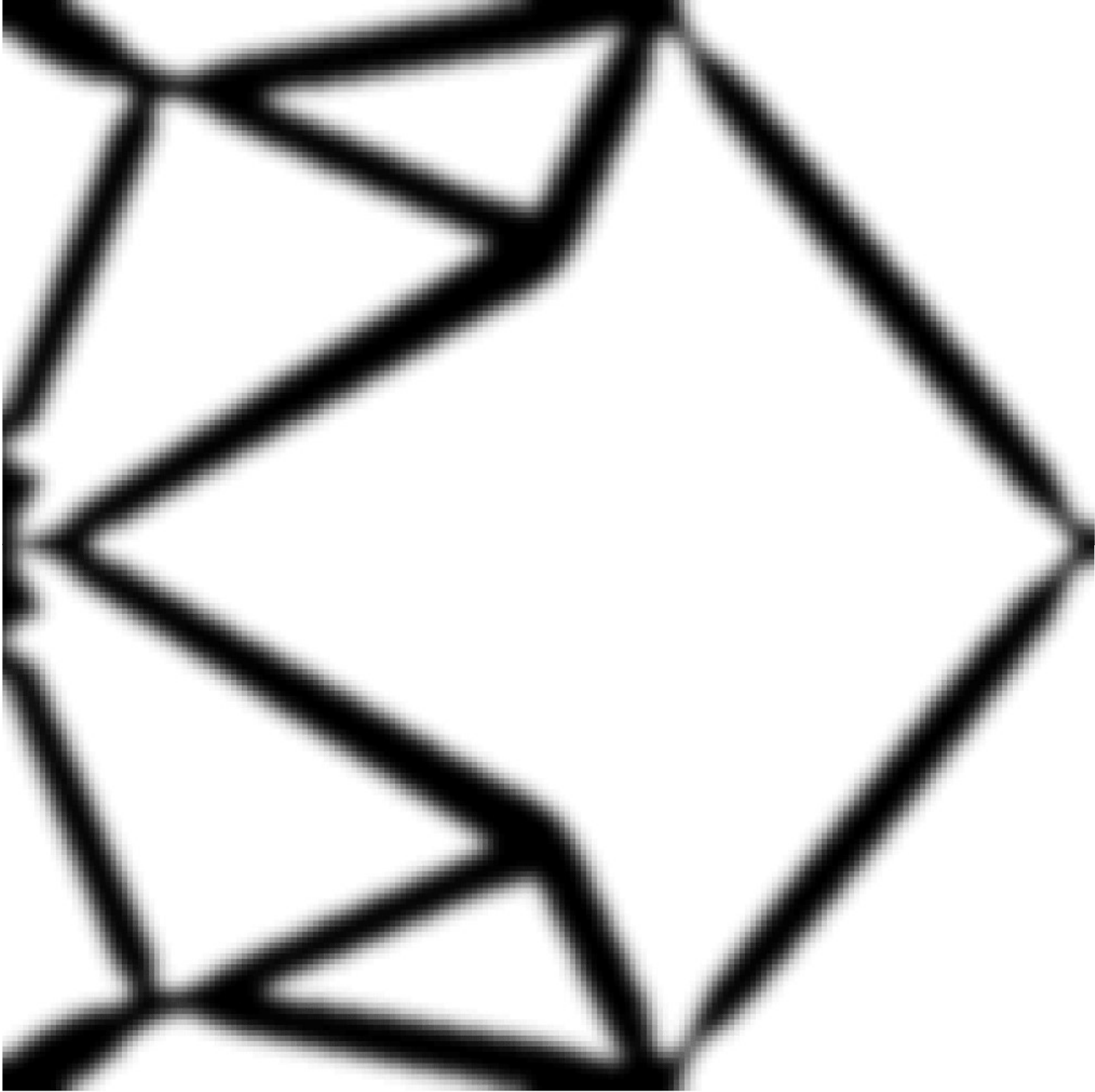} 
        & \includegraphics[width=0.22\textwidth]{inverter300x150_largedisp_top_bottom.pdf}
        & \includegraphics[width=0.22\textwidth]{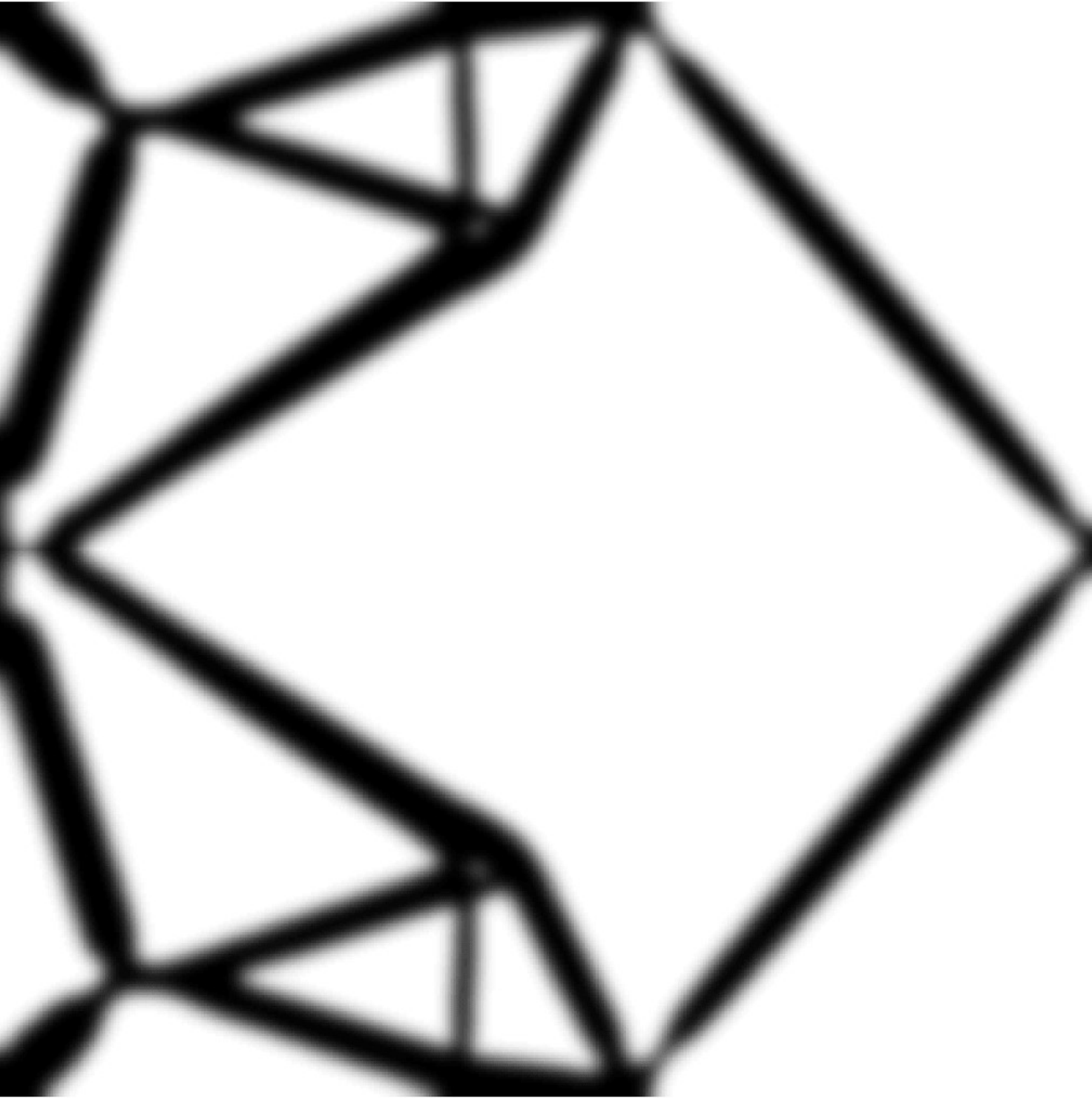} 
        & \includegraphics[width=0.22\textwidth]{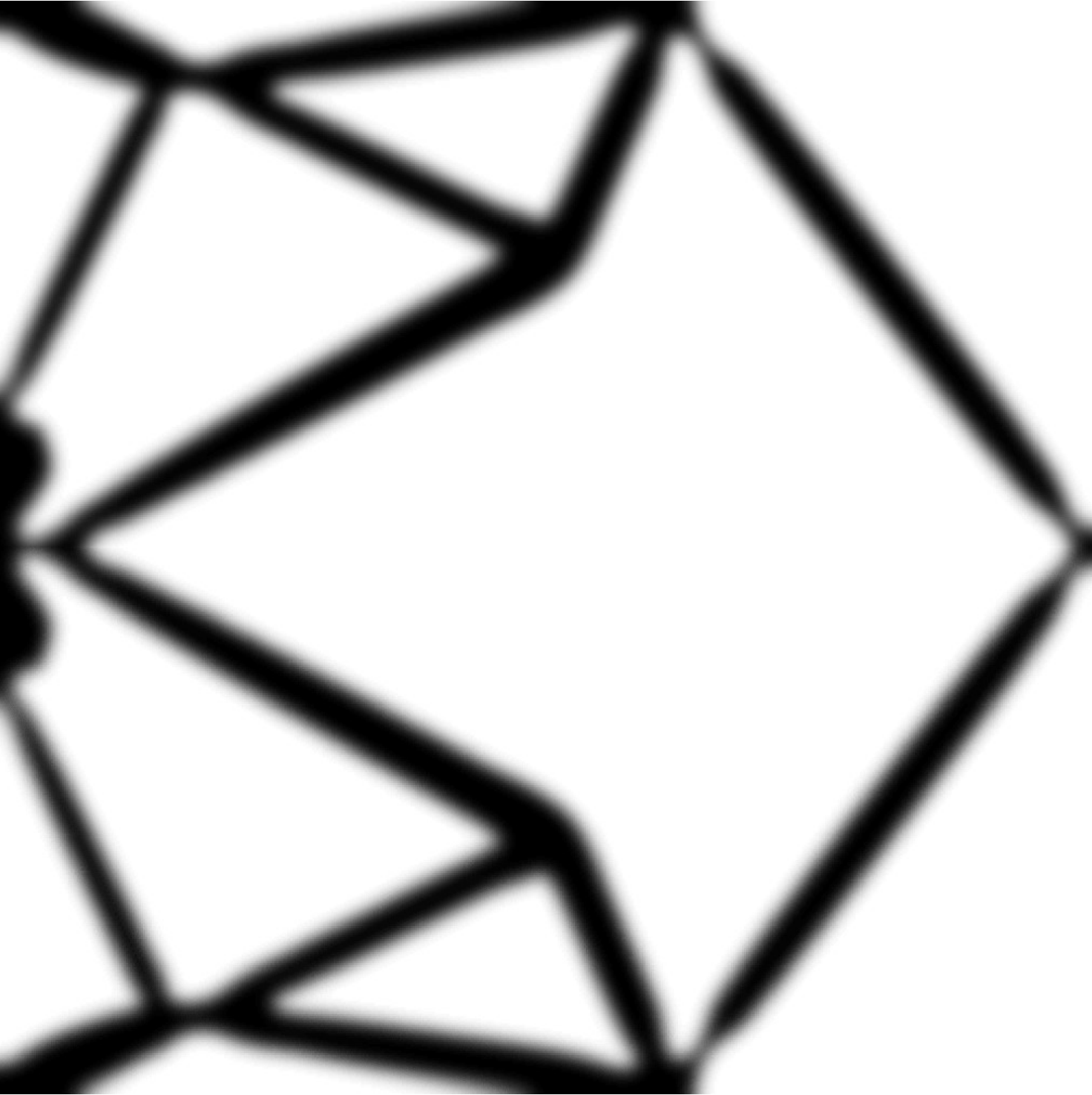} \\[3pt]
          $200 \times 100$ & $300 \times 150$ & $400 \times 200$ & $500 \times 250$
      \end{tabular}
      \caption{Optimal topologies for the inverter.} 
      \label{figure_inverter}
\end{figure}

As we observed for the slender beam, for all of the mesh sizes, the value of the objective function attained by the various strategies after 300 iterations of the SPLP method was almost the same. In this case, the difference did not exceed $0.07\,\%$, which is very reasonable from the practical point of view. Moreover, the number of iterations necessary to obtain the approximate solution of the nonlinear system (\ref{eq01b}) also showed little variation among strategies. As a consequence, the optimal topologies obtained were quite similar, so only the results for the {\tt upK100g} strategy are shown in Fig. \ref{figure_inverter}.

\begin{figure}[htbp!]
   \centering    
       \includegraphics[width=0.85\textwidth]{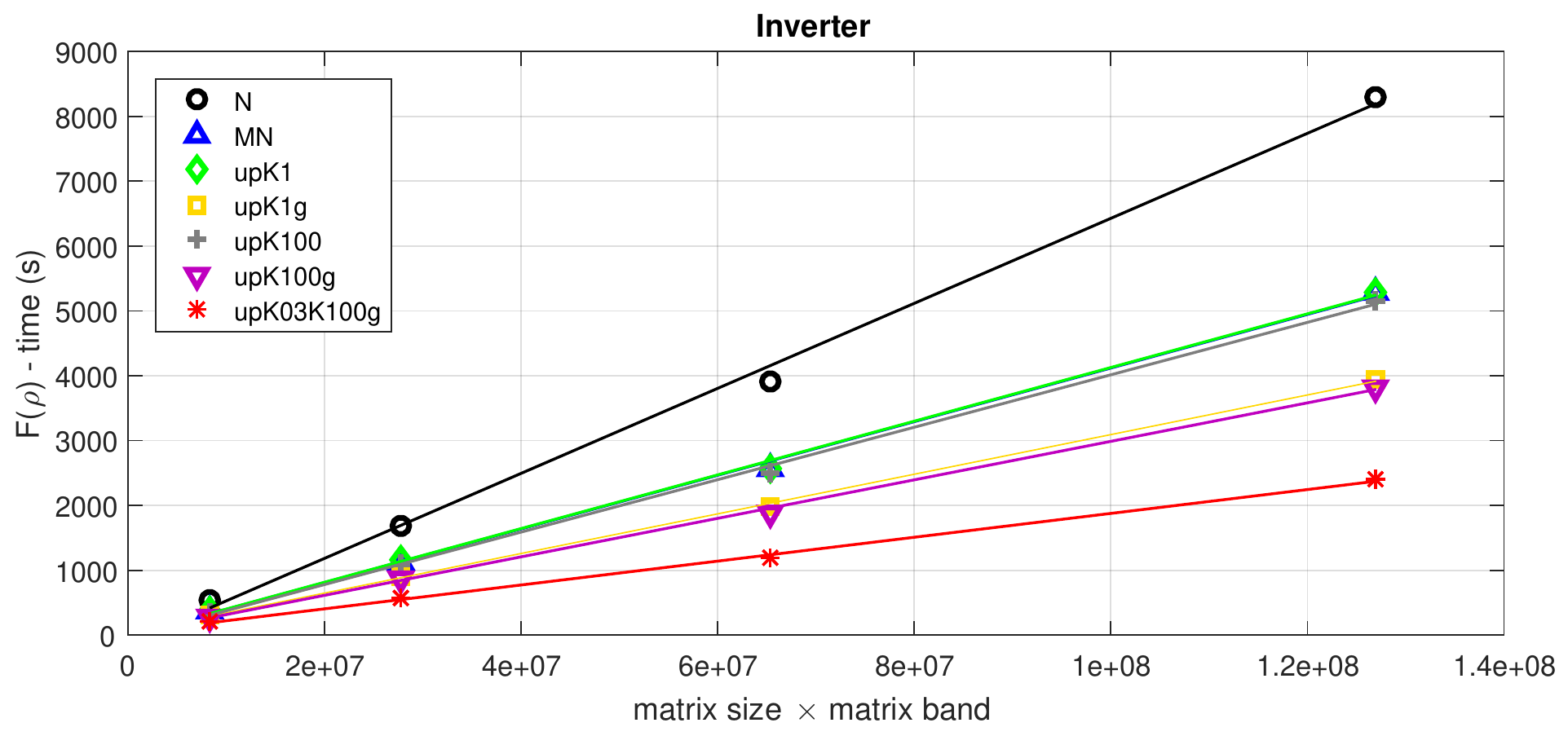}
      \caption{$F(\bsrho)$-time for each strategy according to the dimension of $\KK$ (inverter).} 
      \label{time_fobj_discretization_inverter}
\end{figure}

Figure \ref{time_fobj_discretization_inverter} shows how the time spent to compute $F(\bsrho)$ depends on the problem size. Once again, we observe that our strategies are very efficient in reducing the time consumed in computing $F(\bsrho)$\, as the mesh is refined. The regression lines suggest that the growth rates of the Modified Newton method and strategies {\tt upK1} and {\tt upK100} are very similar. The same occurs with {\tt upK1g} and {\tt upK100g}. As expected, the {\tt upK03K100g} strategy achieved the best performance. Compared to the regression line of Newton's method, the lines of the Modified Newton method and strategies {\tt upK1} and {\tt upK100} showed a 37\,\%\, smaller slope. For strategies {\tt upK1g} and {\tt upK100g}, the slope is $54.1\,\%$\, smaller, and for {\tt upK03K100g} the reduction is about $72.0\,\%$.

\begin{figure}[htbp!]
   \centering 
      \begin{tabular}{cc} 
          \includegraphics[width=0.45\textwidth]{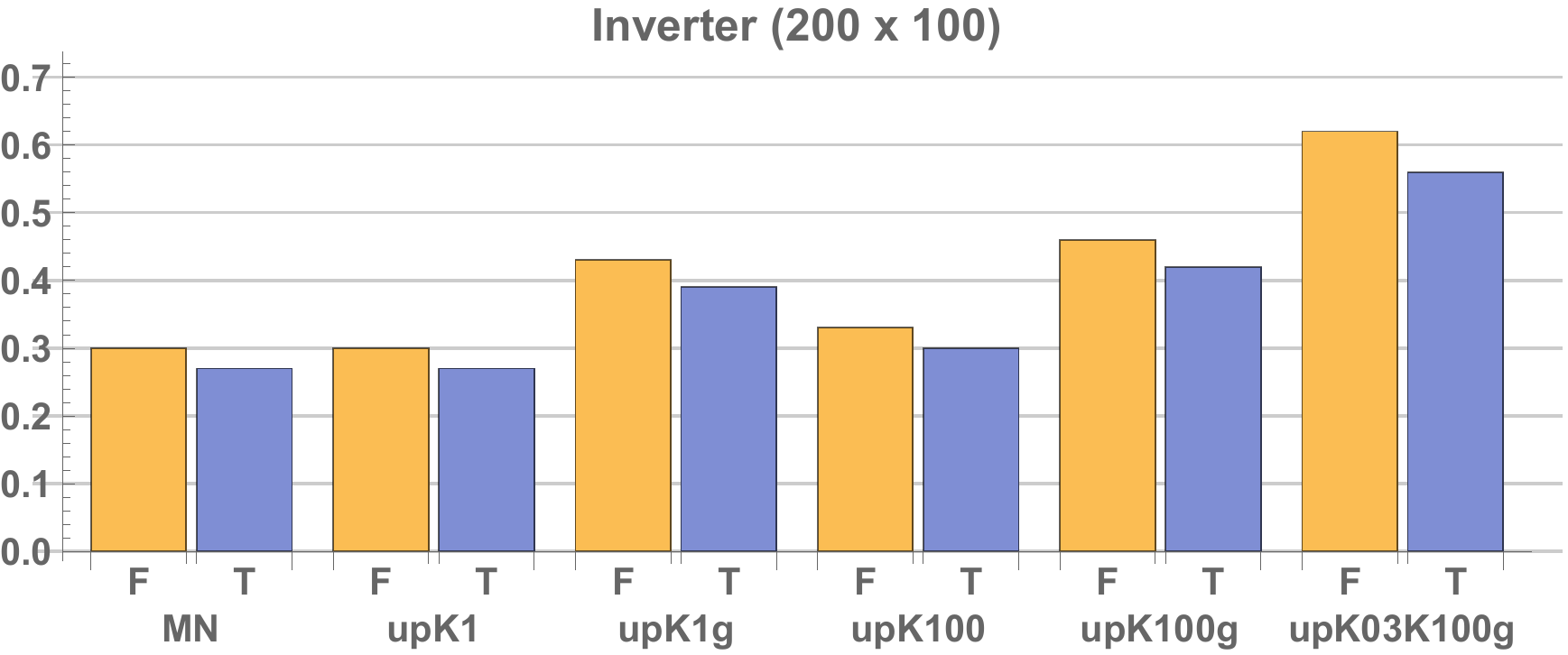}
        & \hspace{0.3cm}\includegraphics[width=0.45\textwidth]{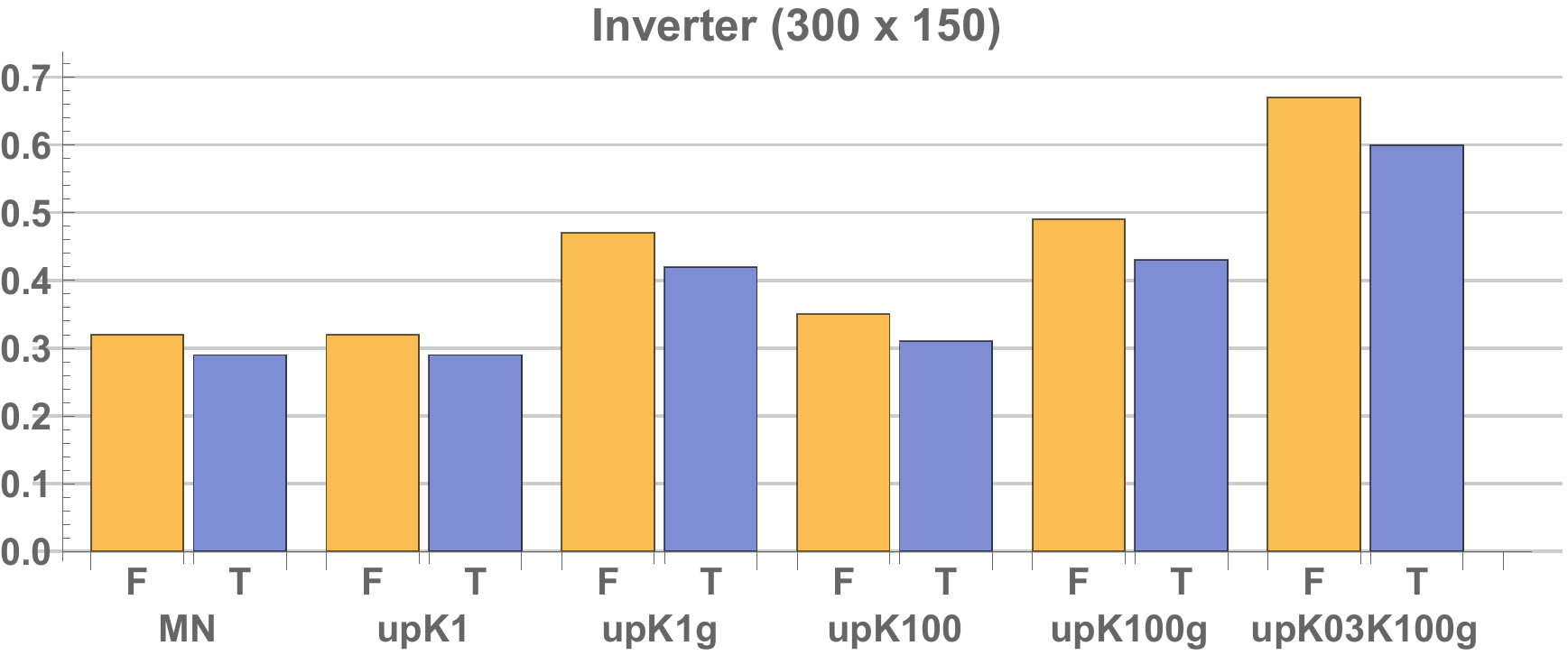} \\[10pt]
          \includegraphics[width=0.45\textwidth]{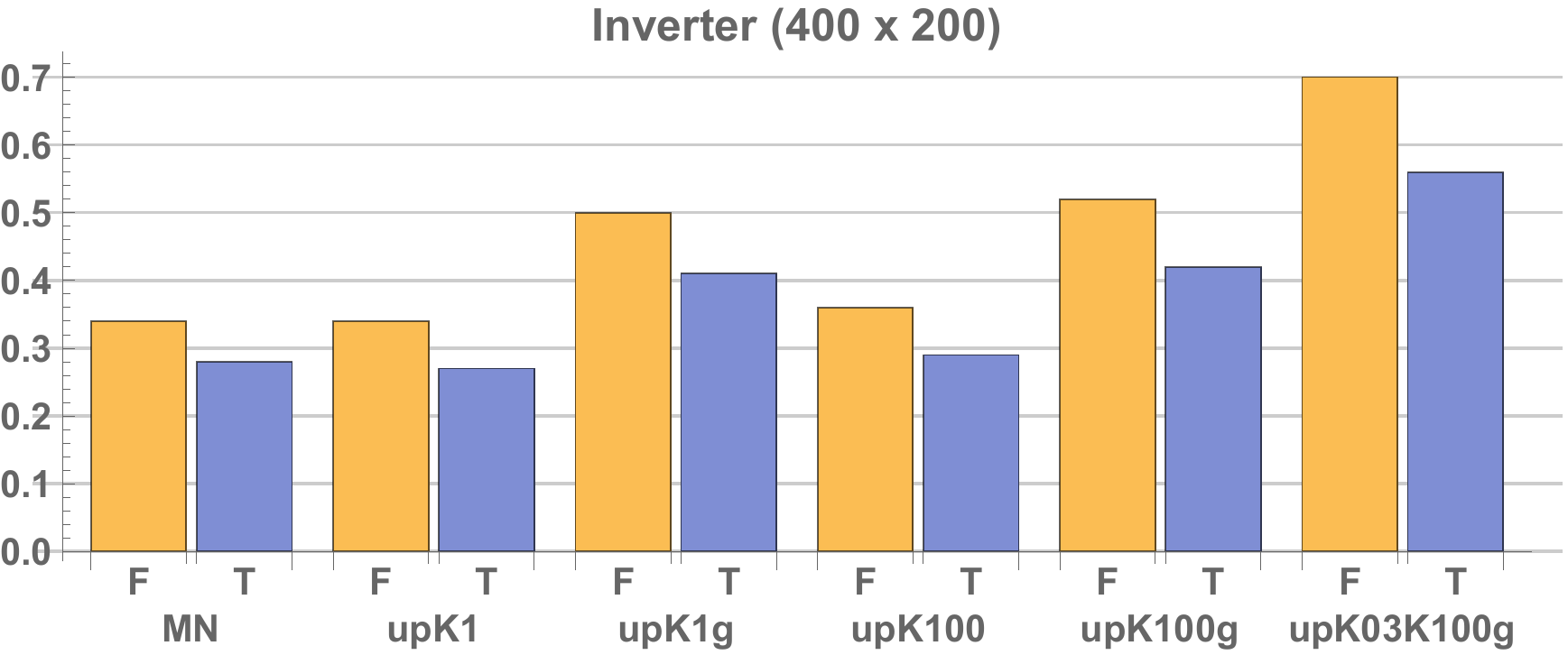}
        & \hspace{0.3cm}\includegraphics[width=0.45\textwidth]{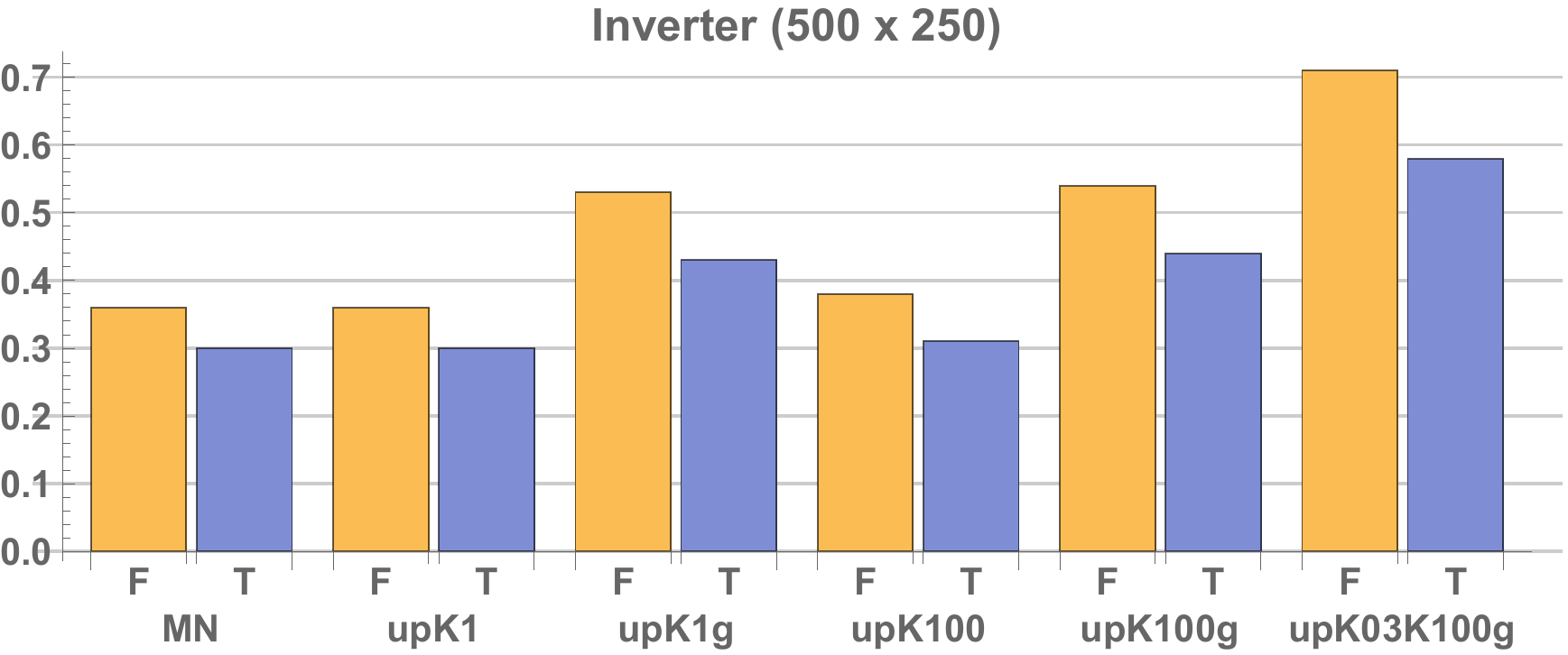} \\ 
      \end{tabular}
\caption{Percentage of decrease in the time demanded to evaluate the objective function (\textsf{F}) and in the total time (\textsf{T}), as compared with Newton's method, for the inverter mechanism.}
 \label{percentage_reduction_fobj_inverter}
\end{figure}

Figure \ref{percentage_reduction_fobj_inverter} presents the percentage of decrease in the time spent to evaluate $F(\bsrho)$ and in the overall time of the algorithm, with respect to Newton's method. 
From the figure, one can notice that the percentages provided by the Modified Newton method and strategies {\tt upK1} and {\tt upK100} are very similar, ranging from $31.0\,\%$ to $37.0\,\%$ when we take into account the four mesh sizes.

Analogously to the results shown in Fig. \ref{percentage_reduction_fobj_sb}, the adoption of the iterative scheme (\ref{eq03}) to obtain an approximate solution for the adjoint system (\ref{eq09}) resulted in a substantial reduction in the time spent with the reanalysis process. In fact, strategies {\tt upK1g} and {\tt upK100g} provided an average decrease that ranged from $44.5\,\%$ (for the $200 \times 100$ mesh) to $53.3\,\%$ (for the $500 \times 250$ one).

As before, strategy {\tt upK03K100g} showed the best results for all of the meshes considered, providing a reduction that varies from $62.2\,\%$ (for the coarser mesh) to $71.2\,\%$ (for the finer one). On average, this strategy spent $46.0\,\%$ less time to assemble $\KK$ and $78.7\,\%$ less time to factorize this matrix, in comparison with Newton's method.

For the Modified Newton method and strategies {\tt upK1} and {\tt upK100}, the overall time spent by the SPLP algorithm was $28.6\,\%$, $28.1\,\%$ and $30.4\,\%$ smaller, respectively, than the time consumed when Newton's method was adopted. For strategies {\tt upK1g} and {\tt upK100g}, we have observed a reduction of $41.0\,\%$ and $42.8\,\%$, respectively. Finally, the {\tt upK03K100g} strategy showed an average decrease of $57.4\,\%$.

%===========================================
\section{Performance of the SPLP method with the nonlinear reanalysis}\label{sec:convsplp}

In this section, we investigate the performance of the reanalysis strategies when the SPLP method is run until convergence. We have seen above that, for a fixed budget of 300 iterations, these strategies can substantially reduce the overall time spent by the algorithm, as compared to Newton's method. Our purpose now is to check if this also happens when we adopt a mathematically well founded stopping criterion, as described below.

The Lagrangian function associated with the topology optimization problem is defined by
    $$\mathcal{L}(\bsrho,\,\theta) \, := \, F(\bsrho) + \theta V(\bsrho),$$
where $F(\bsrho)$ is the objective function, $V(\bsrho)$ is the volume constraint, and $\theta \in \mathbb{R}$ is the Lagrange multiplier associated with this constraint. Denoting the gradient vector of the Lagrangian function by $\nabla_{\bsrho}\mathcal{L}(\bsrho,\,\theta)$  and defining $v(\bsrho,\,\theta) := \bsrho-\nabla_{\bsrho}\mathcal{L}(\bsrho,\,\theta)$, the continuous projected gradient 
onto the set 
    $$X := \left\{\bsrho \in \mathbb{R}^{n_{el}} \,\,\, | \,\,\, \rho_{\min} \leq \rho_{i} \leq 1, \,\,\, i = 1,\,\dots,\, n_{el} \right\}$$
is given by $g_{P}(\bsrho,\,\theta) = P_{X}(v(\bsrho,\,\theta))-\bsrho$, where $P_{X}(v(\bsrho,\,\theta))$ is the orthogonal projection of the vector $v(\bsrho,\,\theta)$ onto $X$.
We consider that the SPLP algorithm has found a good approximation for a stationary point for the topology optimization problem whenever
    \begin{equation}
      \| g_{P}(\bsrho^{(k_{s})},\,\theta^{(k_{s})}) \|_{\infty} < 10^{-3},
			\label{stop_crit_eq}
    \end{equation}
where $\| \, \cdot \, \|_{\infty}$ denotes the max-norm, and $\theta^{(k_{s})}$ is the approximate Lagrange multiplier of the volume constraint at the $k_{s}$-th iteration of the SPLP algorithm.  

The results obtained for the rigid structures and compliant mechanisms using the criterion~\eqref{stop_crit_eq} are shown in Tables \ref{table_convergence_cb_sb} and \ref{table_convergence_inverter_gripper}, respectively. As we observe, no matter the problem, the optimal values of the objective function at the solution were very similar among strategies. On the other hand, with the exception of the gripper, the number of iterations of the SPLP algorithm showed some variation, affecting the overall performance of the algorithm. 

The performance of the various strategies was also affected by the number of iterations required for solving the nonlinear system (\ref{eq01b}). For the slender beam, for example, most of the strategies took more iterations than Newton's method, so the decrease in the total CPU time was less noticeable. An exception to this rule, however, is seen when {\tt upK03K100g} is applied to the gripper. In this case, the total time was reduced by more than 60\,\%, if compared to Newton's method, even though the number of iterations was 31\,\% higher.

\begin{table}[htbp!]
    \centering
    \caption{Results obtained using the SPLP stopping criterion (cantilever and slender beams).} \vspace{0.3cm}
    \label{table_convergence_cb_sb}
        \begin{tabular}{lc@{}c@{}c@{}c@{}c@{}c@{}c@{}}
        \hline\noalign{\smallskip}
        \textbf{Cantilever beam} & \multirow{2}{*}{\tt N} & \multirow{2}{*}{\tt MN} & \multirow{2}{*}{\tt upK1} & \multirow{2}{*}{\tt upK1g} & \multirow{2}{*}{\tt upK100} & \multirow{2}{*}{\tt upK100g} & {\tt upK03} \\
        $\mathbf{400 \times 100}$ &  &  &  &  &  &  & {\tt K100g} \\ 
        \noalign{\smallskip}\hline\noalign{\smallskip} 
        $F(\bsrho)$ - value        & $2019.60$ \ & \ $2020.10$  \ & \  $2020.01$  \ & \  $2019.91$   \ & \  $2020.10$   \ & \  $2019.91$   \ & \  $2020.04$  \\
        \# SPLP iter.         & $1378$    & $1061$     & $1237$     & $1303$      & $1079$      & $1303$      & $1100$     \\
        \# Newton iter.       & $2846$    & $3565$     & $3032$     & $2832$      & $2721$      & $3172$      & $2839$     \\  
        \hline
        \noalign{\smallskip}
        CPU time (s)  &  &  &  &  &  &  &  \\ 
        Total                & $6242.60$ & $3670.59$  & $4216.19$  & $3777.00$   & $3570.34$   & $3623.42$   & $1919.11$  \\
	    $F(\bsrho)$         & $5025.19$ & $2759.84$  & $3102.30$  & $2593.31$   & $2592.30$   & $2430.46$   & $911.78$   \\     
        $\KK$               & $658.27$  & $587.72$   & $639.92$   & $672.96$    & $419.56$    & $512.51$    & $269.09$   \\
        RHS                  & $106.94$  & $116.85$   & $118.00$   & $153.81$    & $106.00$    & $154.65$    & $156.64$   \\
	    Factorizations       & $4209.96$ & $1999.73$  & $2270.70$  & $1630.75$   & $1999.47$   & $1626.52$   & $331.50$   \\
	    Linear systems      & $50.02$   & $55.54$    & $73.68$    & $135.79$    & $67.27$     & $136.78$    & $154.55$   \\
	    $\nabla F(\bsrho)$   & $420.00$  & $306.87$   & $383.42$   & $417.06$    & $338.85$    & $423.35$    & $343.47$   \\
	    SPLP solving         & $426.90$  & $335.56$   & $391.65$   & $397.67$    & $340.59$    & $396.76$    & $358.62$   \\
	    Filtering            & $362.58$  & $260.98$   & $330.89$   & $360.98$    & $291.09$    & $364.73$    & $297.84$   \\
	    Other              & $7.93$    & $7.34$     & $7.93$     & $7.98$      & $7.51$      & $8.12$      & $7.40$     \\
	    \hline\hline\noalign{\smallskip}
                \textbf{Slender beam} & \multirow{2}{*}{\tt N} & \multirow{2}{*}{\tt MN} & \multirow{2}{*}{\tt upK1} & \multirow{2}{*}{\tt upK1g} & \multirow{2}{*}{\tt upK100} & \multirow{2}{*}{\tt upK100g} & {\tt upK03} \\
        $\mathbf{600 \times 75}$ &  &  &  &  &  &  & {\tt K100g} \\ 
        \noalign{\smallskip}\hline\noalign{\smallskip} 
        $F(\bsrho)$ - value           & $137.67$  & $137.64$   & $137.66$   & $137.61$   & $137.66$   & $137.61$  & $137.67$  \\
        \# SPLP iter.            & $541$     & $625$      & $686$      & $625$      & $686$      & $625$     & $541$     \\
        \# Newton iter.         & $1130$    & $1399$     & $1424$     & $1317$     & $1424$     & $1320$    & $1448$    \\ 
        \hline
        \noalign{\smallskip}
        CPU time (s)  &  &  &  &  &  &  &  \\  
        Total           & $2613.13$ & $2217.55$  & $2452.18$  & $2078.82$  & $2491.74$  & $1929.06$ & $1181.73$ \\
	    $F(\bsrho)$        & $2081.28$ & $1651.55$  & $1847.23$  & $1471.68$  & $1871.53$  & $1337.04$ & $647.42$  \\     
        $\KK$               & $303.87$  & $363.87$   & $378.22$   & $374.06$   & $318.36$   & $291.11$  & $164.61$  \\
        RHS                   & $47.29$   & $58.03$    & $64.01$    & $81.30$    & $66.25$    & $80.72$   & $94.82$   \\
	    Factorizations     & $1708.34$ & $1203.17$  & $1367.30$  & $942.58$   & $1447.84$  & $888.14$  & $290.56$  \\
	    Linear systems       & $21.78$   & $26.48$    & $37.70$    & $73.74$    & $39.08$    & $77.07$   & $97.43$   \\
	    $\nabla F(\bsrho)$ s    & $60.66$   & $70.76$    & $77.47$    & $73.39$    & $79.77$    & $71.44$   & $62.05$   \\
	    SPLP solvings        & $424.32$  & $440.87$   & $468.25$   & $477.47$   & $479.31$   & $465.33$  & $424.79$  \\
	    Filterings           & $44.78$   & $52.06$    & $56.95$    & $54.02$    & $58.84$    & $52.98$   & $45.36$   \\
	    Others               & $2.09$    & $2.31$     & $2.28$     & $2.26$     & $2.29$     & $2.27$    & $2.11$    \\	    
        \noalign{\smallskip}\hline
    \end{tabular}
\end{table} 

\begin{table}[htbp!]
    \centering
    \caption{Results obtained using the SPLP stopping criterion (Inverter and Gripper).} \vspace{0.3cm}
    \label{table_convergence_inverter_gripper}
    \begin{tabular}{lc@{}c@{}c@{}c@{}c@{}c@{}c@{}}
        \hline\noalign{\smallskip}\
                        \textbf{Inverter} & \multirow{2}{*}{\tt N} & \multirow{2}{*}{\tt MN} & \multirow{2}{*}{\tt upK1} & \multirow{2}{*}{\tt upK1g} & \multirow{2}{*}{\tt upK100} & \multirow{2}{*}{\tt upK100g} & {\tt upK03} \\
        $\mathbf{300 \times 150}$ &  &  &  &  &  &  & {\tt K100g} \\ 
        \noalign{\smallskip}\hline\noalign{\smallskip} 
        $F(\bsrho)$ - value           & $-8.7438$ \, & $-8.7276$ \,& $-8.7020$ \,& $-8.7278$ \,& $-8.7221$ \,& $-8.7337$ \,& $-8.7248$   \\
        \# SPLP iter.            & $439$     & $367$     & $325$     & $351$     & $350$     & $382$     & $342$       \\
        \# Newton iter.          & $835$     & $1236$    & $929$     & $967$     & $971$     & $1039$    & $1033$      \\
        \hline
        \noalign{\smallskip}
        CPU time (s)  &  &  &  &  &  &  &  \\ 
        Total             & $2660.33$ & $1592.05$ & $1452.41$ & $1258.17$ & $1482.75$ & $1311.61$ & $821.47$    \\
	    $F(\bsrho)$          & $2409.59$ & $1374.04$ & $1246.94$ & $1044.12$ & $1270.01$ & $1070.33$ & $610.67$    \\     
        $\KK$                  & $236.58$  & $230.54$  & $193.87$  & $209.72$  & $148.21$  & $169.34$  & $117.72$    \\
        RHS                  & $40.89$   & $45.46$   & $46.21$   & $61.54$   & $48.01$   & $66.47$   & $73.94$     \\
	    Factorizations         & $2112.25$ & $1075.41$ & $971.24$  & $708.80$  & $1037.80$ & $765.76$  & $330.75$    \\
	    Linear systems           & $19.87$   & $22.63$   & $35.62$   & $64.06$   & $35.99$   & $68.76$   & $88.26$     \\
	    $\nabla F(\bsrho)$    & $51.00$   & $43.70$   & $40.62$   & $42.09$   & $42.01$   & $45.88$   & $41.45$     \\
	    SPLP solving          & $160.26$  & $140.06$  & $132.53$  & $138.58$  & $137.77$  & $158.82$  & $136.88$    \\
	    Filtering              & $37.62$   & $32.34$   & $30.31$   & $31.44$   & $31.09$   & $34.60$   & $30.51$     \\
	    Other                 & $1.86$    & $1.91$    & $2.01$    & $1.94$    & $1.87$    & $1.98$    & $1.96$      \\
	    \hline\hline\noalign{\smallskip}
                        \textbf{Gripper} & \multirow{2}{*}{\tt N} & \multirow{2}{*}{\tt MN} & \multirow{2}{*}{\tt upK1} & \multirow{2}{*}{\tt upK1g} & \multirow{2}{*}{\tt upK100} & \multirow{2}{*}{\tt upK100g} & {\tt upK03} \\
        $\mathbf{320 \times 160}$ &  &  &  &  &  &  & {\tt K100g} \\ 
        \noalign{\smallskip}\hline\noalign{\smallskip} 
        $F(\bsrho)$ - value           & $-3.8542$ & $-3.8504$ & $-3.8567$  & $-3.8546$ & $-3.8596$ & $-3.8559$ & $-3.8531$   \\
        \# SPLP iter.            & $314$     & $308$     & $305$      & $316$     & $308$     & $319$     & $315$       \\
        \# Newton iter.          & $687$     & $904$     & $741$      & $756$     & $749$     & $782$     & $902$       \\ 
       \hline 
       \noalign{\smallskip}
        CPU time (s)  &  &  &  &  &  &  &  \\ 
        Total            & $2642.35$ & $1876.04$ & $1847.92$  & $1562.44$ & $1733.29$ & $1520.73$ & $968.43$    \\
	    $F(\bsrho)$          & $2287.84$ & $1524.55$ & $1503.31$  & $1202.70$ & $1464.60$ & $1155.81$ & $604.06$    \\     
        $\KK$               & $217.67$  & $242.79$  & $222.47$   & $233.43$  & $174.84$  & $185.27$  & $134.99$    \\
        RHS                 & $33.77$   & $40.12$   & $42.44$    & $54.76$   & $43.21$   & $56.22$   & $7.62$      \\
	    Factorizations         & $2019.77$ & $1221.37$ & $1208.39$  & $861.13$  & $1215.83$ & $859.42$  & $367.86$    \\
	    Linear systems          & $16.63$   & $20.27$   & $30.01$    & $53.38$   & $30.72$   & $54.90$   & $93.59$     \\
	    $\nabla F(\bsrho)$     & $107.38$  & $107.79$  & $105.02$   & $111.89$  & $107.22$  & $113.11$  & $113.77$    \\
	    SPLP solving        & $151.46$  & $148.18$  & $146.67$   & $150.21$  & $146.36$  & $152.98$  & $152.73$    \\
	    Filtering          & $89.69$   & $89.19$   & $86.79$    & $90.96$   & $8.86$    & $92.02$   & $91.03$     \\
	    Other                  & $5.98$    & $6.33$    & $6.13$     & $6.68$    & $6.25$    & $6.81$    & $6.84$      \\    
        \noalign{\smallskip}\hline
    \end{tabular}
\end{table} 

\begin{figure}[htbp!]
   \centering 
      \begin{tabular}{cc} 
          \includegraphics[width=0.45\textwidth]{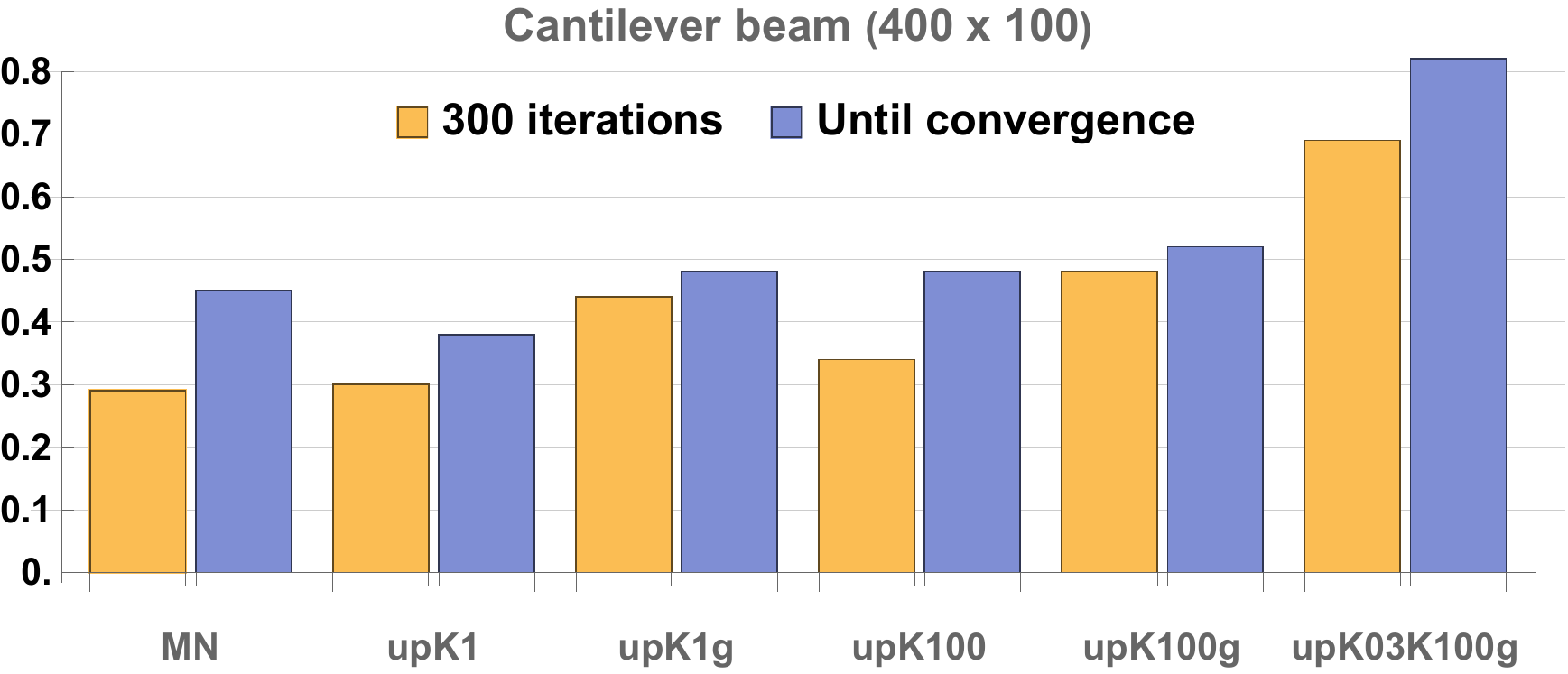}
        & \hspace{0.3cm}\includegraphics[width=0.45\textwidth]{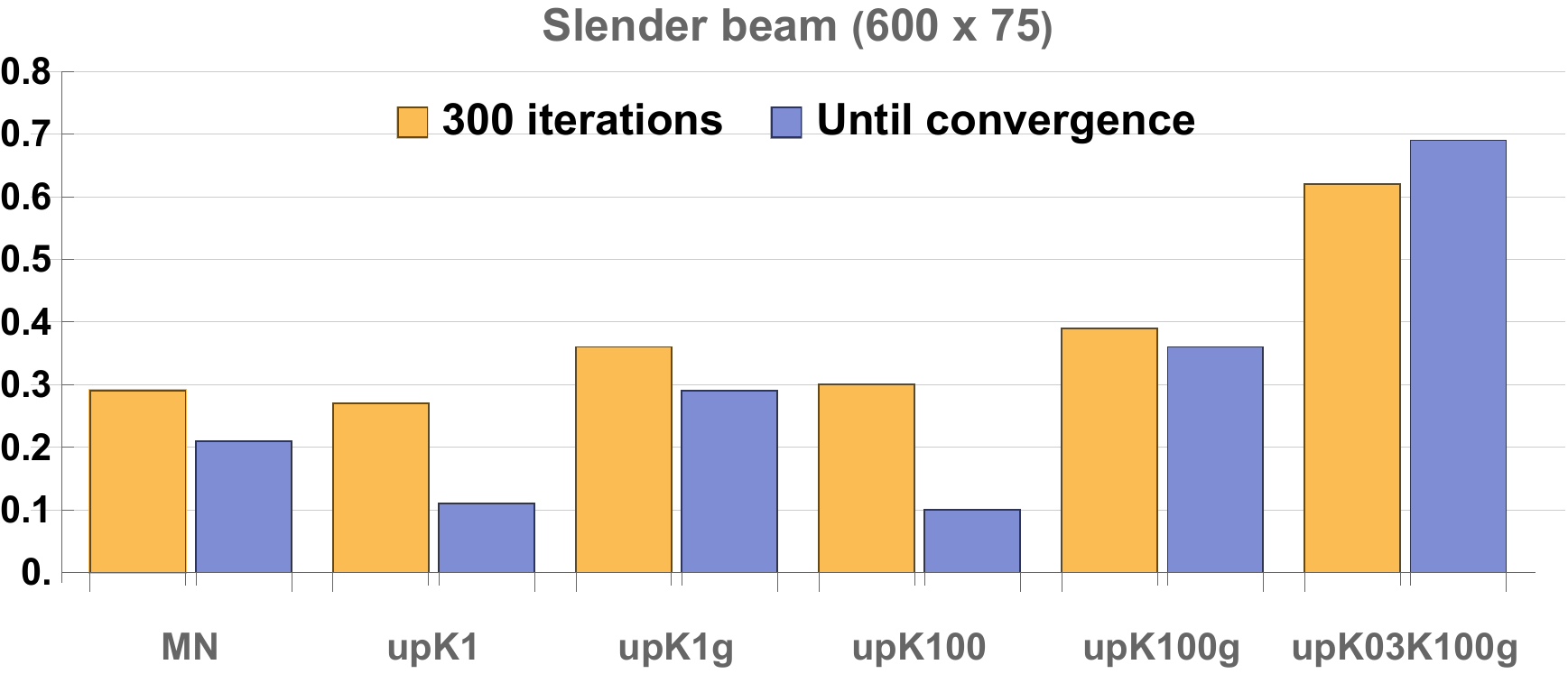} \\[10pt]
          \includegraphics[width=0.45\textwidth]{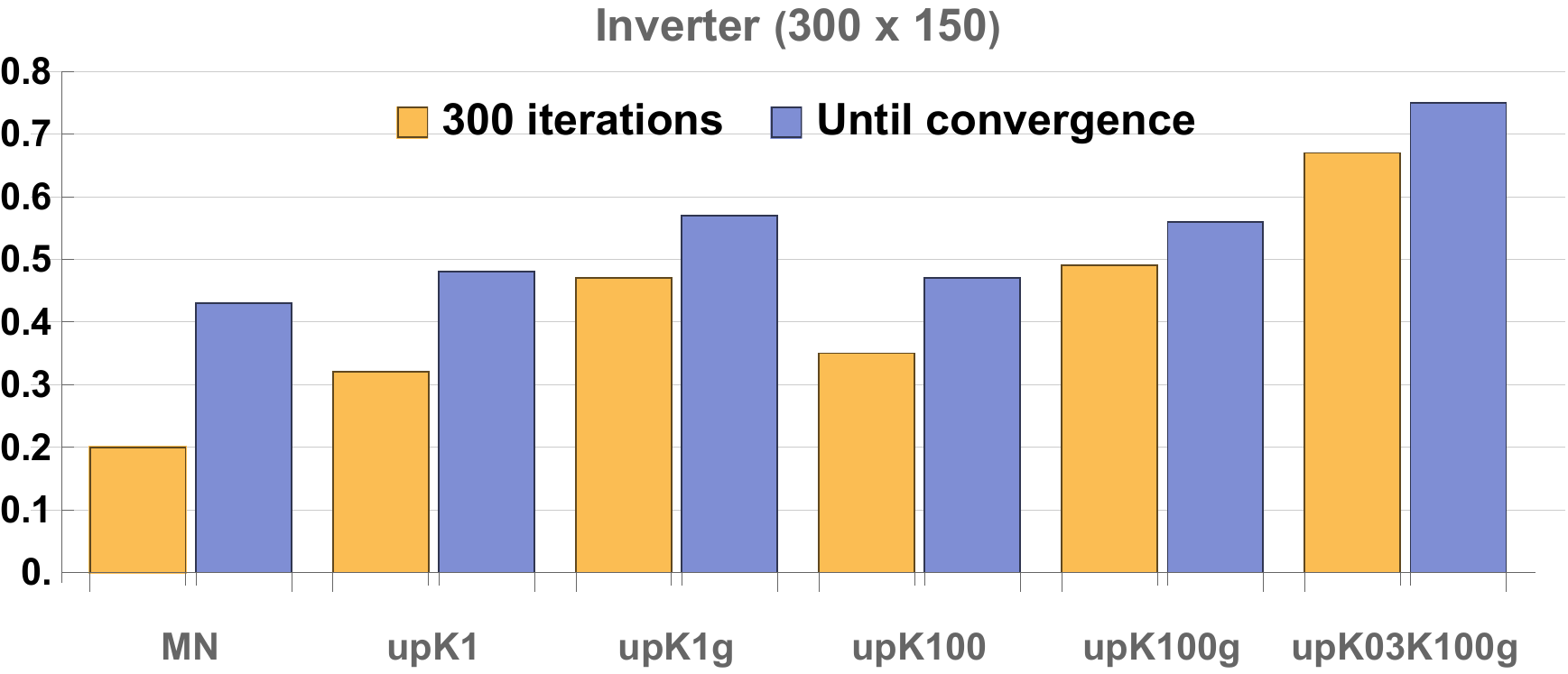}
        & \hspace{0.3cm}\includegraphics[width=0.45\textwidth]{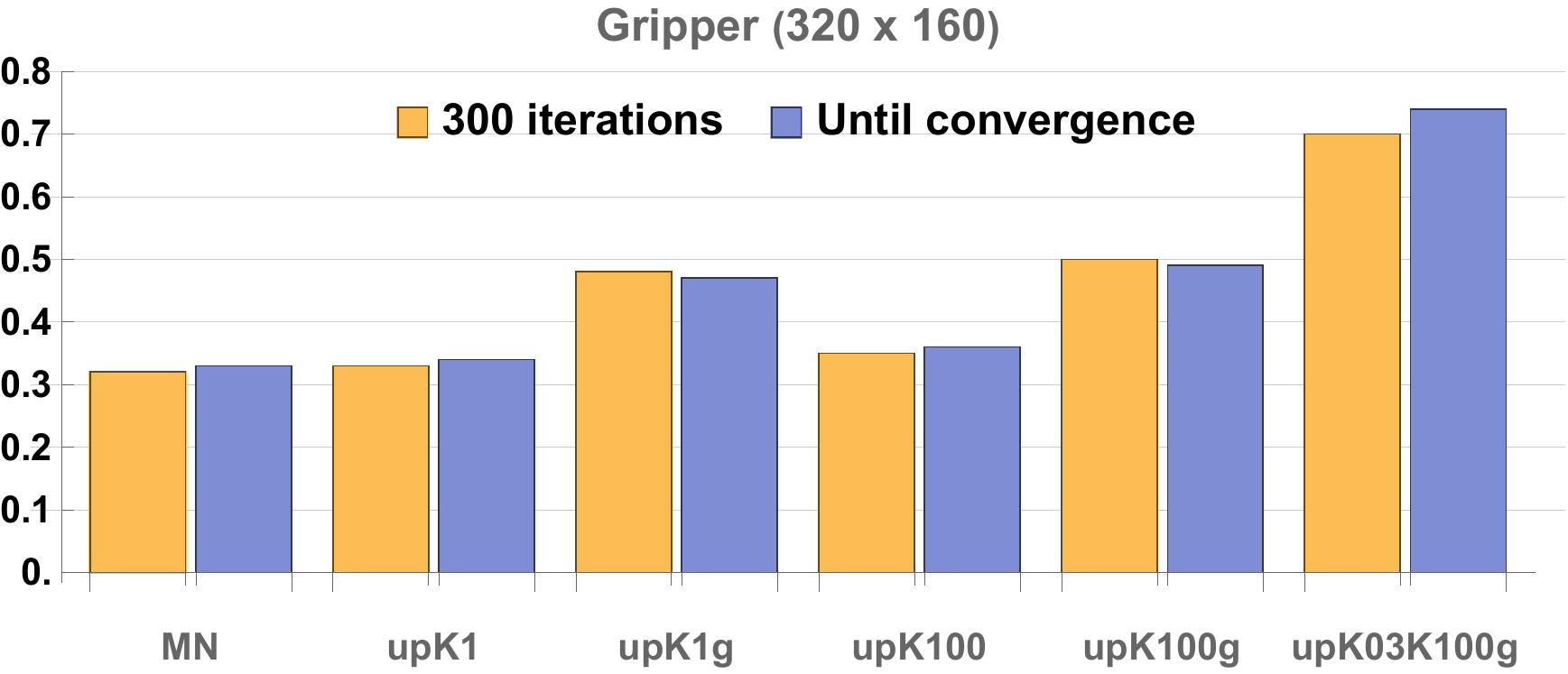} \\ 
      \end{tabular}
\caption{Percentage of decrease in the time required to evaluate the objective function, as compared to Newton's method.}
 \label{percentage_reduction_fobj_full_conv}
\end{figure}

Figure \ref{percentage_reduction_fobj_full_conv} summarizes the differences observed on the time the various strategies spent to compute the objective function, both for the fixed budget of 300 iterations and in the case we require full convergence. The bars of the figure show the reduction provided by each strategy, with respect to Newton's method. From the figure, it is easy to notice that, for the cantilever beam and the inverter, the results obtained when we use the stopping criterion (\ref{stop_crit_eq}) are better than those attained with a fixed budget. For the gripper, there was almost no difference between the results, and for the slender beam the performance of the first five strategies was degraded when full convergence was required, mainly due to an increase on the number of iterations of Newton's method, as mentioned above. The figure also reveals the relevance of applying the ICA scheme to the computation of the objective function gradient. In general, the strategies whose name end with {\tt g} (i.e., those that use ICA for solving (\ref{eq09})) attained the best results, notably {\tt upK03K100g}, that showed a consistent improvement when the stopping criterion was changed.

%===========================================
\section{Investigating the behavior of the norm of matrix \emph{B}}\label{sec:normB}

Let $\BB$ be the matrix defined in (\ref{eq:stilde&matB}). In Section \ref{sec:inca}, we demonstrated that if $\| \BB \| < 1$, then the sequence $\left\{\sk\right\}^{+\infty}_{k=0}$ (defined in (\ref{eq03}) and associated with the ICA method) converges to the solution $\sstar$ of the linear system (\ref{eq02}) of Newton's method. Moreover, we showed that this sequence has linear rate of convergence. 

With the aim of investigating the impact of the variation on the maximum value of $\| \BB \|_{2}$ attained whenever the nonlinear system (\ref{eq01b}) is solved, we have performed an additional test considering the slender beam and the gripper mechanism, using strategies {\tt upK100g} and {\tt upK03K100g}. Since the explicit evaluation of $\BB$ is very expensive, these numerical test was performed in MATLAB (version R2016a) and took into account a coarser mesh, with $160 \times 20 = 3200$ finite elements for the slender beam and $80 \times 40 = 3200$ finite elements for the gripper mechanism, applying the stopping criterion for the SPLP method stated in (\ref{stop_crit_eq}). 

\begin{figure}[htbp!]
   \centering 
      \begin{tabular}{ccc} 
			    \raisebox{1.1cm}{\rotatebox[origin=lb]{90}{{\tt upK100g}}} \hspace{-0.2cm}
        & \includegraphics[width=0.45\textwidth]{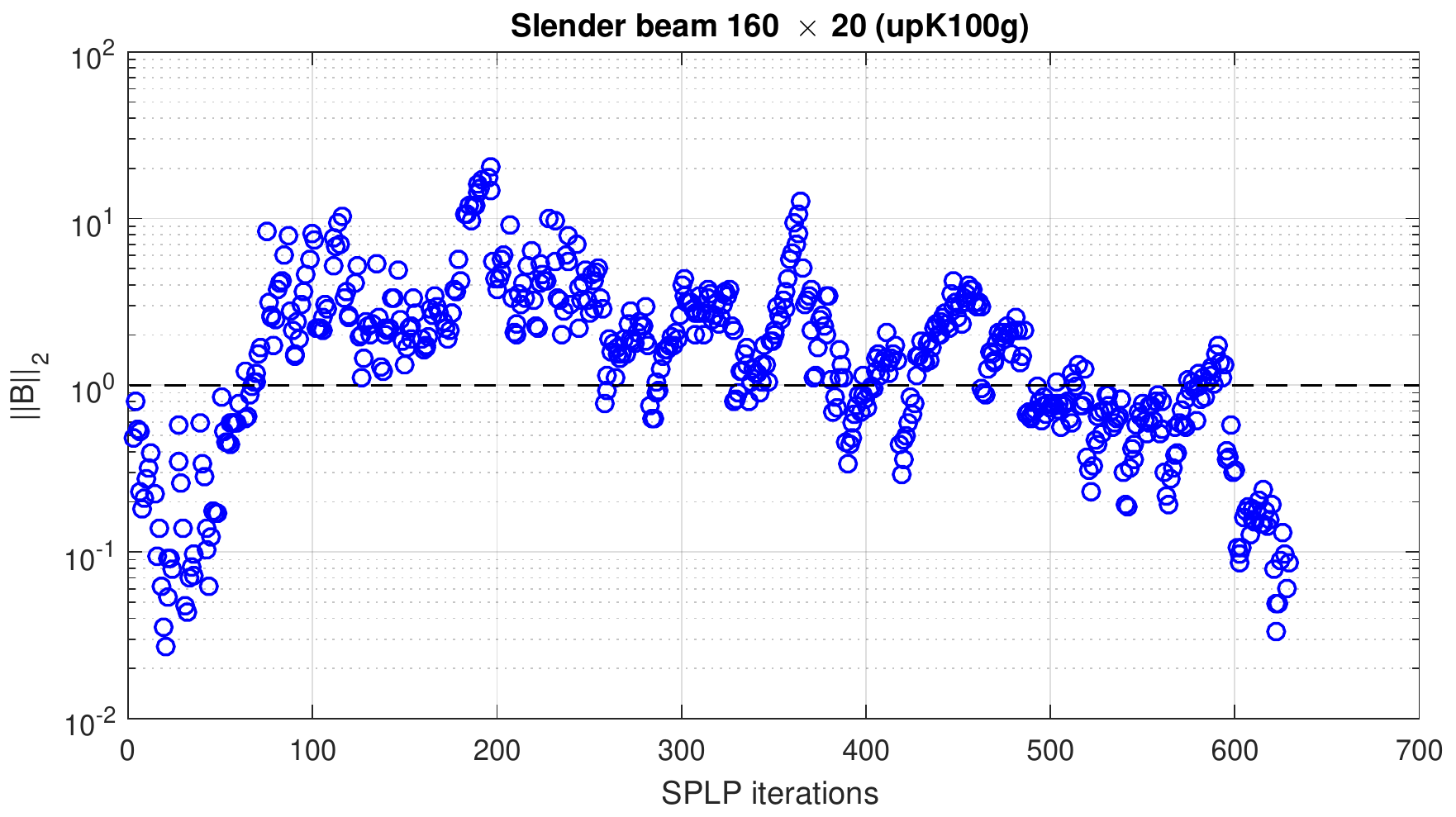}
        & \includegraphics[width=0.45\textwidth]{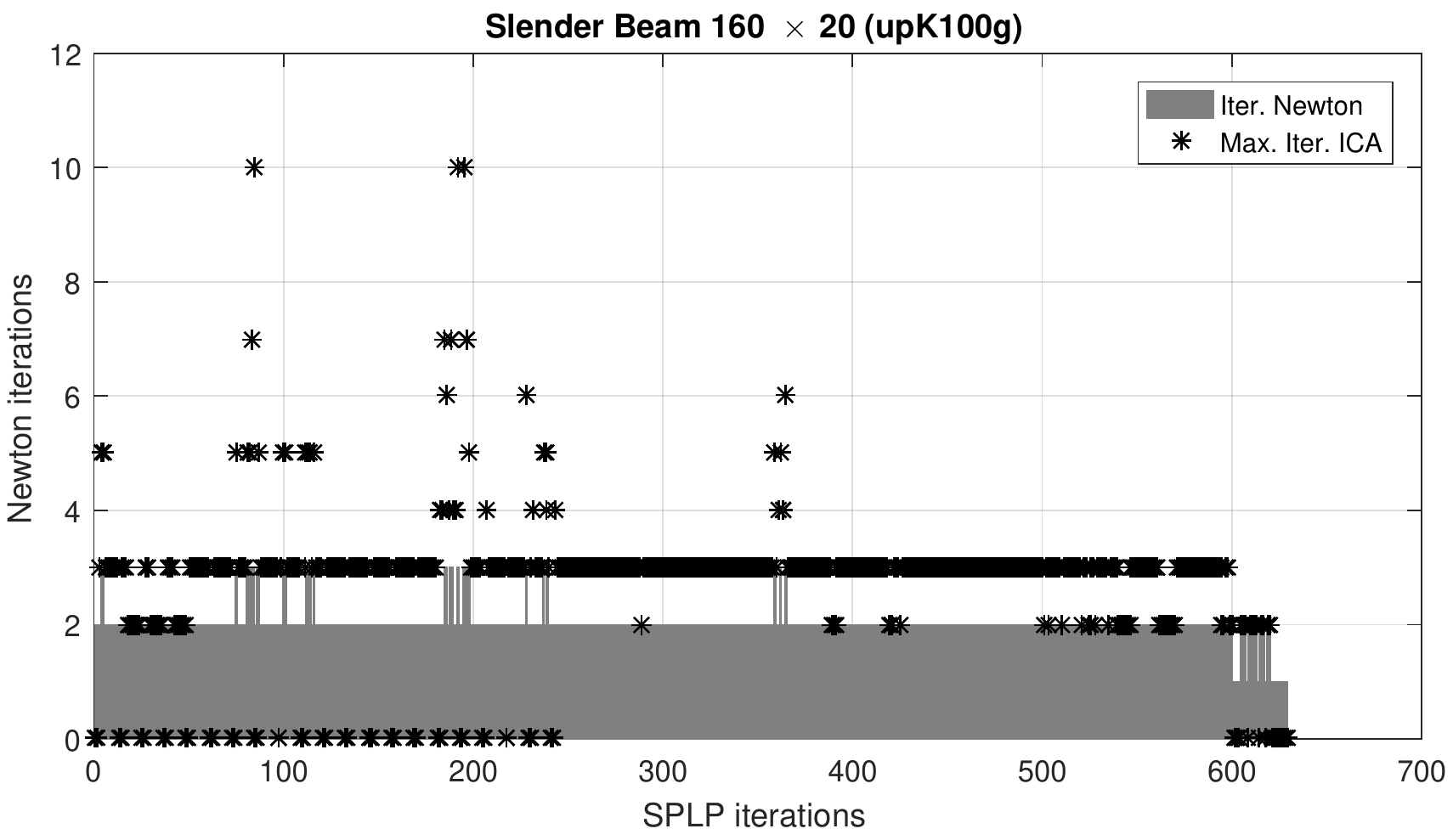} \\[0.3cm]
			    \raisebox{0.9cm}{\rotatebox[origin=lb]{90}{{\tt upK03K100g}}} \hspace{-0.2cm}
        & \includegraphics[width=0.45\textwidth]{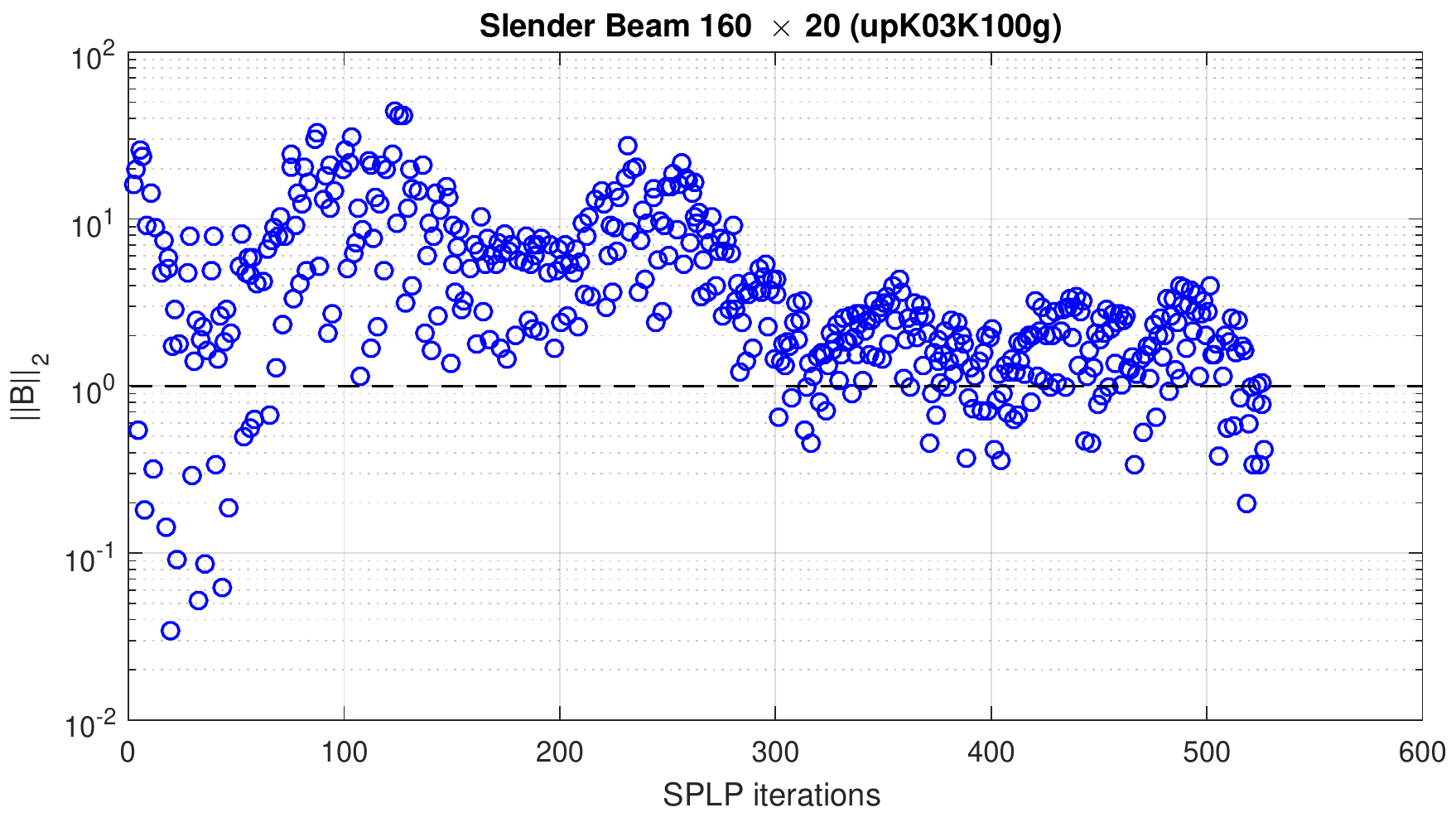} 
        & \includegraphics[width=0.45\textwidth]{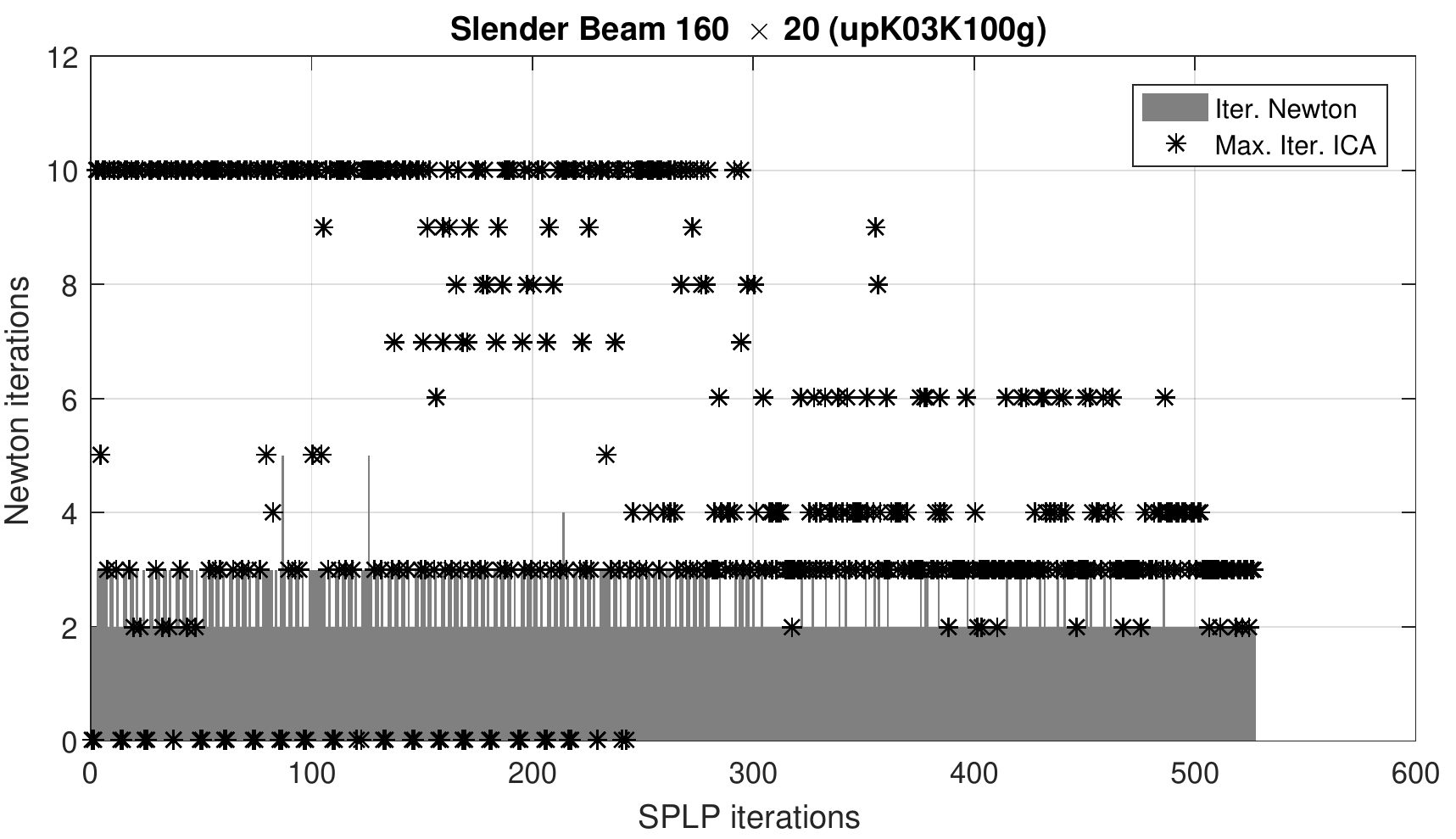}
      \end{tabular}
      \caption{Maximum value of $\| \BB \|_{2}$ (left) and number of iterations of Newton's method (right) for the slender beam $160 \times 20$.} 
      \label{normb_itnewton_sb160x20}
\end{figure}

\begin{figure}[htbp!]
   \centering 
      \begin{tabular}{ccc} 
			    \raisebox{1.1cm}{\rotatebox[origin=lb]{90}{{\tt upK100g}}} \hspace{-0.2cm}
        & \includegraphics[width=0.45\textwidth]{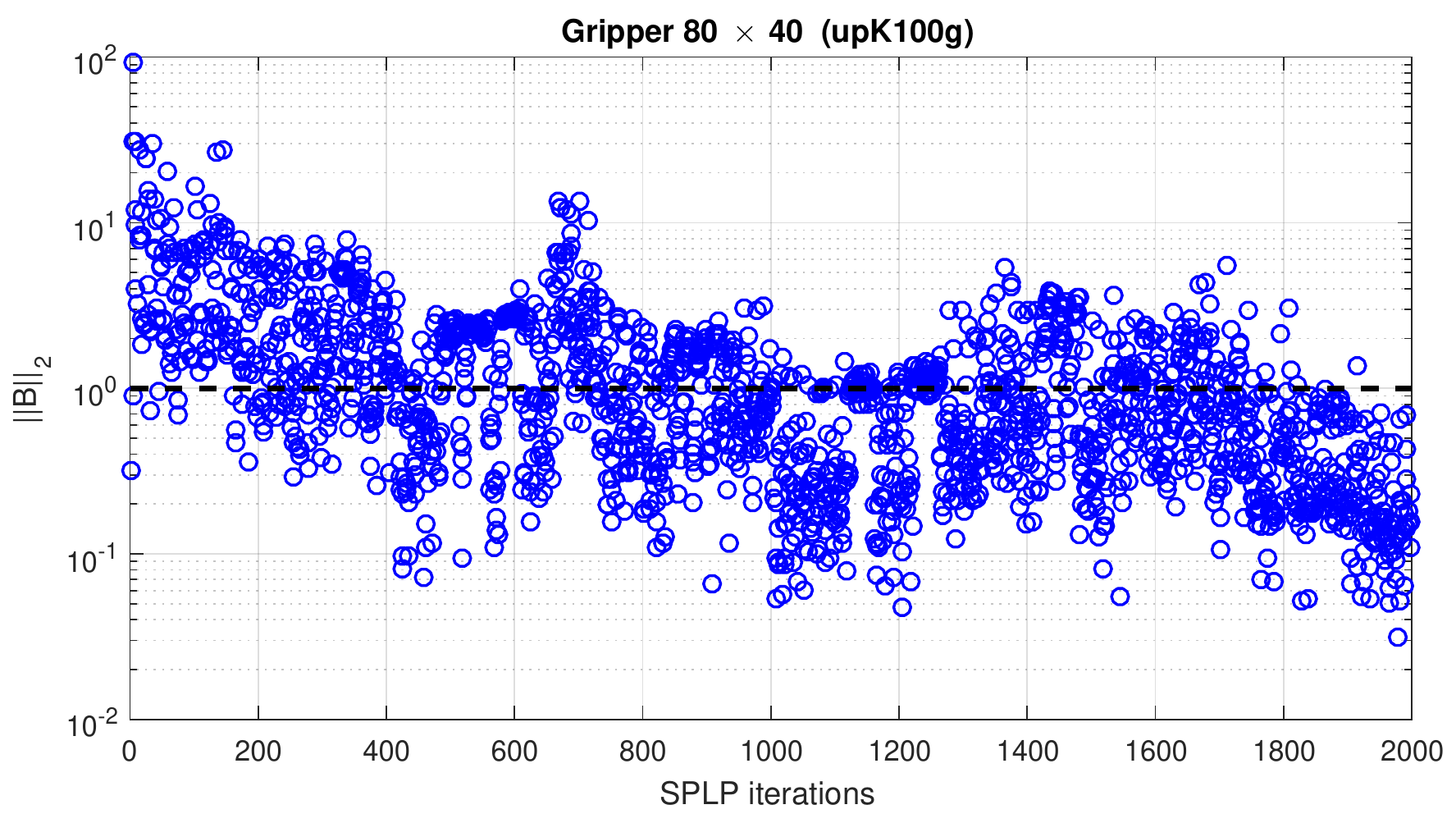}
        & \includegraphics[width=0.45\textwidth]{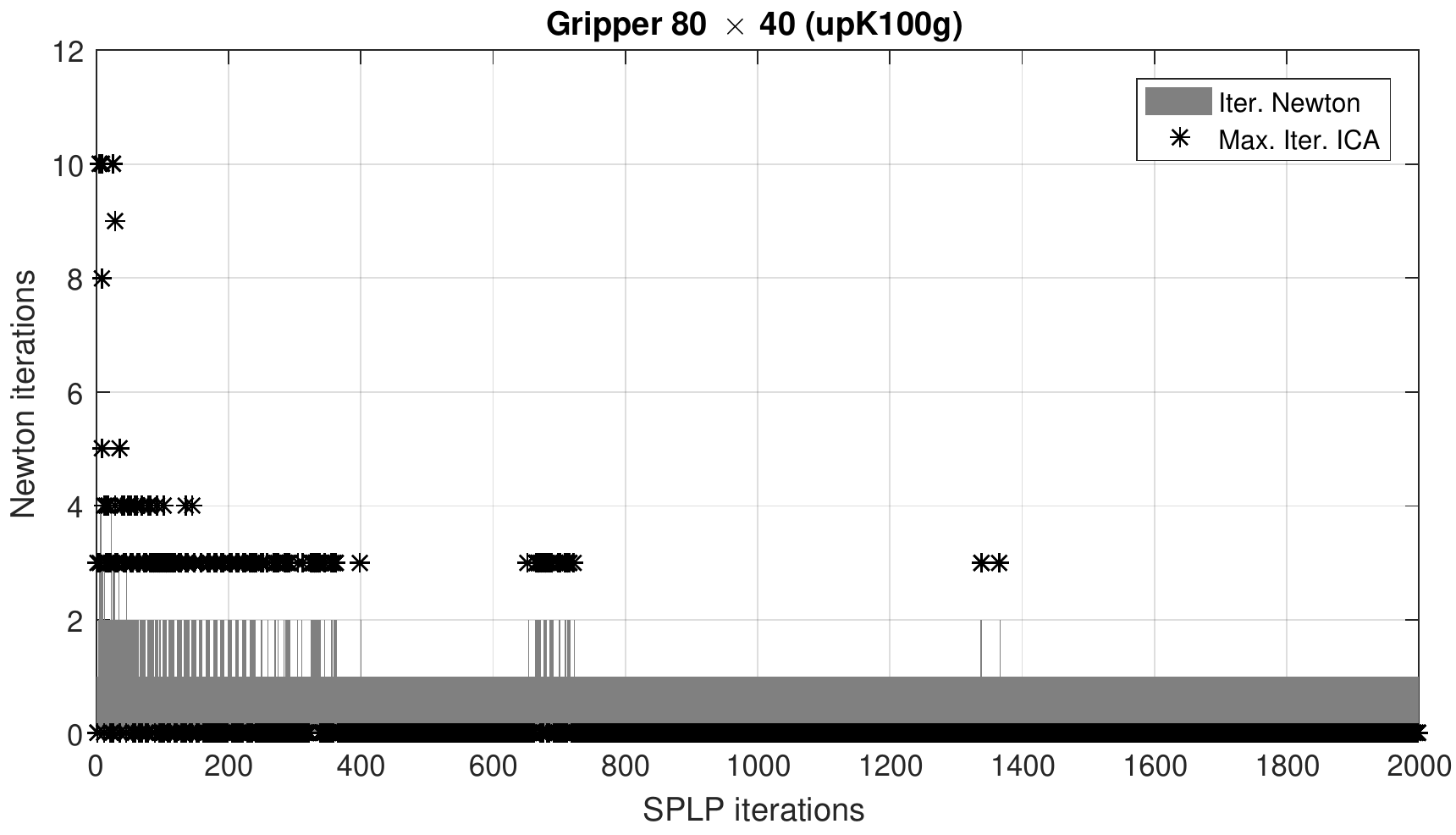} \\[0.3cm]
			    \raisebox{0.9cm}{\rotatebox[origin=lb]{90}{{\tt upK03K100g}}} \hspace{-0.2cm}
        &\includegraphics[width=0.45\textwidth]{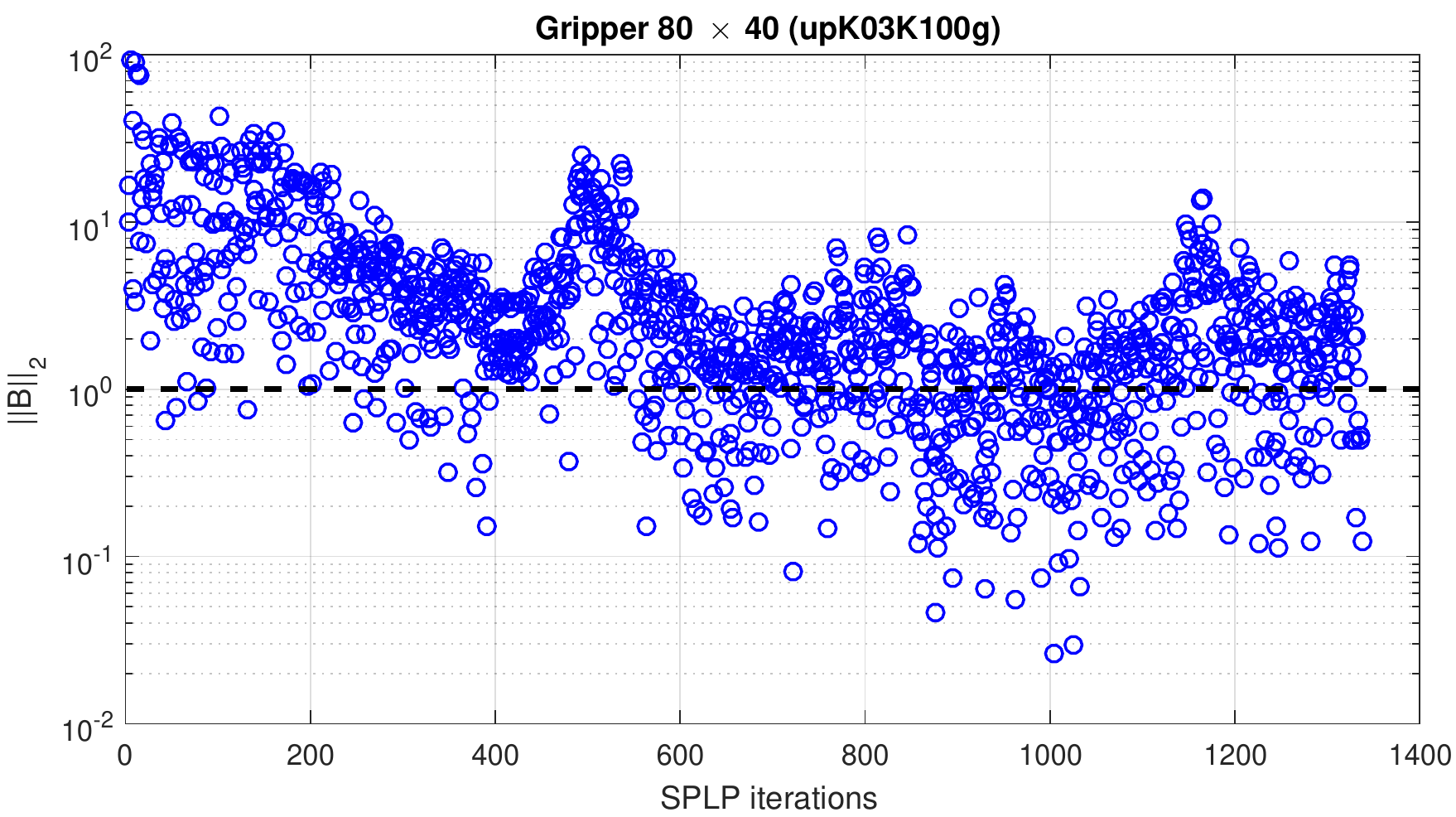} 
        &\includegraphics[width=0.45\textwidth]{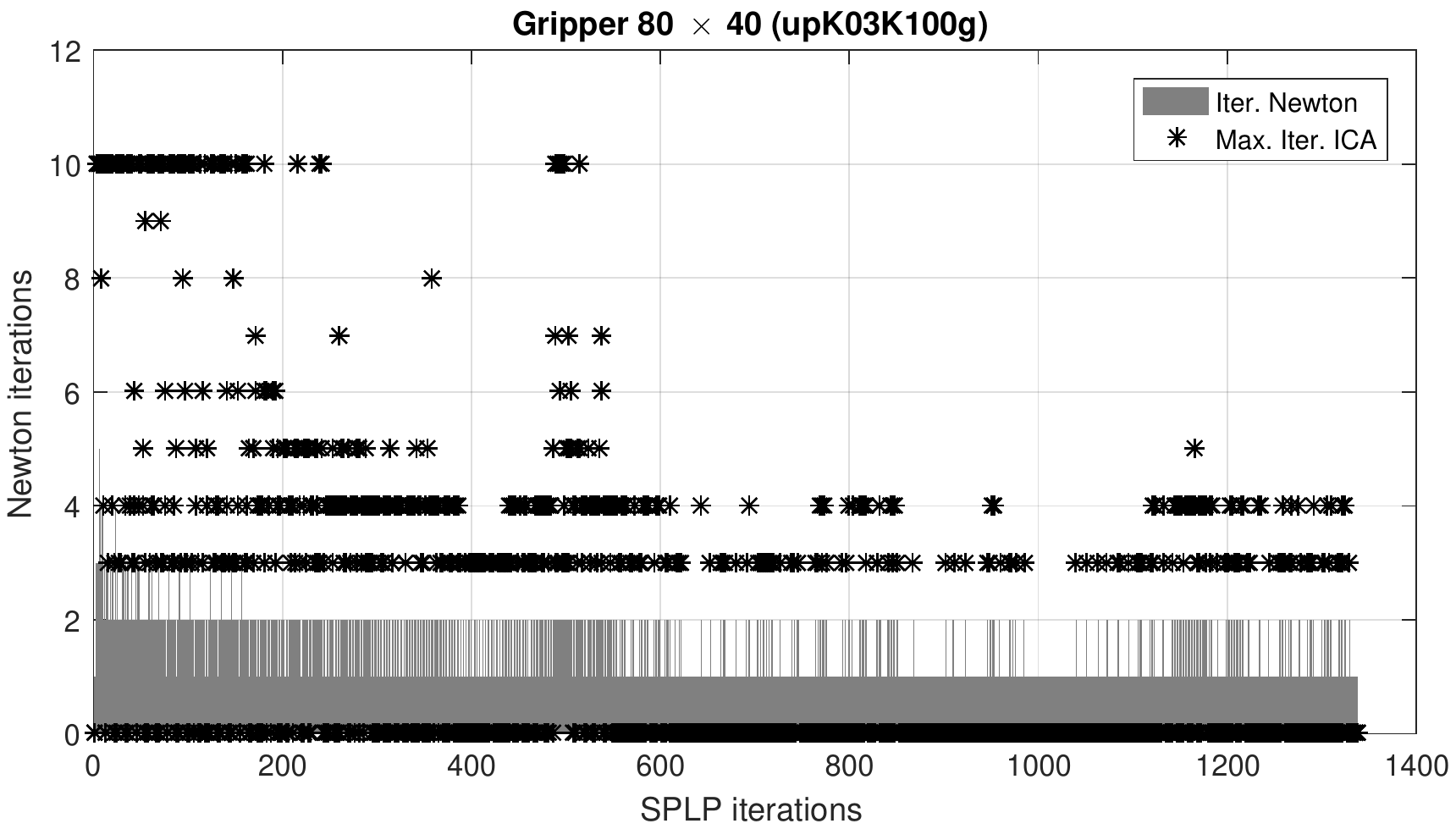}
      \end{tabular}
      \caption{Maximum value of $\| \BB \|_{2}$ (left) and number of iterations of Newton's method (right) for the gripper $80 \times 40$.} 
      \label{normb_itnewton_gripper80x40}
\end{figure}

As in the previous experiments, the optimal values of the objective function found in this new test were very similar for the two structures considered here, regardless of the applied strategy. For each combination of problem and strategy, we have recorded $\| \BB \|_{2}$ and the number of iterations of Newton's method, to see how they vary along time. The results are shown in Figures \ref{normb_itnewton_sb160x20} and \ref{normb_itnewton_gripper80x40}.

We can observe that there is a direct relation between the number of iterations of Newton's method and the value of $\| \BB \|_{2}$, i.e., the chance of performing more than 3 iterations increases with the norm of $B$. Even though, when  the {\tt upK100g} strategy is applied, the number of iterations necessary to find the approximate solution of the linear system (\ref{eq02}) is at most equals to $3$ in $93.5\,\%$ of the iterations for the slender beam and $98.2\,\%$ for the gripper. For the {\tt upK03K100g} strategy, no more than 3 linear system solutions are required in $48.6\,\%$ of the iterations for the slender beam and $73.5\,\%$ for the gripper. 

As shown in Figures \ref{normb_itnewton_sb160x20} and \ref{normb_itnewton_gripper80x40}, when the {\tt upK100g} strategy is used, the maximum value of $\| \BB \|_{2}$ is below $1$ in $41.7\,\%$ of the SPLP iterations for the slender beam, and in $56.3\,\%$ for the gripper. For the {\tt upK03K100g} strategy, this percentage is reduced to $20.3\,\%$ for the slender beam and $26.5\,\%$ for the gripper. However, for both strategies, the approximate solution of the nonlinear system (\ref{eq01b}) was always obtained after at most 4 iterations of Newton's method, no matter the problem under analysis. Therefore, we conclude that our reanalysis strategies are robust.

%===========================================
\section{Final remarks}
\label{sec:final}

Under the geometric nonlinearity assumption, the nonlinear equations that model the equilibrium of topology optimization problems must be solved repeatedly along the solution process, constituting the bottleneck of any iterative scheme. 

In this work, to accommodate the trade-off between speed and accuracy, five strategies have been proposed to steer the frequencies of the factorizations and of the updatings along the solution process. These strategies have provided significant time savings as compared with Newton's method,  without spoiling the convergence behavior of the algorithm. Indeed, strategy {\tt upK03K100g} has saved more than 60\% with the fixed budget of 300 SPLP iterations, and around 70\% when the solver is run up to convergence. The deformed structures and mechanisms obtained by means of this strategy are given in Appendix \ref{sec:final_conf}.

By choosing four different finite element meshes, we have observed that the CPU time spent for evaluating the objective function scales well with the dimensions of the global stiffness matrix.
Moreover, strategy {\tt upK03K100g} produced  the smallest rate of increase in time as the meshes were refined.

Despite the theoretical condition for convergence of the reanalysis, namely  $\| \boldsymbol{B}\|_2 < 1$, not being verified all along the iterations, we have observed that the incidental increasing of such a norm did not have a negative impact upon the results,  corroborating the robustness of the proposed strategies. 

 We have also observed that the use of the ICA scheme for solving the adjoint linear system associated with the sensitivity analysis had a crucial role in accelerating the SPLP algorithm. The combination of this scheme with the factorization of $\KK$ at every 3 iterations and the update of $\DeltaK$ at every 100 iterations of Newton's method turned {\tt upK03K100g} into a fast, yet robust, approach, outperforming the remaining strategies.

Related topics for future research include the mixing of the proposed strategies with further acceleration schemes, such as reduced-order modeling and multigrid methods, as well as handling additional stress constraints into the model (see e.g.\cite{Amir2021,DVT2021,Holmberg2013}). Another room for investigation would be to admit geometric nonlinearities within the minimum weight problem formulation  (cf. \cite{Amir2015,LongEtal2019}), to develop the inherent approximate reanalysis framework by means of iterative combined approximations, and to assess its impact.

\begin{paragraph}
{\bf Acknowledgements.} 
This work was partially supported by \textit{Funda\c{c}\~ao de Amparo \`a Pesquisa do Estado de S\~ao Paulo} 
(grants 2013/07375-0,  2018/24293-0); and \textit{Conselho Nacional de Desenvolvimento Cient\'{\i}fico 
e Tecnol\'ogico} (grant 305010/2020-4). 
\end{paragraph}

\bibliographystyle{abbrv}

\bibliography{ICA_bib}

\newpage

\appendix

\section{Results obtained for 300 iterations of SPLP}\label{sec:tables_300it}

\begin{table}[htbp!]
    \centering
    \caption{Results obtained for 300 iterations of SPLP (Slender beam $200 \times 25$).} \vspace{0.3cm}
    \label{table_100_iter_sb200x25}
        \begin{tabular}{lc@{}c@{}c@{}c@{}c@{}c@{}c@{}}
        \hline\noalign{\smallskip}
                \textbf{Slender beam} & \multirow{2}{*}{\tt N} & \multirow{2}{*}{\tt MN} & \multirow{2}{*}{\tt upK1} & \multirow{2}{*}{\tt upK1g} & \multirow{2}{*}{\tt upK100} & \multirow{2}{*}{ \tt upK100g } & { \tt upK03} \\
        $\mathbf{200 \times 25}$ &  &  &  &  &  &  & { \tt K100g} \\
        \noalign{\smallskip}\hline\noalign{\smallskip} 
        $F(\bsrho)$ - value              & \ $134.81$ \  & \ $134.81$ \  & \ $134.81$  \   & \ $134.81$ \  & \ $134.81$  &  $134.81$   &  $134.81$ \\
        \# Newton iter.             & $712$      & $1420$     & $792$        & $792$      & $793$      & $793$    & $887$    \\
       \hline
        \noalign{\smallskip}
        CPU time (s)  &  &  &  &  &  &  &  \\         
        Total                  & $77.68$    & $82.50$    & $69.11$      & $66.86$    & $62.97$    & $59.27$  & $48.25$  \\
	    $F(\bsrho)$            & $66.93$    & $70.72$    & $58.24$      & $55.85$    & $52.04$    & $48.25$  & $36.26$  \\       
        $\KK$                 & $22.49$    & $36.02$    & $25.20$      & $26.81$    & $18.71$    & $19.59$  & $13.14$  \\
        RHS                & $3.04$     & $5.39$     & $4.32$       & $5.14$     & $4.35$     & $5.18$   & $6.91$   \\
	    Factorizations        & $40.36$    & $27.38$    & $25.88$      & $19.79$    & $26.06$    & $19.34$  & $10.61$  \\
	    Linear systems         & $1.04$     & $1.93$     & $2.84$       & $4.11$     & $2.92$     & $4.14$   & $5.60$   \\
	    $\nabla F(\bsrho)$   & $1.15$     & $1.25$     & $1.17$       & $1.23$     & $1.14$     & $1.28$   & $1.31$   \\
	    SPLP solving          & $9.40$     & $10.32$    & $9.49$       & $9.57$     & $9.59$     & $9.53$   & $10.08$  \\
	    Filtering            & $0.16$     & $0.17$     & $0.17$       & $0.17$     & $0.16$     & $0.17$   & $0.20$   \\
	    Other                & $0.04$     & $0.04$     & $0.04$       & $0.04$     & $0.04$     & $0.04$   & $0.40$   \\	      
       \noalign{\smallskip}\hline
    \end{tabular}
\end{table}

\begin{table}[htbp!]
    \centering
    \caption{Results obtained for 300 iterations of SPLP (Slender Beam $400 \times 50$).} \vspace{0.3cm}
    \label{table_100_iter_sb400x50}
    \begin{tabular}{lc@{}c@{}c@{}c@{}c@{}c@{}c@{}}
        \hline\noalign{\smallskip}
    \textbf{Slender beam} & \multirow{2}{*}{\tt N} & \multirow{2}{*}{\tt MN} & \multirow{2}{*}{\tt upK1} & \multirow{2}{*}{\tt upK1g} & \multirow{2}{*}{ \tt upK100} & \multirow{2}{*}{ \tt upK100g } & { \tt upK03} \\
        $\mathbf{400 \times 50}$ &  &  &  &  &  &  & { \tt K100g} \\
        \noalign{\smallskip}\hline\noalign{\smallskip} 
        $F(\bsrho)$ - value   &  \ $136.88$  \  & \  $136.88$  \  & \ $136.88$  \  & \  $136.88$ \   &  $136.88$  &  $136.88$  &  $136.88$ \\
        \# Newton iter.             & $659$      & $975$      & $727$        & $727$      & $727$       & $727$    & $872$    \\ 
               \hline
        \noalign{\smallskip}
        CPU time (s)  &  &  &  &  &  &  &  \\ 
        Total                & $463.66$   & $395.47$   & $379.82$     & $348.30$   & $364.02$    & $324.34$ & $292.08$ \\
    $F(\bsrho)$            & $366.68$   & $298.04$   & $282.30$     & $248.50$   & $265.01$    & $226.43$ & $187.72$ \\      
        $\KK$                 & $86.36$    & $103.34$   & $89.72$      & $92.11$    & $70.75$     & $71.93$  & $53.17$  \\
        RHS                   & $12.35$    & $16.10$    & $14.07$      & $17.09$    & $14.27$     & $16.70$  & $29.80$  \\
    Factorization          & $262.93$   & $171.85$   & $170.85$     & $127.55$   & $172.28$    & $125.99$ & $74.67$  \\
    Linear system           & $5.04$     & $6.75$     & $7.66$       & $11.75$    & $7.71$      & $11.81$  & $30.08$  \\
    $\nabla F(\bsrho)$      & $10.63$    & $10.70$    & $10.59$      & $10.75$    & $10.61$     & $10.82$  & $11.40$  \\
    SPLP solving           & $79.78$    & $80.15$    & $80.33$      & $82.34$    & $81.70$     & $80.51$  & $85.89$  \\
    Filtering               & $6.07$     & $6.08$     & $6.10$       & $6.21$     & $6.20$      & $6.08$   & $6.54$   \\
    Other                & $0.50$     & $0.50$     & $0.50$       & $0.50$     & $0.50$      & $0.50$   & $0.53$   \\      
        \noalign{\smallskip}\hline
    \end{tabular}
\end{table}

\begin{table}[htbp!]
    \centering
    \caption{Results obtained for 300 iterations of SPLP (Slender beam $600 \times 75$).} \vspace{0.3cm}
    \label{table_100_iter_sb600x75}
     \begin{tabular}{lc@{}c@{}c@{}c@{}c@{}c@{}c@{}}
        \hline\noalign{\smallskip}
                        \textbf{Slender beam} & \multirow{2}{*}{\tt N} & \multirow{2}{*}{\tt MN} & \multirow{2}{*}{\tt upK1} & \multirow{2}{*}{\tt upK1g} & \multirow{2}{*}{\tt upK100} & \multirow{2}{*}{\tt upK100g} & {\tt upK03} \\
        $\mathbf{600 \times 75}$ &  &  &  &  &  &  & {\tt K100g} \\
        \noalign{\smallskip}\hline\noalign{\smallskip} 
        $F(\bsrho)$ - value              & $138.62$   & $138.52$   & $138.68$     & $138.77$   & $138.68$   & $138.77$  & $138.62$ \\
        \# Newton iter.            & $648$      & $758$      & $658$        & $667$      & $658$      & $670$     & $838$    \\
        \hline
        \noalign{\smallskip}
        CPU time (s)  &  &  &  &  &  &  &  \\ 
        Total                & $1622.93$ \ & \ $1262.76$  \ & \ $1293.72$  \ & \ $1201.68$  \ & \ $1263.94$  \ & \ $1170.56$ \ & \ $876.84$ \\
	    $F(\bsrho)$           & $1208.74$  & $853.88$   & $878.07$     & $770.80$   & $846.39$   & $741.24$  & $463.44$ \\       
        $\KK$                 & $189.38$   & $210.10$   & $191.42$     & $199.88$   & $157.05$   & $162.87$  & $106.44$ \\
        RHS                   & $26.97$    & $30.74$    & $29.07$      & $43.71$    & $29.01$    & $46.34$   & $59.35$  \\
	    Factorization         & $979.93$   & $599.03$   & $640.57$     & $480.06$   & $643.33$   & $479.84$  & $232.96$ \\
	    Linear system         & $12.46$    & $14.01$    & $17.01$      & $47.15$    & $17.00$    & $52.19$   & $64.69$  \\
	    $\nabla F(\bsrho)$    & $34.74$    & $35.33$    & $35.09$      & $35.12$    & $35.06$    & $35.37$   & $35.35$  \\
	    SPLP solving          & $351.85$   & $345.64$   & $352.93$     & $367.64$   & $354.74$   & $365.88$  & $350.45$ \\
	    Filtering              & $25.71$    & $25.97$    & $25.73$      & $26.23$    & $25.86$    & $26.16$   & $25.69$  \\
	    Other                  & $1.89$     & $1.94$     & $1.90$       & $1.89$     & $1.89$     & $1.91$    & $1.91$   \\    
        \noalign{\smallskip}\hline
    \end{tabular}
\end{table}

\begin{table}[htbp!]
    \centering
    \caption{Results obtained for 300 iterations of SPLP (Slender beam $800 \times 100$).} \vspace{0.3cm}
    \label{table_100_iter_sb800x100}
         \begin{tabular}{lc@{}c@{}c@{}c@{}c@{}c@{}c@{}}
       \hline\noalign{\smallskip}
     \textbf{Slender beam} & \multirow{2}{*}{\tt N} & \multirow{2}{*}{\tt MN} & \multirow{2}{*}{\tt upK1} & \multirow{2}{*}{\tt upK1g} & \multirow{2}{*}{\tt upK100} & \multirow{2}{*}{\tt upK100g} & {\tt upK03} \\
        $\mathbf{800 \times 100}$ &  &  &  &  &  &  & {\tt K100g} \\
        \noalign{\smallskip}\hline\noalign{\smallskip} 
        $F(\bsrho)$ - value              & $140.22$   & $140.22$   & $140.04$     & $140.23$   & $140.04$     & $140.04$  & $140.04$   \\
        \# Newton iter.             & $649$      & $754$      & $650$        & $656$      & $650$        & $582$     & $830$      \\
          \hline
        \noalign{\smallskip}
        CPU time (s)  &  &  &  &  &  &  &  \\
        Total                  & $4068.97$  \ & \ $3365.41$  \ & \ $3275.92$    \ & \ $3014.94$  \ & \ $3233.85$    \ & \ $2830.65$ \ & \ $2485.81$  \\
	    $F(\bsrho)$          & $2454.18$  & $1746.14$  & $1660.74$    & $1390.40$  & $1610.00$    & $1277.13$ & $861.18$   \\       
        $\KK$                 & $336.94$   & $336.89$   & $349.12$     & $360.41$   & $280.64$     & $254.73$  & $185.56$   \\
        RHS                    & $48.22$    & $53.72$    & $52.21$      & $77.82$    & $51.67$      & $52.11$   & $102.37$   \\
	    Factorizations        & $2045.59$  & $1329.63$  & $1238.12$    & $865.55$   & $1246.52$    & $932.81$  & $462.46$   \\
	    Linear systems        & $23.43$    & $25.90$    & $21.29$      & $86.62$    & $31.17$      & $37.48$   & $109.79$   \\
	    $\nabla F(\bsrho)$   & $261.85$   & $267.79$   & $266.58$     & $269.33$   & $266.78$     & $242.10$  & $274.78$   \\
	    SPLP solving          & $1107.07$  & $1106.43$  & $1102.45$    & $1109.93$  & $1112.00$    & $1085.03$ & $1103.35$  \\
	    Filtering             & $229.20$   & $228.41$   & $229.45$     & $228.28$   & $228.40$     & $206.20$  & $229.25$   \\
	    Other               & $16.67$    & $16.64$    & $16.70$      & $17.00$    & $16.67$      & $20.19$   & $17.25$    \\      
        \noalign{\smallskip}\hline
    \end{tabular}
\end{table} 

\begin{table}[htbp!]
    \centering
    \caption{Results obtained for 300 iterations of SPLP (Inverter $200 \times 100$).} \vspace{0.3cm}
    \label{table_100_iter_inverter200x200}
        \begin{tabular}{lc@{}c@{}c@{}c@{}c@{}c@{}c@{}}
        \hline\noalign{\smallskip}
             \textbf{Inverter} & \multirow{2}{*}{\tt N} & \multirow{2}{*}{\tt MN} & \multirow{2}{*}{\tt upK1} & \multirow{2}{*}{\tt upK1g} & \multirow{2}{*}{\tt upK100} & \multirow{2}{*}{\tt upK100g} & {\tt upK03} \\
        $\mathbf{200 \times 100}$ &  &  &  &  &  &  & {\tt K100g} \\
        \noalign{\smallskip}\hline\noalign{\smallskip} 
        $F(\bsrho)$ - value              &  $-6.9372$  \ & \ $-6.9374$    \ & \ $-6.9375$   \ & \ $-6.9330$  \ & \ $-6.9330$  \ & \ $-6.9373$  \ & \ $-6.9371$ \\
        \# Newton iter.            &  $680$      & $928$       & $850$      & $833$     & $833$     & $853$     & $923$     \\
        \hline
        \noalign{\smallskip}
        CPU time (s)  &  &  &  &  &  &  &  \\
        Total                 &  $589.85$   & $427.84$    & $432.10$   & $362.21$  & $412.60$  & $342.83$  & $257.88$  \\
	    $F(\bsrho)$           &  $536.27$   & $373.98$    & $377.97$   & $307.35$  & $358.16$  & $287.61$  & $202.73$  \\       
        $\KK$                  &  $72.27$    & $82.49$     & $79.33$    & $80.00$   & $56.82$   & $59.24$   & $48.24$   \\
        RHS                   &  $12.26$    & $15.45$     & $17.98$    & $22.45$   & $17.78$   & $23.07$   & $29.35$   \\
	    Factorizations          &  $446.15$   & $268.96$    & $269.07$   & $184.50$  & $272.20$  & $184.17$  & $93.75$   \\
	    Linear systems          &  $5.59$     & $7.08$      & $11.59$    & $20.40$   & $11.36$   & $21.13$   & $31.39$   \\
	    $\nabla F(\bsrho)$     &  $10.68$    & $10.79$     & $10.80$    & $11.08$   & $10.79$   & $11.06$   & $11.07$   \\
	    SPLP solving          &  $36.29$    & $36.40$     & $36.67$    & $37.02$   & $37.02$   & $37.47$   & $37.38$   \\
	    Filtering             &  $6.15$     & $6.20$      & $6.19$     & $6.29$    & $6.16$    & $6.22$    & $6.23$    \\
	    Other                  &  $0.46$     & $0.47$      & $0.47$     & $0.47$    & $0.47$    & $0.47$    & $0.47$    \\    
        \noalign{\smallskip}\hline
    \end{tabular}
\end{table}

\begin{table}[htbp!]
    \centering
    \caption{Results obtained for 300 iterations of SPLP (Inverter $300 \times 150$).} \vspace{0.3cm}
    \label{table_100_iter_inverter300x300}
        \begin{tabular}{lc@{}c@{}c@{}c@{}c@{}c@{}c@{}}
        \hline\noalign{\smallskip}
             \textbf{Inverter} & \multirow{2}{*}{\tt N} & \multirow{2}{*}{\tt MN} & \multirow{2}{*}{\tt upK1} & \multirow{2}{*}{\tt upK1g} & \multirow{2}{*}{\tt upK100} & \multirow{2}{*}{\tt upK100g} & {\tt upK03} \\
        $\mathbf{300 \times 150}$ &  &  &  &  &  &  & {\tt K100g} \\
        \noalign{\smallskip}\hline\noalign{\smallskip} 
        $F(\bsrho)$ - value              & $-8.6969$   \ & \ $-8.6989$   \ & \ $-8.7020$  \ & \ $-8.7035$   \ & \ $-8.6997$  \ & \ $-8.7028$  \ & \ $-8.7046$ \\
        \# Newton iter.            & $696$       & $1065$      & $859$     & $849$      & $848$     & $841$     & $929$     \\
        \hline
        \noalign{\smallskip}
        CPU time (s)  &  &  &  &  &  &  &  \\
        Total                 & $1883.89$   & $1334.16$   & $1344.36$ & $1097.37$  & $1293.52$ & $1066.04$ & $756.60$  \\
	    $F(\bsrho)$.           & $1691.98$   & $1144.86$   & $1149.91$ & $905.20$   & $1102.38$ & $860.02$  & $564.05$  \\       
        $\KK$                  & $161.20$    & $194.80$    & $178.60$  & $179.35$   & $128.35$  & $132.06$  & $106.75$  \\
        RHS                   & $28.87$     & $38.90$     & $42.83$   & $54.16$    & $42.22$   & $54.37$   & $66.69$   \\
	  Factorizations         & $1487.78$   & $891.76$    & $895.51$  & $614.71$   & $899.80$  & $616.35$  & $310.82$  \\
	    Linear systems          & $14.13$     & $19.40$     & $32.97$   & $56.98$    & $31.98$   & $57.24$   & $79.79$   \\
	    $\nabla F(\bsrho)$     & $35.44$     & $36.17$     & $37.69$   & $36.37$    & $36.27$   & $36.61$   & $36.57$   \\
	    SPLP solving          & $128.50$    & $124.56$    & $126.68$  & $126.87$   & $126.22$  & $140.25$  & $127.22$  \\
	    Filtering             & $26.23$     & $26.77$     & $28.17$   & $27.09$    & $26.88$   & $27.29$   & $26.90$   \\
	    Other                 & $1.74$      & $1.80$      & $1.91$    & $1.84$     & $1.77$    & $1.87$    & $1.86$    \\	      
        \noalign{\smallskip}\hline
    \end{tabular}
\end{table}

\begin{table}[htbp!]
    \centering
    \caption{Results obtained for 300 iterations of SPLP (Inverter $400 \times 200$).} \vspace{0.3cm}
    \label{table_100_iter_inverter400x400}
        \begin{tabular}{lc@{}c@{}c@{}c@{}c@{}c@{}c@{}}
        \hline\noalign{\smallskip}
             \textbf{Inverter} & \multirow{2}{*}{\tt N} & \multirow{2}{*}{\tt MN} & \multirow{2}{*}{\tt upK1} & \multirow{2}{*}{\tt upK1g} & \multirow{2}{*}{\tt upK100} & \multirow{2}{*}{\tt upK100g} & {\tt upK03} \\
        $\mathbf{400 \times 200}$ &  &  &  &  &  &  & {\tt K100g} \\
        \noalign{\smallskip}\hline\noalign{\smallskip} 
        $F(\bsrho)$ - value              & $-8.4352$   \ & \ $-8.4342$  \ & \ $-8.4356$ \ & \ $-8.4352$   \ & \ $-8.4355$ \ & \ $-8.4399$ \ & \ $-8.4362$  \\
        \# Newton iter.             & $700$       & $930$       & $797$     & $783$       & $784$     & $804$     & $903$      \\
        \hline
        \noalign{\smallskip}
        CPU time (s)  &  &  &  &  &  &  &  \\
        Total                 & $4827.72$ & $3482.96$  &  $3501.80$ & $2845.28$  &  $3408.99$  &  $2823.02$  &  $2123.78$  \\
	    $F(\bsrho)$           & $3924.52$   & $2575.07$   & $2580.67$ & $1959.74$   & $2502.04$ & $1893.70$ & $1191.78$  \\       
        $\KK$                 & $286.57$    & $324.00$    & $302.85$  & $308.58$    & $227.14$  & $236.92$  & $190.77$   \\
        RHS                    & $51.31$     & $62.80$     & $68.74$   & $86.44$     & $67.57$   & $87.90$   & $111.88$   \\
	  Factorization          & $3559.94$   & $2155.70$   & $2158.38$ & $1476.52$   & $2157.75$ & $1478.63$ & $753.57$   \\
	   Linear system         & $26.70$     & $32.57$     & $50.70$   & $88.20$     & $49.58$   & $90.25$   & $135.56$   \\
	   $\nabla F(\bsrho)$     & $268.66$    & $273.07$    & $272.17$  & $264.29$    & $273.24$  & $283.94$  & $286.48$   \\
	   SPLP solving           & $385.21$    & $383.65$    & $398.97$  & $378.80$    & $383.87$  & $388.66$  & $390.09$   \\
	    Filtering             & $233.72$    & $235.29$    & $234.31$  & $226.75$    & $234.06$  & $239.92$  & $238.60$   \\
	    Other                 & $15.61$     & $15.88$     & $15.68$   & $15.70$     & $15.78$   & $17.40$   & $16.83$    \\     
        \noalign{\smallskip}\hline
    \end{tabular}
\end{table}

\begin{table}[htbp!]
    \centering
    \caption{Results obtained for 300 iterations of SPLP (Inverter $500 \times 250$).} \vspace{0.3cm}
    \label{table_100_iter_inverter500x500}
        \begin{tabular}{lc@{}c@{}c@{}c@{}c@{}c@{}c@{}}
        \hline\noalign{\smallskip}
         \textbf{Inverter} & \multirow{2}{*}{\tt N} & \multirow{2}{*}{\tt MN} & \multirow{2}{*}{\tt upK1} & \multirow{2}{*}{\tt upK1g} & \multirow{2}{*}{\tt upK100} & \multirow{2}{*}{\tt upK100g} & {\tt upK03} \\
        $\mathbf{500 \times 250}$ &  &  &  &  &  &  & {\tt K100g} \\
        \noalign{\smallskip}\hline\noalign{\smallskip} 
        $F(\bsrho)$ - value              & $-7.7597$    \ & \ $-7.7599$    \ & \ $-7.7600$  \ & \ $-7.7598$    \ & \ $-7.7597$  \ & \ $-7.7595$  \ & \ $-7.7596$   \\
        \# Newton iter.            & $673$       & $823$       & $741$     & $745$       & $750$     & $747$     & $864$       \\   
        \hline
        \noalign{\smallskip}
        CPU time (s)  &  &  &  &  &  &  &  \\
        Total                 & $10214.38$   &  $7165.38$  & $7186.63$  &  $5867.43$   &  $7060.42$  &  $5703.46$  &  $4334.30$   \\
	    $F(\bsrho)$           & $8306.03$   & $5285.67$   & $5295.14$ & $3944.25$   & $5157.13$ & $3818.24$ & $2394.37$   \\       
        $\KK$                & $549.73$    & $594.25$    & $570.72$  & $582.55$    & $440.70$  & $450.47$  & $298.57$    \\
        RHS                    & $90.72$     & $91.38$     & $101.18$  & $128.63$    & $102.14$  & $127.72$  & $169.83$    \\
	    Factorizations        & $7618.48$   & $4547.31$   & $4549.31$ & $3102.94$   & $4538.28$ & $3111.46$ & $1694.96$   \\
	    Linear systems          & $47.10$     & $52.73$     & $73.72$   & $130.13$    & $76.01$   & $128.59$  & $231.01$    \\
	    $\nabla F(\bsrho)$   & $669.89$    & $662.19$    & $671.27$  & $690.44$    & $672.98$  & $659.28$  & $717.43$    \\
	    SPLP solving          & $597.88$    & $577.06$    & $576.98$  & $588.24$    & $584.80$  & $579.47$  & $578.53$    \\
	    Filtering           & $601.68$    & $601.35$    & $604.05$  & $603.05$    & $605.36$  & $604.38$  & $599.50$    \\
	    Other                 & $38.90$     & $39.11$     & $39.09$   & $41.45$     & $39.15$   & $42.09$   & $44.47$     \\      
        \noalign{\smallskip}\hline
    \end{tabular}
\end{table}

\newpage

\section{Final configurations}\label{sec:final_conf}

\begin{figure}[htbp!]
   \centering 
      \begin{tabular}{cc} 
          \includegraphics[width=0.4\textwidth]{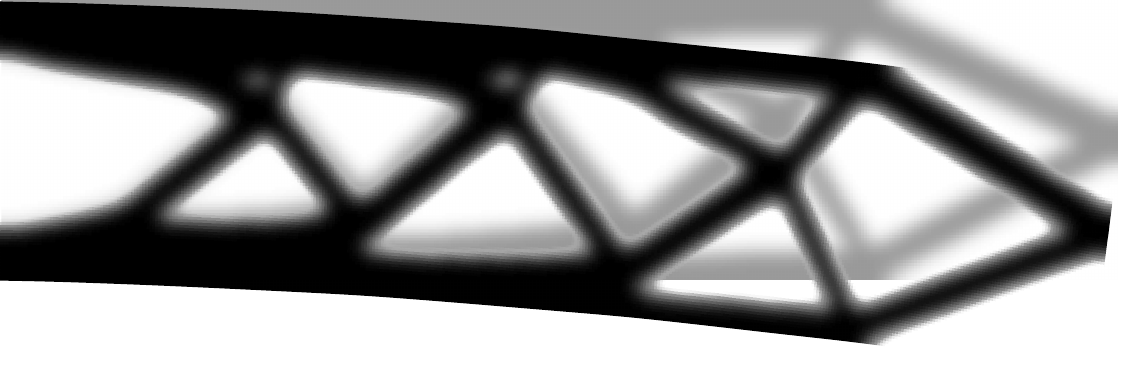} &
           \raisebox{0.3cm}{\includegraphics[width=0.5\textwidth]{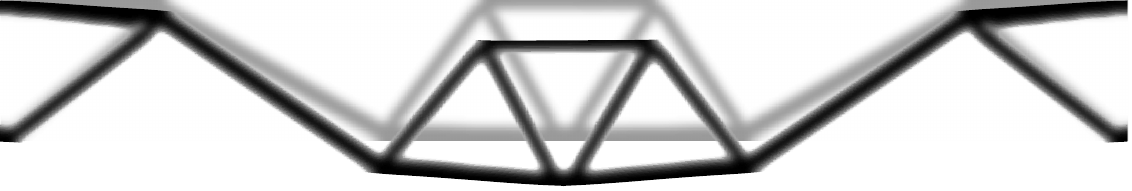}} \\
          cantilever beam ($400 \times 100$) & slender beam ($800 \times 100$)\\[0.3cm]
           \raisebox{0.1cm}{\includegraphics[width=0.4\textwidth]{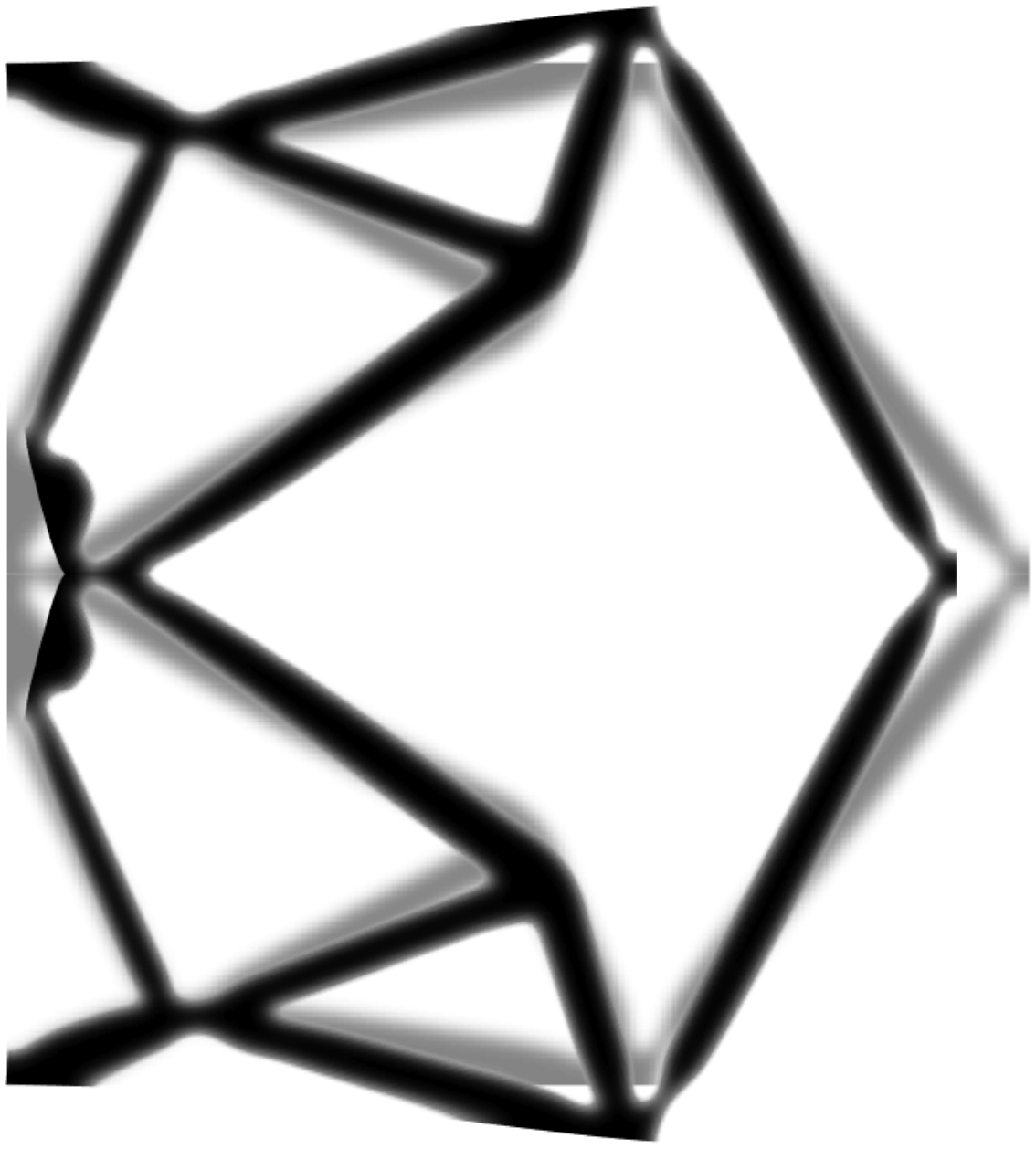}} &          
           \includegraphics[width=0.45\textwidth]{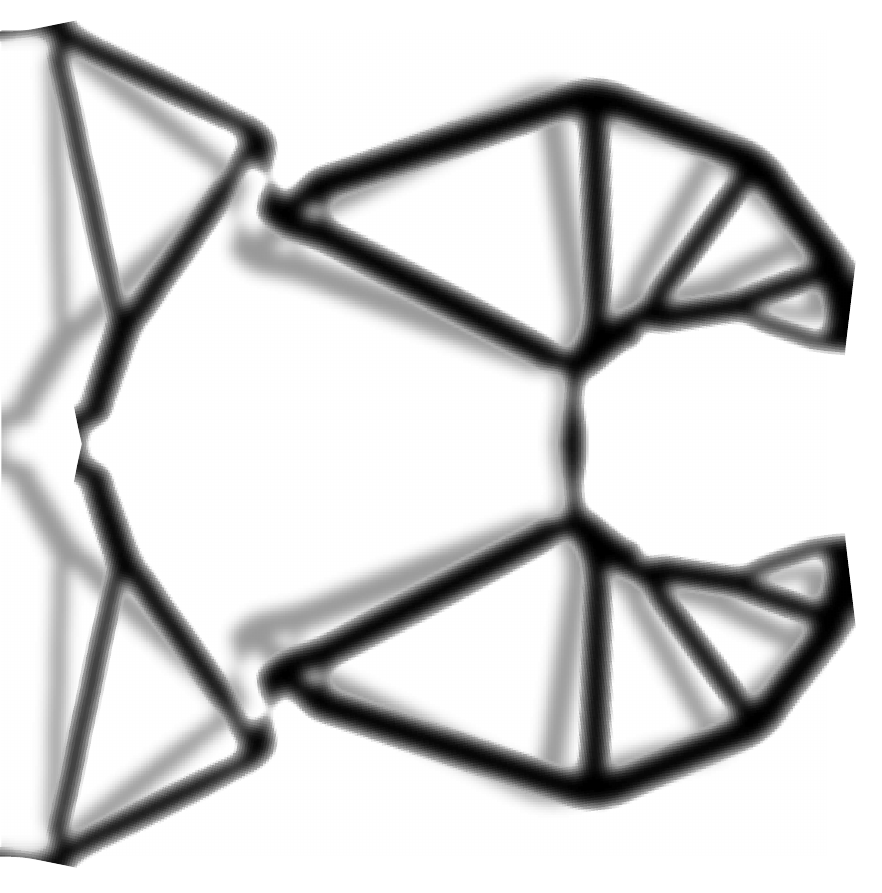} \\
           inverter ($500 \times 250$) & gripper ($320 \times 160$) \\
      \end{tabular}
      \caption{Final configurations of the instances with large displacements (reached by the {\tt upK03K100g}) strategy). The darker structures and mechanisms highlight the attained deformations with respect to the lighter gray ones.} 
      \label{fig:deformations}
\end{figure}

\end{document}